\input amstex
\documentstyle{amsppt}

\def\mod{\text{\rm mod}\, }

\def\i{\text{\rm i}}
\def\ii{\text{\rm ii}}
\def\iii{\text{\rm iii}}
\def\iv{\text{\rm iv}}
\def\v{\text{\rm v}}
\def\vi{\text{\rm vi}}

\def\bN{\Bbb N}

\font\boldtitlefont=cmb10 scaled\magstep1\font\boldfont=cmb10
\font\bigmath=cmbxti10 scaled\magstep2 \font\smmath=cmr10
 
\font\scten=cmti10
\topmatter\magnification=\magstep1
\hsize14truecm
\vsize21.6truecm
\leftheadtext {Hu and Ye}
\rightheadtext {Canonical Basis for Type ${\bigmath A}_{\bigmath
4}\;\;({\hbox {\rm II}})$}
\topmatter
\title
{\boldtitlefont  CANONICAL BASIS FOR TYPE ${\hbox{\bigmath
A}}_{\hbox{\boldfont 4}}\; ({\hbox{\boldtitlefont II}})$\\
$-$ Polynomial  elements in \\ one variable $^\dag$}
\endtitle
\thanks
$^\dag$ Supported in part by the National Natural Science
Foundation of China (10271088)\newline .\quad $^{*}$ {Supported by
the Natural Science Foundation of Henan Province (0311010100).}
\endthanks
\author
{Hu Yuwang$^{1}$~ and ~Ye Jiachen$^{2}$}
\endauthor
\affil
{\scten $^1$Department of Mathematics, Xinyang Teachers College,\\
Henan 464000, People's Republic of China\\
E-Mail: hsh0412\@sina.com\\
\smallskip
$^2$Department of Applied Mathematics, Tongji University,\\
Shanghai 200092, People's Republic of China\\
E-Mail: jcye\@mail.tongji.edu.cn\\}
\endaffil
\endtopmatter
\noindent {\bf Abstract:}\quad {\smmath All the $62$ monomial
elements in the canonical basis $\bold B$ of the quantized
enveloping algebra for type $A_4$ have been determined in [2].
According to Lusztig's idea [7], the elements in the canonical
basis $\bold B$  consist of monomials and linear combinations of
monomials (for convenience, we call them polynomials). In this
note, we compute all the $144$ polynomial elements in one variable
in the canonical basis $\bold B$ of the quantized enveloping
algebra for type $A_4$ based on our joint note [2]. We conjecture
that there are other polynomial elements in two or three variables
in the canonical basis $\bold B$, which include independent
variables and dependent variables. Moreover, it is conjectured
that there are no polynomial elements in the canonical basis
$\bold B$ with four or more variables.}
\medskip
\noindent {\bf 2000 MR Subject Classification:} 17B37\quad 20G42\quad 81R50
\vskip1cm
\document

\head {1. Polynomial elements in one variable in the canonical
basis of $U^+$}\endhead

We shall freely use the notations in [2] without further comments.

\par \vskip0.1cm

{\bf 1.1.} In [2], we have determined the $62$ monomial elements
in the canonical basis $\bold B$ of the quantized enveloping
algebra for type $A_4$. Each monomial element corresponds to a
region which consists of six independent inequalities. Here
independence of six inequalities implies that we can't deduce one
inequality from the others. Moreover, the regions can't be
described with less than six inequalities, and the interiors of
any two regions are disjoint. These regions of monomial elements
have ``nice" forms, which will help us compute polynomial elements
in the canonical basis $\bold B$.

All the $144$ polynomial elements in one variable in the canonical
basis $\bold B$ will be determined in this note based on our joint
work [2]. We have been keeping on the computations of polynomial
elements in several variables in the canonical basis $\bold B$.
And we do have found more than thirty polynomial elements in two
independent variables $u$ and $w$ in the canonical basis $\bold
B$. We conjecture that there are other polynomial elements in two
or three variables in the canonical basis $\bold B$, which include
independent variables and dependent variables. Moreover, it is
conjectured that there are no polynomial elements in the canonical
basis $\bold B$ with four or more variables.

It should be mentioned here that so-called \lq\lq independent
variables\rq\rq ~is that the summing in the polynomial is
independent of the order of variables, and \lq\lq dependent
variables\rq\rq ~implies that the summing in the polynomial is
dependent of the order of variables.

\par\vskip0.1cm

{\bf 1.2.} Let us compute polynomial elements in one variable in
the canonical basis $\bold B$ of the quantized enveloping algebra
for type $A_4$. According to Lusztig's idea, these polynomial
elements should be linear combinations of the monomial elements.
We shall describe how one can compute them and begin with monomial
element {\bf 1.}(1) in [2] as an example. Let
$A=(a_{1},a_{2},a_{3},a_{4},a_{5},a_{6},a_{7},a_{8},a_{9},a_{10})\in
{\Bbb N}^{10}  $, we have known from [2] P. 242 that in order to
make the linear form $l(x_1, x_2, \cdots , x_{14})$ and the unit
form $q(x_1, x_2, \cdots , x_{14})$ to be non-negative for any
$x=(x_1, x_2, \cdots , x_{14})\in \Bbb{N}^{14}$, the following
inequalities
$$\align
& \quad a_5+a_6+a_7\geq a_1+a_2+a_3,\qquad a_7+a_9\geq a_3+a_6,
\cr & \quad a_7+a_9+a_{10}\geq a_3+a_6+a_8, \qquad a_6+a_7\geq
a_2+a_3,\cr (\sharp)&\quad
 a_7+a_8+a_9\geq a_2+a_3+a_5 ,\qquad a_2+a_5\geq
a_6+a_8,\cr& \quad a_8+a_9\geq a_5+a_6,\qquad a_9+a_{10}\geq
a_6+a_8,\cr & \quad a_1\geq a_5,\quad a_2\geq a_6,\quad a_5\geq
a_8, \quad a_{10}\geq a_8, \quad a_9\geq a_6,
\endalign $$
equivalently,
$$\align
& (*) \quad\;\; a_5\geq a_8,\qquad \;\;a_{10}\geq a_8,\qquad
\;\;a_2\geq a_6,\qquad \;\;a_1\geq a_5,\cr &\quad \quad
\;\;a_8+a_9\geq a_5+a_6, \quad \quad \;\;a_5+a_6+a_7\geq
a_1+a_2+a_3
\endalign $$
must hold, and then the monomial
$$\align
&(**)\qquad
e_{2}^{(a_3)}e_{3}^{(a_2+a_3)}e_{4}^{(a_1+a_2+a_3)}e_{2}^{(a_6)}
e_{3}^{(a_5+a_6)}e_{2}^{(a_8)}e_{1}^{(a_4+a_7+a_9+a_{10})}\times\\
&\qquad \qquad \times e_{2}^{(a_4+a_7+a_9)}
e_{3}^{(a_4+a_7)}e_{4}^{(a_4)}
\endalign $$
belongs to $\bold B$. Because it is impossible that we deduce any
one from the other inequalities in $(*)$, the six inequalities in
$(*)$ are independent. $(*)$ is called the region of the monomial
$(**)$, and the six independent inequalities in $(*)$ are called
the defining inequalities of the region $(*)$.  \par \vskip0.1cm

We now consider the linear combination of the monomial $(**)$ that
could become a new member in $\bold B$. Because there is an
one-to-one correspondence between ${\Bbb N}^{10}$ and $\bold B$,
we have to observe the regions similar to  $(*)$. Firstly, we
reverse the first defining inequality $a_5 \geq a_8$ in $(*)$, and
get a new inequality $a_8 \geq a_5$. In order to make the other
inequalities hold in $(\sharp)$ but $a_5 \geq a_8$, we have to
replace the remaining five defining inequalities in $(*)$ by the
following five defining inequalities

$$\align
& \qquad \;\;a_{10}\geq a_8,\qquad \;\;a_2+a_5\geq a_6+a_8,\qquad
\;\;a_1\geq a_5,\cr &\quad \quad \;\;a_9\geq a_6, \quad \quad
\;\;a_5+a_6+a_7\geq a_1+a_2+a_3.
\endalign $$
The five new defining inequalities together with $ a_8\geq a_5$
form a new region, which corresponds to a linear combination of
the monomial $(**)$. Similarly, we can deal with other two cases
that the second defining inequality and the third defining
inequality in $(*)$ are reversed, respectively.  Secondly, we
reverse the fourth defining inequality $a_1 \geq a_5$ in $(*)$,
and get a new inequality $a_5 \geq a_1$. In order to make the
other inequalities hold in $(\sharp)$ but $a_1 \geq a_5$, we have
to replace the remaining five defining inequalities in $(*)$ by
the following five defining inequalities
$$\align
& \qquad \;\;a_6+a_7\geq a_2+a_3,\qquad \;\;a_8+a_9\geq
a_5+a_6,\cr &\qquad  \quad a_2\geq a_6,\quad \quad \;\;a_1\geq
a_8, \quad \quad \;\;a_{10}\geq a_8.
\endalign $$
The five new defining inequalities together with $ a_5\geq a_1$
form a new region, which corresponds to a linear combination of
the monomial $(**)$. Similarly, we can deal with other two cases
that the fifth defining inequality and the sixth defining
inequality in $(*)$ are reversed, respectively.  By many
computations, we have noticed that the coefficient of the linear
combination of the monomial $(**)$ is closely related to the
reversed inequality. Therefore, corresponding to the six defining
inequalities in $(*)$, we get the following six possibilities,
respectively.
$$\align
& \qquad(\i)\;\; a_8\geq a_5,\qquad (\ii)\;\;a_8\geq a_{10},\qquad
(\iii)\;\;a_6\geq a_2,\qquad (\iv)\;\;a_5\geq a_1,\cr &\qquad
(\v)\;\;a_5+a_6\geq a_8+a_9, \quad \quad (\vi)\;\;a_1+a_2+a_3\geq
a_5+a_6+a_7,
\endalign $$
Then we get six polynomial elements in one variable, each
corresponds to one of the above six cases. Also, the region of
each of the six polynomial elements in one variable consists of
six independent inequalities. In this way, corresponding to each
of the 62 monomial elements, we compute all polynomial elements in
one variable and their regions.

Observing the regions of every six polynomial elements in one
variable which correspond to one of the 62 monomials, we find that
at least three regions may have already occurred in $\Cal S$, the
set of the $62$ regions determined by the $62$ monomial elements,
and there may be at most another two regions, which are the same
as those of polynomial elements in one variable we already
computed before. They will not make contribution to the canonical
basis $\bold B$.

In the above example, the regions corresponding to $(\i), (\ii),
(\iii)$ have already occurred in $\Cal S$, so the corresponding
three polynomial elements in one variable don't make contribution
to the canonical basis $\bold B$, although they are combinations
of the monomials $(**)$. In the other hand, the regions
corresponding to $(\iv),(\v),(\vi)$ don't occur in $\Cal S$, and
the corresponding three polynomial elements in one variable become
new members in the canonical basis $\bold B$ because this is the
first consideration for non-monomial case.

Repeating the above procedure for all $62$ monomials one by one,
we obtain all the $144$ polynomial elements in one variable, which
belong to the canonical basis $\bold B$. The main result
concerning with the polynomial elements in one variable in the
canonical basis $\bold B$ of the quantized enveloping algebra for
type $A_4$ is the following theorem.

\proclaim {Theorem 1.3} Let
$A=(a_1,a_2,a_3,a_4,a_5,a_6,a_7,a_8,a_9,a_{10})\in {\bN}^{10}$.
Then \par {\rm 1.}\; corresponding to $62$ equivalence classes for
$\sim $, the $72$ polynomial elements in one variable in the
canonical basis $\bold B$
$$\theta(A)=\theta(a_1,a_2,a_3,a_4,a_5,a_6,a_7,a_8,a_9,a_{10})\in
{\bold B}$$ are given by the following

$$\align {\bold 1.}\quad
&(1).\;\;\sum_{0\leq u\leq
a_4+a_7}\!\!\!(-1)^u\left[\!\smallmatrix a_5+a_6-a_8-a_9-1+u\cr u
\endsmallmatrix
\!\right]e_2^{(a_3)}e_3^{(a_2+a_3)}e_4^{(a_1+a_2+a_3)}e_2^{(a_6)}\cr
&\qquad\qquad\times
e_3^{(a_5+a_6+u)}e_2^{(a_8)}e_1^{(a_4+a_7+a_9+a_{10})}e_2^{(a_4+a_7+a_9)}
e_3^{(a_4+a_7-u)}e_4^{(a_4)}\cr &\qquad\qquad\text{\rm if}\qquad
a_5+a_6+a_7\geq a_1+a_2+a_3,\quad a_5+a_6\geq a_8+a_9,\cr
&\qquad\qquad\qquad \qquad  a_9\geq a_6,\quad a_2\geq a_6,\quad
a_1\geq a_5,\quad a_{10}\geq a_8.\cr&\cr &(2).\;\;\sum_{0\leq
u\leq a_4}\!\!(-1)^u \!\left[\!\smallmatrix
a_1+a_2+a_3-a_5-a_6-a_7-1+u\cr u
\endsmallmatrix
\!\right]e_2^{(a_3)}e_3^{(a_2+a_3)}e_4^{(a_1+a_2+a_3+u)}\cr &
\qquad\qquad\times
e_2^{(a_6)}e_3^{(a_5+a_6)}e_2^{(a_8)}e_1^{(a_4+a_7+a_9+a_{10})}e_2^{(a_4+a_7+a_9)}
e_3^{(a_4+a_7)}e_4^{(a_4-u)}\cr &\quad\quad\text{\rm if}\quad
a_1+a_2+a_3\geq a_5+a_6+a_7,\quad a_8+a_9\geq a_5+a_6,\quad
a_2\geq a_6,\cr & \qquad\qquad\qquad\qquad a_6+a_7\geq
a_2+a_3,\quad  a_5\geq a_8,\quad a_{10}\geq a_8.\cr
&(3).\;\;\sum_{0\leq u\leq a_2+a_3}\!\!(-1)^u\left[\smallmatrix
a_5-a_1-1+u\cr u
\endsmallmatrix
\right]e_2^{(a_3)}e_3^{(a_2+a_3-u)}e_4^{(a_1+a_2+a_3)}e_2^{(a_6)}\cr
&\qquad\qquad\times
e_3^{(a_5+a_6+u)}e_2^{(a_8)}e_1^{(a_4+a_7+a_9+a_{10})}e_2^{(a_4+a_7+a_9)}
e_3^{(a_4+a_7)}e_4^{(a_4)}\cr &\qquad\qquad\text{\rm if}\quad
a_6+a_7\geq a_2+a_3,\quad a_8+a_9\geq a_5+a_6,\quad a_2\geq
a_6,\cr &\qquad\qquad \qquad\qquad\qquad a_5\geq a_1\geq a_8,\quad
a_{10}\geq a_8. \cr&\cr {\bold 2.}\quad &(1).\;\;\sum_{0\leq u\leq
a_4+a_7+a_9}\!\!\!(-1)^u\left[\smallmatrix a_{10}-a_8-1+u\cr u
\endsmallmatrix
\right]e_3^{(a_2)}e_2^{(a_3+a_6)}e_1^{(a_4+a_7+a_9-u)}e_3^{(a_3)}\cr
&\qquad\qquad\quad\times
e_2^{(a_4+a_7)}e_3^{(a_4)}e_4^{(a_1+a_2+a_3+a_4)}e_3^{(a_5+a_6+a_7)}e_2^{(a_8+a_9)}
e_1^{(a_{10}+u)}\cr &\qquad\qquad\text{\rm if}\qquad
a_1+a_2+a_3\geq a_5+a_6+a_7,\quad a_5+a_6\geq a_8+a_9,\cr
&\quad\qquad \qquad\qquad\qquad a_9\geq a_6\geq a_2,\quad
a_{10}\geq a_8,\quad a_7\geq a_3.\cr &\cr&(2).\;\;\sum_{0\leq
u\leq a_4+a_7}\!\!\!(-1)^u\left[\smallmatrix
a_8+a_9-a_5-a_6-1+u\cr u
\endsmallmatrix
\right]e_3^{(a_2)}e_2^{(a_3+a_6)}e_1^{(a_4+a_7+a_9)}e_3^{(a_3)}\cr
&\quad\quad\quad\times
e_2^{(a_4+a_7-u)}e_3^{(a_4)}e_4^{(a_1+a_2+a_3+a_4)}e_3^{(a_5+a_6+a_7)}e_2^{(a_8+a_9+u)}
e_1^{(a_{10})}\cr &\qquad\qquad\text{\rm if}\qquad a_1+a_2+a_3\geq
a_5+a_6+a_7,\quad a_8+a_9\geq a_5+a_6,\cr &\quad\qquad
\qquad\qquad\qquad a_5\geq a_8\geq a_{10},\quad a_6\geq a_2,\quad
a_7\geq a_3.\cr&\cr &(3).\;\;\sum_{0\leq u\leq
a_4+a_7}\!\!\!(-1)^u\left[\smallmatrix a_6-a_9-1+u\cr u
\endsmallmatrix
\right]e_3^{(a_2)}e_2^{(a_3+a_6+u)}e_1^{(a_4+a_7+a_9)}e_3^{(a_3)}\cr
&\qquad\qquad\times
e_2^{(a_4+a_7-u)}e_3^{(a_4)}e_4^{(a_1+a_2+a_3+a_4)}e_3^{(a_5+a_6+a_7)}e_2^{(a_8+a_9)}
e_1^{(a_{10})}\cr &\qquad\qquad\text{\rm if}\qquad a_1+a_2+a_3\geq
a_5+a_6+a_7,\quad a_5\geq a_8\geq a_{10},\cr &\qquad\qquad\qquad
\qquad\quad a_6\geq a_9\geq a_2,\quad a_7\geq a_3.\cr &\cr
\endalign$$$$\align {\bold 3.}\quad &(1).\;\;\sum_{0\leq u\leq
a_3+a_4+a_6}\!\!\!\!\!\!\!\!\!\!(-1)^u\left[\smallmatrix
a_8-a_5-1+u\cr u
\endsmallmatrix
\right]e_1^{(a_4)}e_3^{(a_2)}e_2^{(a_3+a_4+a_6-u)}e_4^{(a_1+a_2)}\cr
&\qquad\qquad\times
e_1^{(a_7+a_9)}e_3^{(a_3+a_4+a_5+a_6)}e_2^{(a_7+a_8+a_9+u)}e_4^{(a_3+a_4)}e_1^{(a_{10})}
e_3^{(a_7)}\cr &\qquad\qquad\text{\rm if}\qquad a_5+a_6\geq
a_1+a_2,\quad a_3+a_6\geq a_7+a_9,\quad a_1\geq a_5,\cr
&\qquad\qquad \qquad\qquad\qquad\qquad  a_8\geq a_5\geq
a_{10},\quad a_9 \geq a_6.\cr&\cr
 &(2).\;\;\sum_{0\leq u
\leq a_7}\!(-1)^u\left[\smallmatrix a_5+a_6-a_8-a_9-1+u\cr u
\endsmallmatrix
\right]e_1^{(a_4)}e_3^{(a_2)}e_2^{(a_3+a_4+a_6)}e_4^{(a_1+a_2)}\cr
&\qquad\qquad\times
e_1^{(a_7+a_9)}e_3^{(a_3+a_4+a_5+a_6+u)}e_2^{(a_7+a_8+a_9)}e_4^{(a_3+a_4)}e_1^{(a_{10})}
e_3^{(a_7-u)}\cr &\qquad\qquad\text{\rm if}\qquad a_5+a_6\geq
a_8+a_9\geq a_1+a_2,\quad a_3+a_6\geq a_7+a_9,\cr
&\qquad\qquad\qquad \qquad\qquad\qquad  a_8\geq a_{10},\quad
a_1\geq a_5,\quad a_9\geq a_6.\cr &(3).\;\;\sum_{0\leq u\leq
a_2}\!(-1)^u\left[\smallmatrix a_5-a_1-1+u\cr u
\endsmallmatrix
\right]e_1^{(a_4)}e_3^{(a_2-u)}e_2^{(a_3+a_4+a_6)}e_4^{(a_1+a_2)}e_1^{(a_7+a_9)}\cr
&\qquad\qquad\qquad\qquad\times
e_3^{(a_3+a_4+a_5+a_6+u)}e_2^{(a_7+a_8+a_9)}e_4^{(a_3+a_4)}e_1^{(a_{10})}e_3^{(a_7)}\cr
&\qquad\qquad\text{\rm if}\qquad  a_8+a_9\geq a_5+a_6,\quad
a_3+a_6\geq a_7+a_9,\cr &\qquad\qquad\qquad\qquad   a_8\geq
a_1\geq a_5\geq a_{10},\quad a_6\geq a_2.\cr&\cr {\bold 4.}\quad
&(1).\;\;\sum_{0\leq u\leq
a_2+a_3+a_4}\!\!\!\!\!\!\!\!(-1)^u\left[\smallmatrix a_5-a_1-1+u
\cr u
\endsmallmatrix
\right]e_1^{(a_4)}e_2^{(a_3+a_4)}e_3^{(a_2+a_3+a_4-u)}e_2^{(a_6)}\cr
&\qquad\qquad\times
e_1^{(a_7+a_9)}e_4^{(a_1+a_2+a_3+a_4)}e_2^{(a_7)}e_3^{(a_5+a_6+a_7+u)}e_2^{(a_8+a_9)}
e_1^{(a_{10})}\cr &\qquad\qquad\text{\rm if}\qquad a_3+a_6\geq
a_7+a_9,\quad a_1+a_6\geq a_8+a_9,\quad a_2\geq a_6,\cr
&\qquad\qquad\qquad \qquad a_8\geq a_{10},\quad a_5\geq a_1,\quad
a_9\geq a_6 .\cr&\cr
  &(2).\;\;\sum_{0\leq u\leq
a_7+a_9}\!\!\!\!\!\!\!\!(-1)^u\left[\smallmatrix a_{10}-a_8-1+u
\cr u
\endsmallmatrix
\right]e_1^{(a_4)}e_2^{(a_3+a_4)}e_3^{(a_2+a_3+a_4)}e_2^{(a_6)}e_1^{(a_7+a_9-u)}\cr
&\qquad\qquad\qquad\qquad\qquad\times
e_4^{(a_1+a_2+a_3+a_4)}e_2^{(a_7)}e_3^{(a_5+a_6+a_7)}e_2^{(a_8+a_9)}e_1^{(a_{10}+u)}\cr
&\qquad\qquad\text{\rm if}\qquad a_3+a_6+a_8\geq
a_7+a_9+a_{10},\quad a_5+a_6\geq a_8+a_9,\cr &\qquad\qquad
\qquad\qquad a_{10}\geq a_8,\quad a_9\geq a_6,\quad a_2\geq
a_6,\quad a_1\geq a_5.\cr&\cr &(3).\;\;\sum_{0\leq u\leq
a_4}\!\!\!(-1)^u\left[\smallmatrix a_7+a_9-a_3-a_6-1+u\cr u
\endsmallmatrix
\right]e_1^{(a_4-u)}e_2^{(a_3+a_4)}e_3^{(a_2+a_3+a_4)}e_2^{(a_6)}\cr
&\qquad\qquad\qquad\times
e_1^{(a_7+a_9+u)}e_4^{(a_1+a_2+a_3+a_4)}e_2^{(a_7)}e_3^{(a_5+a_6+a_7)}e_2^{(a_8+a_9)}
e_1^{(a_{10})}\cr &\qquad\qquad\text{\rm if}\qquad a_7+a_9\geq
a_3+a_6,\quad a_5+a_6\geq a_8+a_9,\quad a_1\geq a_5,\cr
&\qquad\qquad\qquad \qquad a_8\geq a_{10},\quad a_2\geq a_6,\quad
a_3\geq a_7 .\cr&\cr \endalign $$$$\align {\bold 5.}\quad
&(1).\;\;\sum_{0\leq u\leq a_3+a_4}\!\!\!(-1)^u\left[\smallmatrix
a_5+a_6-a_1-a_2-1+u\cr u
\endsmallmatrix
\right]e_1^{(a_4)}e_3^{(a_2)}e_2^{(a_3+a_4+a_6)}e_1^{(a_7+a_9)}\cr
&\qquad\qquad\times
e_3^{(a_3+a_4-u)}e_4^{(a_1+a_2+a_3+a_4)}e_2^{(a_7)}e_3^{(a_5+a_6+a_7+u)}e_2^{(a_8+a_9)}
e_1^{(a_{10})}\cr &\qquad\qquad\text{\rm if}\qquad a_5+a_6\geq
a_1+a_2\geq a_8+a_9,\quad a_3+a_6\geq a_7+a_9,\cr
&\qquad\qquad\qquad \qquad  a_8\geq a_{10},\quad a_9\geq a_6,\quad
a_1\geq a_5 .\cr &\cr&(2).\;\;\sum_{0\leq u\leq
a_7+a_9}\!\!\!(-1)^u \left[\smallmatrix a_{10}-a_8-1+u\cr u
\endsmallmatrix
\right]e_1^{(a_4)}e_3^{(a_2)}e_2^{(a_3+a_4+a_6)}e_1^{(a_7+a_9-u)}\cr
&\qquad\qquad\times
e_3^{(a_3+a_4)}e_4^{(a_1+a_2+a_3+a_4)}e_2^{(a_7)}e_3^{(a_5+a_6+a_7)}e_2^{(a_8+a_9)}
e_1^{(a_{10}+u)}\cr &\qquad\qquad\text{\rm if}\qquad
a_3+a_6+a_8\geq a_7+a_9+a_{10},\quad a_9\geq a_6\geq a_2,\cr
&\qquad\qquad \qquad\qquad  a_1+a_2\geq a_5+a_6\geq a_8+a_9,\quad
a_{10}\geq a_8.\cr &(3).\;\;\sum_{0\leq u\leq
a_7}\!\!\!(-1)^u\left[\smallmatrix a_6-a_9-1+u\cr u
\endsmallmatrix
\right]e_1^{(a_4)}e_3^{(a_2)}e_2^{(a_3+a_4+a_6+u)}e_1^{(a_7+a_9)}e_3^{(a_3+a_4)}\cr
&\qquad\qquad\qquad\qquad\times
e_4^{(a_1+a_2+a_3+a_4)}e_2^{(a_7-u)}e_3^{(a_5+a_6+a_7)}e_2^{(a_8+a_9)}e_1^{(a_{10})}\cr
&\qquad\qquad\text{\rm if}\qquad a_1+a_2\geq a_5+a_6,\quad a_6\geq
a_9\geq a_2,\quad a_3\geq a_7,\cr &\qquad\qquad\qquad
\qquad\qquad\qquad a_5\geq a_8\geq a_{10}.\cr &\cr
 {\bold 6.}\quad
&(1).\;\;\sum_{0\leq u\leq a_6}(-1)^u\left[\smallmatrix
a_8-a_5-1+u\cr u
\endsmallmatrix
\right]e_1^{(a_4)}e_2^{(a_3+a_4)}e_3^{(a_2+a_3+a_4)}e_2^{(a_6-u)}\cr
&\qquad\quad\times
e_4^{(a_1+a_2+a_3+a_4)}e_3^{(a_5+a_6)}e_1^{(a_7+a_9)}e_2^{(a_7+a_8+a_9+u)}e_3^{(a_7)}
e_1^{(a_{10})}\cr &\qquad\text{\rm if}\quad a_3+a_6\geq
a_7+a_9,\quad a_2+a_5 \geq a_6+a_8,\quad a_1\geq a_5,\cr
&\qquad\qquad\qquad\qquad a_8\geq a_5\geq a_{10},\quad a_9\geq
a_6.\cr&\cr &(2).\;\;\sum_{0\leq u\leq
a_4}\!\!\!(-1)^u\left[\smallmatrix a_7+a_9-a_3-a_6-1+u\cr u
\endsmallmatrix
\right]e_1^{(a_4-u)}e_2^{(a_3+a_4)}e_3^{(a_2+a_3+a_4)}e_2^{(a_6)}\cr
&\qquad\qquad\times
e_4^{(a_1+a_2+a_3+a_4)}e_3^{(a_5+a_6)}e_1^{(a_7+a_9+u)}e_2^{(a_7+a_8+a_9)}e_3^{(a_7)}
e_1^{(a_{10})}\cr &\quad\text{\rm if}\quad a_3+a_5+a_6\geq
a_7+a_8+a_9,\quad a_7+a_9\geq a_3+a_6,\quad a_1\geq a_5,\cr
&\qquad\qquad\qquad \qquad a_8+a_9\geq a_5+a_6,\quad a_8\geq
a_{10},\quad a_2\geq a_6.\cr&\cr &(3).\;\;\sum_{0\leq u\leq
a_2+a_3+a_4}\!\!\!(-1)^u\left[\smallmatrix a_5-a_1-1+u\cr u
\endsmallmatrix
\right]e_1^{(a_4)}e_2^{(a_3+a_4)}e_3^{(a_2+a_3+a_4-u)}e_2^{(a_6)}\cr
&\qquad\qquad\qquad\times
e_4^{(a_1+a_2+a_3+a_4)}e_3^{(a_5+a_6+u)}e_1^{(a_7+a_9)}e_2^{(a_7+a_8+a_9)}e_3^{(a_7)}
e_1^{(a_{10})}\cr &\qquad\qquad\text{\rm if}\qquad a_8+a_9\geq
a_5+a_6,\quad a_3+a_6\geq a_7+a_9,\cr &\qquad\qquad\qquad\qquad
a_5\geq a_1\geq a_8\geq a_{10},\quad  a_2\geq a_6.\cr
&\cr\endalign$$

$$\align
 {\bold 7.}\quad &(1).\;\;\sum_{0\leq u\leq
a_7}\!\!\!(-1)^u\left[\smallmatrix a_5+a_6-a_8-a_9-1+u\cr u
\endsmallmatrix
\right]e_3^{(a_2)}e_2^{(a_3+a_6)}e_1^{(a_4+a_7+a_9)}e_2^{(a_4)}\cr
&\qquad\qquad\times
e_4^{(a_1+a_2)}e_3^{(a_3+a_4+a_5+a_6+u)}e_2^{(a_7+a_8+a_9)}e_4^{(a_3+a_4)}e_1^{(a_{10})}
e_3^{(a_7-u)}\cr &\qquad\qquad\text{\rm if}\qquad a_5+a_6\geq
a_8+a_9\geq a_1+a_2,\quad a_7+a_9\geq a_3+a_6,\cr
&\qquad\qquad\qquad \qquad  a_1\geq a_5,\quad a_3\geq a_7,\quad
a_8\geq a_{10} .\cr &\cr &(2).\;\;\sum_{0\leq u\leq
a_2}\!\!\!(-1)^u \left[\smallmatrix a_5-a_1-1+u\cr u
\endsmallmatrix
\right]e_3^{(a_2-u)}e_2^{(a_3+a_6)}e_1^{(a_4+a_7+a_9)}e_2^{(a_4)}e_4^{(a_1+a_2)}\cr
&\qquad\qquad\qquad\times
e_3^{(a_3+a_4+a_5+a_6+u)}e_2^{(a_7+a_8+a_9)}e_4^{(a_3+a_4)}e_1^{(a_{10})}e_3^{(a_7)}\cr
&\qquad\text{\rm if}\qquad a_1+a_3+a_6\geq a_7+a_8+a_9,\quad
a_7+a_9\geq a_3+a_6,\quad a_5\geq a_1,\cr &\qquad\qquad\qquad
a_8+a_9\geq a_5+a_6,\quad a_6\geq a_2,\quad a_8\geq a_{10} .\cr
&(3).\;\;\sum_{0\leq u\leq a_4+a_7+a_9}\!\!\!(-1)^u
\left[\smallmatrix a_{10}-a_8-1+u\cr u
\endsmallmatrix
\right]e_3^{(a_2)}e_2^{(a_3+a_6)}e_1^{(a_4+a_7+a_9-u)}e_2^{(a_4)}\cr
&\qquad\qquad\qquad\times
e_4^{(a_1+a_2)}e_3^{(a_3+a_4+a_5+a_6)}e_2^{(a_7+a_8+a_9)}e_4^{(a_3+a_4)}e_1^{(a_{10}+u)}
e_3^{(a_7)}\cr
&\qquad\text{\rm if}\qquad a_3+a_5+a_6\geq a_7+a_8+a_9,\quad
a_7+a_9\geq a_3+a_6,\quad a_1\geq a_5,\cr &\qquad\qquad\qquad
a_8+a_9\geq a_5+a_6\geq a_1+a_2,\quad a_{10}\geq a_8.\cr &\cr
 {\bold 8.}\quad &(1).\sum_{0\leq u\leq a_4}\!\!\!(-1)^u
\left[\smallmatrix a_7+a_9+a_{10}-a_3-a_6-a_8-1+u\cr u
\endsmallmatrix
\right]e_3^{(a_2)}e_1^{(a_4-u)}e_4^{(a_1+a_2)}\cr &\quad \times
e_2^{(a_3+a_4+a_6)}e_3^{(a_3+a_4+a_5+a_6)}e_4^{(a_3+a_4)}e_2^{(a_8)}e_1^{(a_7+a_9+a_{10}+u)}
e_2^{(a_7+a_9)}e_3^{(a_7)}\cr &\qquad\text{\rm if}\qquad
a_7+a_9+a_{10}\geq a_3+a_6+a_8,\quad a_3+a_6\geq a_7+a_9,\cr
&\qquad\qquad\qquad a_8+a_9\geq a_5+a_6\geq a_1+a_2,\quad a_1\geq
a_5\geq a_8.\cr &\cr&(2).\;\;\sum_{0\leq u\leq
a_2}(-1)^u\left[\smallmatrix a_5-a_1-1+u\cr u
\endsmallmatrix
\right]e_3^{(a_2-u)}e_1^{(a_4)}e_4^{(a_1+a_2)}e_2^{(a_3+a_4+a_6)}\cr
&\quad\quad\times
e_3^{(a_3+a_4+a_5+a_6+u)}e_4^{(a_3+a_4)}e_2^{(a_8)}e_1^{(a_7+a_9+a_{10})}e_2^{(a_7+a_9)}
e_3^{(a_7)}\cr &\qquad\text{\rm if}\qquad a_3+a_6+a_8\geq
a_7+a_9+a_{10},\quad a_8+a_9\geq a_5+a_6,\cr &\qquad\qquad \qquad
a_5\geq a_1\geq a_8,\quad a_6\geq a_2,\quad a_{10}\geq a_8
.\cr&\cr &(3).\;\;\sum_{0\leq u\leq
a_7}\!\!\!(-1)^u\left[\smallmatrix a_5+a_6-a_8-a_9-1+u\cr u
\endsmallmatrix
\right]e_3^{(a_2)}e_1^{(a_4)}e_4^{(a_1+a_2)}e_2^{(a_3+a_4+a_6)}\cr
&\qquad\quad\times
e_3^{(a_3+a_4+a_5+a_6+u)}e_4^{(a_3+a_4)}e_2^{(a_8)}e_1^{(a_7+a_9+a_{10})}e_2^{(a_7+a_9)}
e_3^{(a_7-u)}\cr &\qquad\text{\rm if}\qquad a_3+a_6+a_8\geq
a_7+a_9+a_{10},\quad a_1\geq a_5,\quad a_{10}\geq a_8,\cr &
\qquad\qquad\qquad a_5+a_6\geq a_8+a_9\geq a_1+a_2,\quad a_9\geq
a_6.\cr &\cr\endalign $$

$$\align {\bold 9.}\quad &(1).\!\!\sum_{0\leq
u\leq a_6+a_7}\!\!\!(-1)^u \left[\smallmatrix a_8-a_5-1+u\cr u
\endsmallmatrix
\right]e_1^{(a_4)}e_2^{(a_3+a_4)}e_3^{(a_2+a_3+a_4)}e_1^{(a_7)}e_2^{(a_6+a_7-u)}\cr
&\qquad\qquad\qquad\quad\quad\times
e_1^{(a_9)}e_4^{(a_1+a_2+a_3+a_4)}e_3^{(a_5+a_6+a_7)}e_2^{(a_8+a_9+u)}e_1^{(a_{10})}\cr
&\qquad\qquad\qquad\text{\rm if}\qquad a_2+a_5\geq a_6+a_8,\quad
a_1\geq a_5\geq a_{10},\cr &\qquad\qquad \qquad\qquad\qquad
a_8\geq a_5,\quad a_3\geq a_7,\quad a_6\geq a_9.\cr &\cr
&(2).\sum_{0\leq u\leq a_2+a_3+a_4}\!\!\!(-1)^u \left[\smallmatrix
a_5-a_1-1+u\cr u
\endsmallmatrix
\right]e_1^{(a_4)}e_2^{(a_3+a_4)}e_3^{(a_2+a_3+a_4-u)}e_1^{(a_7)}\cr
&\qquad\qquad\times
e_2^{(a_6+a_7)}e_1^{(a_9)}e_4^{(a_1+a_2+a_3+a_4)}e_3^{(a_5+a_6+a_7+u)}e_2^{(a_8+a_9)}
e_1^{(a_{10})}\cr
&\qquad\qquad\text{\rm if}\qquad a_5\geq a_1\geq a_8\geq
 a_{10},\quad a_2\geq a_6\geq a_9,\quad a_3\geq a_7.\cr
&(3).\sum_{0\leq u\leq a_3+a_4}\!\!\!(-1)^u\left[\smallmatrix
a_6-a_2-1+u\cr u
\endsmallmatrix
\right]e_1^{(a_4)}e_2^{(a_3+a_4-u)}e_3^{(a_2+a_3+a_4)}e_1^{(a_7)}\cr
&\qquad\qquad\times
e_2^{(a_6+a_7+u)}e_1^{(a_9)}e_4^{(a_1+a_2+a_3+a_4)}e_3^{(a_5+a_6+a_7)}e_2^{(a_8+a_9)}
e_1^{(a_{10})}\cr &\qquad\qquad\qquad\text{\rm if}\qquad
a_1+a_2\geq a_5+a_6,\quad a_5\geq a_8\geq a_{10},\cr &\qquad\qquad
\qquad\qquad\qquad a_6\geq a_2\geq a_9,\quad a_3\geq a_7.\cr &\cr
{\bold {10.}}\quad &(1).\;\;\sum_{0\leq u\leq a_3+a_4}\!\!(-1)^u
\left[\smallmatrix a_6-a_2-1+u\cr u
\endsmallmatrix
\right]e_1^{(a_4)}e_2^{(a_3+a_4-u)}e_3^{(a_2+a_3+a_4)}e_1^{(a_7)}\cr
&\qquad\quad\times
e_2^{(a_6+a_7+u)}e_4^{(a_1+a_2+a_3+a_4)}e_3^{(a_5+a_6+a_7)}e_2^{(a_8)}e_1^{(a_9+a_{10})}
e_2^{(a_9)}\cr &\quad\qquad\qquad\text{\rm if}\quad a_1+a_2\geq
a_5+a_6,\quad a_2+a_8\geq a_9+a_{10},\cr &\qquad\qquad\quad \qquad
a_{10}\geq a_8,\quad a_5\geq a_8,\quad a_6\geq a_2,\quad a_3\geq
a_7.\cr &\cr&(2).\;\;\sum_{0\leq u\leq a_7}\!\!(-1)^u
\left[\smallmatrix a_9+a_{10}-a_6-a_8-1+u\cr u
\endsmallmatrix
\right]e_1^{(a_4)}e_2^{(a_3+a_4)}e_3^{(a_2+a_3+a_4)}e_1^{(a_7-u)}\cr
&\qquad\qquad\qquad\times
e_2^{(a_6+a_7)}e_4^{(a_1+a_2+a_3+a_4)}e_3^{(a_5+a_6+a_7)}e_2^{(a_8)}e_1^{(a_9+a_{10}+u)}
e_2^{(a_9)}\cr &\quad\qquad\qquad\text{\rm if}\quad
a_3+a_6+a_8\geq a_7+a_9+a_{10},\quad a_9+a_{10}\geq a_6+a_8,\cr
&\qquad\qquad \qquad\quad\qquad a_1\geq a_5\geq a_8,\quad a_2\geq
a_6\geq a_9 .\cr&\cr &(3).\;\;\sum_{0\leq u\leq
a_2+a_3+a_4}\!\!(-1)^u\left[\smallmatrix a_5-a_1-1+u\cr u
\endsmallmatrix
\right]e_1^{(a_4)}e_2^{(a_3+a_4)}e_3^{(a_2+a_3+a_4-u)}e_1^{(a_7)}\cr
&\qquad\qquad\qquad\times
e_2^{(a_6+a_7)}e_4^{(a_1+a_2+a_3+a_4)}e_3^{(a_5+a_6+a_7+u)}e_2^{(a_8)}e_1^{(a_9+a_{10})}
e_2^{(a_9)}\cr &\quad\qquad\qquad\text{\rm if}\quad a_6+a_8\geq
a_9+a_{10},\quad a_5\geq a_1\geq a_8,\cr &\qquad\qquad\qquad
\quad\qquad a_{10}\geq a_8,\quad a_3\geq a_7,\quad a_2\geq a_6
.\cr &\cr\endalign $$

$$\align{\bold {11.}}\quad
&(1).\;\;\sum_{0\leq u\leq a_3+a_6}\!\!(-1)^u\left[\smallmatrix
a_5-a_8-1+u\cr u
\endsmallmatrix
\right]e_4^{(a_1)}e_3^{(a_2+a_5+u)}e_2^{(a_3+a_6+a_8)}e_4^{(a_2)}\cr
&\qquad\qquad\times
e_3^{(a_3+a_6-u)}e_1^{(a_4+a_7+a_9+a_{10})}e_2^{(a_4+a_7+a_9)}e_3^{(a_4)}e_4^{(a_3+a_4)}
e_3^{(a_7)}\cr &\quad\qquad\qquad\qquad\text{\rm if}\quad
a_7+a_8+a_9\geq a_3+a_5+a_6,\quad a_5\geq a_8,\cr
&\qquad\qquad\qquad\quad \qquad a_{10}\geq a_8\geq a_1,\quad
a_6\geq a_2,\quad a_3\geq a_7 .\cr &\cr&(2).\;\;\sum_{0\leq u\leq
a_3+a_4}\!\!(-1)^u \left[\smallmatrix a_2-a_6-1+u\cr u
\endsmallmatrix
\right]e_4^{(a_1)}e_3^{(a_2+a_5)}e_2^{(a_3+a_6+a_8)}e_4^{(a_2+u)}\cr
&\qquad\qquad\times
e_3^{(a_3+a_6)}e_1^{(a_4+a_7+a_9+a_{10})}e_2^{(a_4+a_7+a_9)}e_3^{(a_4)}e_4^{(a_3+a_4-u)}
e_3^{(a_7)}\cr &\quad\qquad\text{\rm if}\quad a_7+a_9\geq
a_3+a_6,\quad a_6+a_8\geq a_2+a_5,\quad a_3\geq a_7,\cr
&\qquad\qquad\qquad \quad\qquad a_{10}\geq a_8,\quad  a_2\geq
a_6,\quad a_5\geq a_1 .\cr &(3).\;\;\sum_{0\leq u\leq
a_4+a_7+a_9}\!\!(-1)^u\left[\smallmatrix a_8-a_{10}-1+u\cr u
\endsmallmatrix
\right]e_4^{(a_1)}e_3^{(a_2+a_5)}e_2^{(a_3+a_6+a_8+u)}e_4^{(a_2)}\cr
&\qquad\qquad\qquad\times
e_3^{(a_3+a_6)}e_1^{(a_4+a_7+a_9+a_{10})}e_2^{(a_4+a_7+a_9-u)}e_3^{(a_4)}e_4^{(a_3+a_4)}
e_3^{(a_7)}\cr
&\quad\qquad\qquad\qquad\text{\rm if}\quad a_7+a_9\geq
a_3+a_6,\quad a_6\geq a_2,\cr &\qquad\qquad\qquad\quad\qquad
a_8\geq a_{10}\geq a_5\geq a_1,\quad  a_3\geq a_7.\cr &\cr
 {\bold {12.}}\quad &(1).\sum_{0\leq u\leq a_4+a_7}\!\!\!\!\!\!(-1)^u
\left[\smallmatrix a_5+a_6-a_8-a_9-1+u\cr u
\endsmallmatrix
\right]e_3^{(a_2)}e_2^{(a_3+a_6)}e_1^{(a_4+a_7+a_9)}e_4^{(a_1+a_2)}\cr
&\qquad\qquad\times
e_3^{(a_3+a_5+a_6+u)}e_2^{(a_4+a_7+a_8+a_9)}e_4^{(a_3)}e_3^{(a_4+a_7-u)}e_4^{(a_4)}
e_1^{(a_{10})}\cr &\qquad\qquad\qquad\text{\rm if}\quad
a_5+a_6\geq a_8+a_9\geq a_1+a_2,\quad a_7\geq a_3,\cr
&\qquad\qquad\qquad\qquad a_8\geq a_{10},\quad  a_1\geq a_5,\quad
a_9\geq a_6.\cr&\cr &(2).\;\;\sum_{0\leq u\leq
a_2}\!\!(-1)^u\left[\smallmatrix a_5-a_1-1+u\cr u
\endsmallmatrix
\right]e_3^{(a_2-u)}e_2^{(a_3+a_6)}e_1^{(a_4+a_7+a_9)}e_4^{(a_1+a_2)}\cr
&\qquad\qquad\times
e_3^{(a_3+a_5+a_6+u)}e_2^{(a_4+a_7+a_8+a_9)}e_4^{(a_3)}e_3^{(a_4+a_7)}e_4^{(a_4)}
e_1^{(a_{10})}\cr &\qquad\qquad\qquad\text{\rm if}\quad
a_8+a_9\geq a_5+a_6,\quad a_7\geq a_3,\cr
&\qquad\qquad\qquad\qquad a_5\geq a_1\geq a_8\geq a_{10},\quad
a_6\geq a_2.\cr &\cr &(3).\;\;\sum_{0\leq u\leq a_3+a_6}\!\!(-1)^u
\left[\smallmatrix a_8-a_5-1+u\cr u
\endsmallmatrix
\right]e_3^{(a_2)}e_2^{(a_3+a_6-u)}e_1^{(a_4+a_7+a_9)}e_4^{(a_1+a_2)}\cr
&\qquad\qquad\times
e_3^{(a_3+a_5+a_6)}e_2^{(a_4+a_7+a_8+a_9+u)}e_4^{(a_3)}e_3^{(a_4+a_7)}e_4^{(a_4)}
e_1^{(a_{10})}\cr &\qquad\qquad\qquad\text{\rm if}\quad
a_5+a_6\geq a_1+a_2,\quad a_1\geq a_5\geq a_{10},\cr
&\qquad\qquad\qquad \qquad a_8\geq a_5,\quad  a_7\geq a_3,\quad
a_9\geq a_6.\cr &\cr \endalign $$$$\align{\bold {13.}}\quad
&(1).\;\;\sum_{0\leq u\leq
a_4+a_6+a_7}\!\!(-1)^u\left[\smallmatrix a_8-a_5-1+u\cr u
\endsmallmatrix
\right]e_2^{(a_3)}e_3^{(a_2+a_3)}e_4^{(a_1+a_2+a_3)}e_1^{(a_4+a_7)}\cr
&\qquad\qquad\times
e_2^{(a_4+a_6+a_7-u)}e_3^{(a_4+a_5+a_6+a_7)}e_1^{(a_9)}e_2^{(a_8+a_9+u)}e_1^{(a_{10})}
e_4^{(a_4)}\cr &\qquad\qquad\qquad\text{\rm if}\quad
a_5+a_6+a_7\geq a_1+a_2+a_3,\quad  a_2+a_5\geq a_6+a_8,\cr
&\qquad\qquad\qquad \qquad a_1\geq a_5\geq a_{10},\quad  a_6\geq
a_9,\quad a_8\geq a_5 .\cr &\cr&(2).\;\;\sum_{0\leq u\leq
a_3}\!\!(-1)^u \left[\smallmatrix a_6-a_2-1+u\cr u
\endsmallmatrix
\right]e_2^{(a_3-u)}e_3^{(a_2+a_3)}e_4^{(a_1+a_2+a_3)}e_1^{(a_4+a_7)}\cr
&\qquad\qquad\times
e_2^{(a_4+a_6+a_7+u)}e_3^{(a_4+a_5+a_6+a_7)}e_1^{(a_9)}e_2^{(a_8+a_9)}e_1^{(a_{10})}
e_4^{(a_4)}\cr &\qquad\qquad\qquad\text{\rm if}\quad
a_5+a_6+a_7\geq a_1+a_2+a_3,\quad  a_1+a_2\geq a_5+a_6,\cr
&\qquad\qquad\qquad \qquad a_5\geq a_8\geq a_{10},\quad  a_6\geq
a_2\geq a_9.\cr &(3).\;\;\sum_{0\leq u\leq a_2+a_3}\!\!(-1)^u
\left[\smallmatrix a_5-a_1-1+u\cr u
\endsmallmatrix
\right]e_2^{(a_3)}e_3^{(a_2+a_3-u)}e_4^{(a_1+a_2+a_3)}e_1^{(a_4+a_7)}\cr
&\qquad\qquad\times
e_2^{(a_4+a_6+a_7)}e_3^{(a_4+a_5+a_6+a_7+u)}e_1^{(a_9)}e_2^{(a_8+a_9)}e_1^{(a_{10})}
e_4^{(a_4)}\cr
&\qquad\qquad\qquad\text{\rm if}\quad a_6+a_7\geq
a_2+a_3,\quad a_2\geq a_6\geq a_9,\cr &\qquad\qquad\qquad
\qquad a_5\geq a_1\geq a_8\geq a_{10}.\cr &\cr
 {\bold {14.}}\quad &(1).\sum_{0\leq u\leq
a_4+a_7+a_9}\!\!\!\!\!\!\!\!\!\!\!(-1)^u\left[\smallmatrix
a_{10}-a_8-1+u\cr u
\endsmallmatrix
\right]e_3^{(a_2)}e_2^{(a_3+a_6)}e_1^{(a_4+a_7+a_9-u)}e_2^{(a_4)}\cr
&\qquad\qquad\times
e_3^{(a_3+a_4)}e_4^{(a_1+a_2+a_3+a_4)}e_3^{(a_5+a_6)}e_2^{(a_7+a_8+a_9)}e_1^{(a_{10}+u)}
e_3^{(a_7)}\cr &\quad\text{\rm if}\quad a_3+a_5+a_6\geq
a_7+a_8+a_9,\quad a_7+a_9\geq a_3+a_6,\quad a_6\geq a_2,\cr
&\qquad\qquad \quad   a_8+a_9\geq a_5+a_6,\quad a_1+a_2\geq
a_5+a_6 ,\quad a_{10}\geq a_8.\cr &\cr&(2).\;\;\sum_{0\leq u\leq
a_4}\!\!\!(-1)^u\left[\smallmatrix a_7+a_8+a_9-a_3-a_5-a_6-1+u \cr
u
\endsmallmatrix
\right]e_3^{(a_2)}e_2^{(a_3+a_6)}e_1^{(a_4+a_7+a_9)}\cr &\qquad
\quad\times
e_2^{(a_4-u)}e_3^{(a_3+a_4)}e_4^{(a_1+a_2+a_3+a_4)}e_3^{(a_5+a_6)}e_2^{(a_7+a_8+a_9+u)}
e_1^{(a_{10})}e_3^{(a_7)}\cr &\quad\text{\rm if}\quad
a_7+a_8+a_9\geq a_3+a_5+a_6,\quad a_1+a_2\geq a_5+a_6,\quad
a_3\geq a_7,\cr &\qquad\qquad \qquad\quad   a_5\geq a_8\geq
a_{10},\quad a_6\geq a_2.\cr&\cr {\bold
{15.}}\quad&(1).\sum_{0\leq u\leq a_4}\!\!\!(-1)^u
\left[\smallmatrix a_7+a_9+a_{10}-a_3-a_6-a_8-1+u\cr u
\endsmallmatrix
\right]e_3^{(a_2)}e_1^{(a_4-u)}e_2^{(a_3+a_4+a_6)}\cr &\quad
\quad\times
e_3^{(a_3+a_4)}e_4^{(a_1+a_2+a_3+a_4)}e_3^{(a_5+a_6)}e_2^{(a_8)}e_1^{(a_7+a_9+a_{10}+u)}
e_2^{(a_7+a_9)}e_3^{(a_7)}\cr &\quad\text{\rm if}\quad
a_7+a_9+a_{10}\geq a_3+a_6+a_8,\quad a_1+a_2\geq a_5+a_6,\quad
a_6\geq a_2,\cr &\qquad\qquad\qquad \quad a_8+a_9\geq
a_5+a_6,\quad a_3+a_6\geq a_7+a_9,\quad a_5\geq a_8.\cr &\cr
\endalign $$$$\align &(2).\;\;\sum_{0\leq u\leq
a_3+a_4+a_6}\!\!\!(-1)^u \left[\smallmatrix a_8-a_5-1+u\cr u
\endsmallmatrix
\right]e_3^{(a_2)}e_1^{(a_4)}e_2^{(a_3+a_4+a_6-u)}e_3^{(a_3+a_4)}\cr
&\qquad\qquad\qquad\times
e_4^{(a_1+a_2+a_3+a_4)}e_3^{(a_5+a_6)}e_2^{(a_8+u)}e_1^{(a_7+a_9+a_{10})}e_2^{(a_7+a_9)}
e_3^{(a_7)}\cr &\qquad\text{\rm if}\qquad a_3+a_6+a_8\geq
a_7+a_9+a_{10},\quad a_1+a_2\geq a_5+a_6,\cr &\qquad\qquad
\qquad\quad a_9\geq a_6\geq a_2,\quad a_{10}\geq a_8\geq
a_5.\cr&\cr  {\bold {16.}}\quad &(1).\sum_{0\leq u\leq
a_3+a_4}\!\!\!(-1)^u \left[\smallmatrix a_5+a_6-a_1-a_2-1+u\cr u
\endsmallmatrix
\right]e_3^{(a_2)}e_2^{(a_3+a_6)}e_1^{(a_4+a_7+a_9)}e_2^{(a_4)}\cr
&\qquad\quad\times
e_3^{(a_3+a_4-u)}e_4^{(a_1+a_2+a_3+a_4)}e_2^{(a_7)}e_3^{(a_5+a_6+a_7+u)}e_2^{(a_8+a_9)}
e_1^{(a_{10})}\cr &\qquad\text{\rm if}\qquad a_5+a_6\geq
a_1+a_2\geq a_8+a_9,\quad a_7+a_9\geq a_3+a_6,\cr
&\qquad\qquad\qquad \quad a_8\geq a_{10},\quad a_3\geq a_7,\quad
a_1\geq a_5.\cr &(2).\sum_{0\leq u\leq a_4+a_7+a_9}\!\!\!(-1)^u
\left[\smallmatrix a_{10}-a_8-1+u\cr u
\endsmallmatrix
\right]e_3^{(a_2)}e_2^{(a_3+a_6)}e_1^{(a_4+a_7+a_9-u)}e_2^{(a_4)}\cr
&\qquad\quad\times
e_3^{(a_3+a_4)}e_4^{(a_1+a_2+a_3+a_4)}e_2^{(a_7)}e_3^{(a_5+a_6+a_7)}e_2^{(a_8+a_9)}
e_1^{(a_{10}+u)}\cr &\qquad\text{\rm if}\qquad a_1+a_2\geq
a_5+a_6\geq a_8+a_9,\quad a_7+a_9\geq a_3+a_6,\cr
&\qquad\qquad\qquad \quad a_{10}\geq a_8,\quad a_3\geq a_7,\quad
a_6\geq a_2.\cr&\cr {\bold {17.}}\quad&(1).\sum_{0\leq u\leq
a_4}\!\!\!(-1)^u \left[\smallmatrix a_1+a_2+a_3-a_5-a_6-a_7-1+u\cr
u
\endsmallmatrix
\right]e_2^{(a_3)}e_3^{(a_2+a_3)}e_4^{(a_1+a_2+a_3+u)}\cr &\qquad
\quad\times
e_3^{(a_5)}e_2^{(a_6+a_8)}e_1^{(a_4+a_7+a_9+a_{10})}e_3^{(a_6)}e_2^{(a_4+a_7+a_9)}
e_3^{(a_4+a_7)}e_4^{(a_4-u)}\cr &\quad\text{\rm if}\quad
a_1+a_2+a_3\geq a_5+a_6+a_7,\quad a_2+a_5\geq a_6+a_8,\quad
a_9\geq a_6,\cr &\qquad\qquad\qquad \quad a_6+a_7\geq
a_2+a_3,\quad  a_{10}\geq a_8\geq  a_5.\cr
&\cr&(2).\;\;\sum_{0\leq u\leq a_4+a_7+a_9}\!\!\!(-1)^u
\left[\smallmatrix a_8-a_{10}-1+u\cr u
\endsmallmatrix
\right]e_2^{(a_3)}e_3^{(a_2+a_3)}e_4^{(a_1+a_2+a_3)}e_3^{(a_5)}\cr
&\qquad\times
e_2^{(a_6+a_8+u)}e_1^{(a_4+a_7+a_9+a_{10})}e_3^{(a_6)}e_2^{(a_4+a_7+a_9-u)}e_3^{(a_4+a_7)}
e_4^{(a_4)}\cr &\qquad\text{\rm if}\qquad a_5+a_6+a_7\geq
a_1+a_2+a_3,\quad a_2+a_5\geq a_6+a_8,\cr &\qquad\qquad\qquad\quad
a_1\geq a_5,\quad a_8\geq a_{10}\geq a_5,\quad a_9\geq a_6.\cr
&\cr  {\bold {18.}}\quad &(1).\;\;\sum_{0\leq u \leq
a_3}\!\!(-1)^u\left[\smallmatrix a_6-a_2-1+u\cr u
\endsmallmatrix
\right]e_2^{(a_3-u)}e_3^{(a_2+a_3)}e_4^{(a_1+a_2+a_3)}e_1^{(a_4+a_7)}\cr
&\qquad\qquad\qquad\times
e_2^{(a_4+a_6+a_7+u)}e_3^{(a_4+a_5+a_6+a_7)}e_2^{(a_8)}e_1^{(a_9+a_{10})}e_4^{(a_4)}
e_2^{(a_9)}\cr
&\quad\text{\rm if}\quad a_5+a_6+a_7\geq a_1+a_2+a_3,\quad
a_2+a_8\geq a_9+a_{10},\quad a_5\geq a_8,\cr &\qquad\qquad
\qquad a_1+a_2\geq a_5+a_6,\quad a_{10}\geq a_8,\quad a_6\geq a_2
.\cr &\cr \endalign $$$$\align &(2).\;\;\sum_{0\leq u\leq a_2+a_3}\!\!(-1)^u
\left[\smallmatrix a_5-a_1-1+u\cr u
\endsmallmatrix
\right]e_2^{(a_3)}e_3^{(a_2+a_3-u)}e_4^{(a_1+a_2+a_3)}e_1^{(a_4+a_7)}\cr
&\qquad\qquad\qquad\times
e_2^{(a_4+a_6+a_7)}e_3^{(a_4+a_5+a_6+a_7+u)}e_2^{(a_8)}e_1^{(a_9+a_{10})}e_4^{(a_4)}
e_2^{(a_9)}\cr &\qquad\qquad\text{\rm if}\qquad a_6+a_7\geq
a_2+a_3,\quad a_6+a_8\geq a_9+a_{10},\cr &\qquad\qquad\qquad\qquad
a_5\geq a_1\geq a_8,\quad a_2\geq a_6,\quad a_{10}\geq a_8.\cr
&\cr {\bold {19.}}\quad &(1).\sum_{0\leq u \leq
a_6+a_7}\!\!\!\!\!\!\!\!(-1)^u\left[\smallmatrix a_5-a_8-1+u \cr u
\endsmallmatrix
\right]e_2^{(a_3)}e_4^{(a_1)}e_1^{(a_4+a_7)}e_2^{(a_4)}e_3^{(a_2+a_3+a_4+a_5+u)}\cr
&\qquad\qquad\qquad\times
e_2^{(a_6+a_7+a_8)}e_1^{(a_9+a_{10})}e_4^{(a_2+a_3+a_4)}e_3^{(a_6+a_7-u)}e_2^{(a_9)}\cr
&\quad\qquad\qquad\text{\rm if}\quad a_2+a_3\geq a_6+a_7,\quad
a_6+a_8\geq a_9+a_{10},\quad a_5\geq a_8,\cr &
\qquad\quad\qquad\qquad\qquad  a_{10}\geq a_8\geq a_1,\quad
a_7\geq a_3.\cr &(2).\;\;\sum_{0\leq u\leq a_9}\!\!(-1)^u
\left[\smallmatrix a_8-a_{10}-1+u\cr u
\endsmallmatrix
\right]e_2^{(a_3)}e_4^{(a_1)}e_1^{(a_4+a_7)}e_2^{(a_4)}e_3^{(a_2+a_3+a_4+a_5)}\cr
&\qquad\qquad\qquad\times
e_2^{(a_6+a_7+a_8+u)}e_1^{(a_9+a_{10})}e_4^{(a_2+a_3+a_4)}e_3^{(a_6+a_7)}e_2^{(a_9-u)}\cr
&\quad\qquad\qquad\text{\rm if}\quad a_2+a_3+a_5\geq
a_6+a_7+a_8,\quad a_6\geq a_9,\cr &\qquad\quad\qquad\qquad \qquad
a_8\geq  a_{10}\geq a_5\geq a_1,\quad a_7\geq a_3.\cr &\cr {\bold
{20.}}\quad &(1).\sum_{0\leq u\leq a_4+a_7+a_9}\!\!(-1)^u
\left[\smallmatrix a_{10}-a_8-1+u\cr u
\endsmallmatrix
\right]e_2^{(a_3)}e_3^{(a_2+a_3)}e_2^{(a_6)}e_1^{(a_4+a_7+a_9-u)}\cr
&\quad\quad\quad\times
e_2^{(a_4+a_7)}e_3^{(a_4)}e_4^{(a_1+a_2+a_3+a_4)}e_3^{(a_5+a_6+a_7)}e_2^{(a_8+a_9)}
e_1^{(a_{10}+u)}\cr &\qquad\text{\rm if}\quad a_1+a_2+a_3\geq
a_5+a_6+a_7,\quad a_6+a_7\geq a_2+a_3,\quad a_2\geq a_6,\cr
&\qquad\qquad\qquad a_5+a_6\geq a_8+a_9,\quad a_{10}\geq a_8,\quad
a_9\geq a_6.\cr &\cr &(2).\sum_{0\leq u\leq
a_4+a_7}\!\!(-1)^u\left[\smallmatrix a_8+a_9-a_5-a_6-1+u\cr u
\endsmallmatrix
\right]e_2^{(a_3)}e_3^{(a_2+a_3)}e_2^{(a_6)}e_1^{(a_4+a_7+a_9)}\cr
&\qquad\quad\quad\times
e_2^{(a_4+a_7-u)}e_3^{(a_4)}e_4^{(a_1+a_2+a_3+a_4)}e_3^{(a_5+a_6+a_7)}e_2^{(a_8+a_9+u)}
e_1^{(a_{10})}\cr &\qquad\text{\rm if}\quad a_1+a_2+a_3\geq
a_5+a_6+a_7,\quad a_6+a_7\geq a_2+a_3,\quad a_2\geq a_6,\cr
&\qquad\qquad\qquad a_8+a_9\geq a_5+a_6,\quad a_5\geq a_8\geq
a_{10}.\cr &\cr {\bold {21.}}\quad &(1).\!\!\!\sum_{0\leq u\leq
a_4}\!\!\!\!(-1)^u\!\left[\!\smallmatrix
a_1+a_2+a_3-a_5-a_6-a_7-1+u\cr u
\endsmallmatrix
\!\right]\!\!
e_3^{(a_2)}e_2^{(a_3+a_6)}e_3^{(a_3)}e_4^{(a_1+a_2+a_3+u)}\cr
&\quad\qquad\qquad\times
e_3^{(a_5+a_6)}e_2^{(a_8)}e_1^{(a_4+a_7+a_9+a_{10})}e_2^{(a_4+a_7+a_9)}e_3^{(a_4+a_7)}
e_4^{(a_4-u)}\cr
&\qquad\text{\rm if}\quad a_1+a_2+a_3\geq a_5+a_6+a_7,\quad
a_8+a_9\geq a_5+a_6,\quad a_5\geq a_8,\cr &\qquad\qquad\qquad
a_7\geq a_3,\quad a_6\geq a_2,\quad a_{10}\geq a_8.\cr &\cr
\endalign $$$$\align &(2).\sum_{0\leq u\leq a_4+a_7}\!\!(-1)^u\left[\smallmatrix
a_5+a_6-a_8-a_9-1+u\cr u
\endsmallmatrix
\right]e_3^{(a_2)}e_2^{(a_3+a_6)}e_3^{(a_3)}e_4^{(a_1+a_2+a_3)}\cr
&\quad\quad\times
e_3^{(a_5+a_6+u)}e_2^{(a_8)}e_1^{(a_4+a_7+a_9+a_{10})}e_2^{(a_4+a_7+a_9)}e_3^{(a_4+a_7-u)}
e_4^{(a_4)}\cr &\qquad\qquad\qquad\text{\rm if}\quad
a_5+a_6+a_7\geq a_1+a_2+a_3,\quad a_{10}\geq a_8,\cr
&\qquad\qquad\qquad\qquad a_1+a_2\geq a_5+a_6\geq a_8+a_9,\quad
a_9\geq a_6\geq a_2.\cr &\cr  {\bold {22.}}\quad
&(1).\;\;\sum_{0\leq u\leq a_3}\!\!(-1)^u\left[\smallmatrix
a_5+a_6-a_1-a_2-1+u\cr u
\endsmallmatrix
\right]e_3^{(a_2)}e_2^{(a_3+a_6)}e_1^{(a_4+a_7+a_9)}e_3^{(a_3-u)}\cr
&\qquad\qquad\times
e_4^{(a_1+a_2+a_3)}e_2^{(a_4+a_7)}e_3^{(a_4+a_5+a_6+a_7+u)}e_2^{(a_8+a_9)}e_4^{(a_4)}
e_1^{(a_{10})}\cr &\qquad\qquad\qquad\text{\rm if}\quad
a_5+a_6\geq a_1+a_2 \geq a_8+a_9,\quad a_7\geq a_3,\cr
&\qquad\qquad\qquad\qquad a_9\geq a_6,\quad a_8\geq a_{10},\quad
a_1\geq a_5.\cr &(2).\;\;\sum_{0\leq u\leq
a_4+a_7}\!\!(-1)^u\left[\smallmatrix a_6-a_9-1+u\cr u
\endsmallmatrix
\right]e_3^{(a_2)}e_2^{(a_3+a_6+u)}e_1^{(a_4+a_7+a_9)}e_3^{(a_3)}\cr
&\qquad\qquad\times
e_4^{(a_1+a_2+a_3)}e_2^{(a_4+a_7-u)}e_3^{(a_4+a_5+a_6+a_7)}e_2^{(a_8+a_9)}e_4^{(a_4)}
e_1^{(a_{10})}\cr &\qquad\qquad\qquad\text{\rm if}\quad
a_5+a_6+a_7\geq a_1+a_2+a_3,\quad a_1+a_2\geq a_5+a_6,\cr
&\qquad\qquad\qquad \qquad a_5\geq  a_8\geq a_{10},\quad a_6\geq
a_9\geq a_2.\cr &\cr {\bold {23.}}\quad &(1).\sum_{0\leq u\leq
a_2+a_3}\!\!(-1)^u \left[\smallmatrix a_5-a_1-1+u\cr u
\endsmallmatrix
\right]e_2^{(a_3)}e_3^{(a_2+a_3-u)}e_4^{(a_1+a_2+a_3)}e_2^{(a_6)}\cr
&\quad\quad\times
e_1^{(a_4+a_7+a_9)}e_3^{(a_5+a_6+u)}e_2^{(a_4+a_7+a_8+a_9)}e_3^{(a_4+a_7)}e_1^{(a_{10})}
e_4^{(a_4)}\cr &\qquad\qquad\qquad\text{\rm if}\quad a_6+a_7\geq
a_2+a_3,\quad a_8+a_9\geq a_5+a_6,\cr &\qquad\qquad\qquad \qquad
a_5\geq a_1\geq  a_8\geq a_{10},\quad a_2\geq a_6.\cr &\cr
&(2).\;\;\sum_{0\leq u\leq a_6}\!\!(-1)^u\left[\smallmatrix
a_8-a_5-1+u\cr u
\endsmallmatrix
\right]e_2^{(a_3)}e_3^{(a_2+a_3)}e_4^{(a_1+a_2+a_3)}e_2^{(a_6-u)}\cr
&\quad\quad\times
e_1^{(a_4+a_7+a_9)}e_3^{(a_5+a_6)}e_2^{(a_4+a_7+a_8+a_9+u)}e_3^{(a_4+a_7)}e_1^{(a_{10})}
e_4^{(a_4)}\cr &\quad\quad\text{\rm if}\quad a_5+a_6+a_7\geq
a_1+a_2+a_3,\quad a_2+a_5\geq a_6+a_8,\cr &\qquad\qquad\qquad
\qquad a_1\geq a_5\geq a_{10},\quad a_9\geq a_6,\quad a_8\geq a_5
.\cr&\cr  {\bold {24.}}\quad &(1).\!\!\!\sum_{0\leq u\leq
a_3+a_4}\!\!\!\!(-1)^u\!\!\left[\!\smallmatrix
a_1+a_2-a_5-a_6-1+u\cr u
\endsmallmatrix
\!\right]\!e_3^{(a_2)}e_2^{(a_3+a_6)}e_4^{(a_1+a_2+u)}e_3^{(a_3+a_5+a_6)}\cr
&\qquad\qquad\qquad\times
e_2^{(a_8)}e_1^{(a_4+a_7+a_9+a_{10})}e_2^{(a_4+a_7+a_9)}e_3^{(a_4)}e_4^{(a_3+a_4-u)}
e_3^{(a_7)}\cr &\qquad\quad\text{\rm if}\quad a_7+a_8+a_9\geq
a_3+a_5+a_6,\quad a_1+a_2\geq a_5+a_6,\cr &\qquad\qquad\qquad
a_3\geq a_7,\quad a_5\geq a_8,\quad  a_{10}\geq a_8,\quad a_6\geq
a_2.\cr&\cr \endalign $$$$\align &(2).\;\;\sum_{0\leq u\leq
a_2}\!\!(-1)^u \left[\smallmatrix a_5-a_1-1+u\cr u
\endsmallmatrix
\right]e_3^{(a_2-u)}e_2^{(a_3+a_6)}e_4^{(a_1+a_2)}e_3^{(a_3+a_5+a_6+u)}\cr
&\qquad\qquad\qquad\times
e_2^{(a_8)}e_1^{(a_4+a_7+a_9+a_{10})}e_2^{(a_4+a_7+a_9)}e_3^{(a_4)}e_4^{(a_3+a_4)}
e_3^{(a_7)}\cr &\qquad\qquad\qquad\text{\rm if}\quad
a_7+a_8+a_9\geq a_3+a_5+a_6,\quad a_5\geq a_1\geq a_8,\cr
&\qquad\qquad\qquad \qquad\qquad a_3\geq a_7,\quad  a_{10}\geq
a_8,\quad a_6\geq a_2.\cr&\cr {\bold {25.}}\quad
&(1).\!\!\!\sum_{0\leq u\leq a_4+a_7}\!\!\!\!\!\!\!(-1)^u\!
\left[\!\smallmatrix a_9+a_{10}-a_6-a_8-1+u\cr u
\endsmallmatrix
\!\right]\!e_2^{(a_3)}e_3^{(a_2+a_3)}e_1^{(a_4+a_7-u)}e_2^{(a_4+a_6+a_7)}\cr
&\qquad\qquad\times
e_3^{(a_4)}e_4^{(a_1+a_2+a_3+a_4)}e_3^{(a_5+a_6+a_7)}e_2^{(a_8)}e_1^{(a_9+a_{10}+u)}
e_2^{(a_9)}\cr &\qquad\quad\text{\rm if}\quad a_1+a_2+a_3\geq
a_5+a_6+a_7,\quad a_9+a_{10}\geq a_6+a_8,\cr &\qquad\qquad \qquad
a_6+a_7\geq a_2+a_3,\quad a_2\geq  a_6\geq a_9,\quad a_5\geq
a_8.\cr &(2).\;\;\sum_{0\leq u\leq a_3}\!\!(-1)^u
\left[\smallmatrix a_6-a_2-1+u\cr u
\endsmallmatrix
\right]e_2^{(a_3-u)}e_3^{(a_2+a_3)}e_1^{(a_4+a_7)}e_2^{(a_4+a_6+a_7+u)}\cr
&\qquad\qquad\times
e_3^{(a_4)}e_4^{(a_1+a_2+a_3+a_4)}e_3^{(a_5+a_6+a_7)}e_2^{(a_8)}e_1^{(a_9+a_{10})}
e_2^{(a_9)}\cr &\qquad\quad\text{\rm if}\quad a_1+a_2+a_3\geq
a_5+a_6+a_7,\quad a_2+a_8\geq a_9+a_{10},\cr &\qquad\qquad \qquad
a_7\geq a_3,\quad a_6\geq  a_2,\quad a_{10}\geq a_8,\quad a_5\geq
a_8.\cr&\cr  {\bold {26.}}\quad &(1).\!\!\!\sum_{0\leq u\leq
a_9}\!\!(-1)^u\left[\!\smallmatrix a_8-a_{10}-1+u\cr u
\endsmallmatrix
\!\right]\!e_1^{(a_4)}e_2^{(a_3+a_4)}e_3^{(a_2+a_3+a_4)}e_4^{(a_1+a_2+a_3+a_4)}\cr
&\qquad\qquad\times
e_3^{(a_5)}e_1^{(a_7)}e_2^{(a_6+a_7+a_8+u)}e_1^{(a_9+a_{10})}e_3^{(a_6+a_7)}e_2^{(a_9-u)}\cr
&\qquad\qquad\qquad\text{\rm if}\quad a_2+a_5\geq a_6+a_8,\quad
a_8\geq a_{10}\geq a_5,\cr &\qquad\qquad\qquad \qquad  a_1\geq
a_5,\quad a_3\geq a_7,\quad a_6\geq a_9.\cr&\cr &(2).\sum_{0\leq
u\leq a_3+a_4}\!\!(-1)^u\left[\smallmatrix a_6+a_8-a_2-a_5-1+u\cr
u
\endsmallmatrix
\right]\!e_1^{(a_4)}e_2^{(a_3+a_4-u)}e_3^{(a_2+a_3+a_4)}\cr
&\qquad\times
e_4^{(a_1+a_2+a_3+a_4)}e_3^{(a_5)}e_1^{(a_7)}e_2^{(a_6+a_7+a_8+u)}e_1^{(a_9+a_{10})}
e_3^{(a_6+a_7)}e_2^{(a_9)}\cr &\qquad\quad\text{\rm if}\quad
a_6+a_8\geq a_2+a_5\geq a_9+a_{10},\quad a_1\geq a_5,\cr
&\qquad\qquad\qquad a_{10}\geq a_8,\quad a_3\geq a_7,\quad a_2\geq
a_6.\cr&\cr {\bold {27.}}\quad &(1).\;\;\sum_{0\leq u\leq
a_3+a_4+a_6}\!\!(-1)^u\left[\smallmatrix a_5-a_8-1+u\cr u
\endsmallmatrix
\right]e_4^{(a_1)}e_3^{(a_2+a_5+u)}e_2^{(a_3+a_6+a_8)}\cr &\qquad
\times
e_1^{(a_4+a_7+a_9+a_{10})}e_2^{(a_4)}e_4^{(a_2)}e_3^{(a_3+a_4+a_6-u)}e_4^{(a_3+a_4)}
e_2^{(a_7+a_9)}e_3^{(a_7)}\cr
&\qquad\text{\rm if}\quad a_7+a_9+a_{10}\geq
a_3+a_6+a_8,\quad a_3+a_6\geq a_7+a_9,\cr &\quad\qquad\qquad
 a_8+a_9\geq a_5+a_6,\quad a_5\geq a_8\geq a_1,\quad a_6\geq
a_2.\cr&\cr \endalign $$$$\align &(2).\!\!\!\!\!\sum_{0\leq u\leq
a_7}\!\!(-1)^u \left[\smallmatrix a_6-a_9-1+u\cr u
\endsmallmatrix
\right]e_4^{(a_1)}e_3^{(a_2+a_5)}e_2^{(a_3+a_6+a_8)}e_1^{(a_4+a_7+a_9+a_{10})}\cr
&\qquad\qquad\times
e_2^{(a_4)}e_4^{(a_2)}e_3^{(a_3+a_4+a_6+u)}e_4^{(a_3+a_4)}e_2^{(a_7+a_9)}e_3^{(a_7-u)}\cr
&\qquad\qquad\qquad\text{\rm if}\quad a_7+a_9+a_{10}\geq
a_3+a_6+a_8,\quad a_3\geq a_7,\cr &\qquad\qquad\qquad\qquad
\qquad a_6\geq a_9\geq a_2,\quad a_8\geq  a_5\geq a_1.\cr&\cr
 {\bold {28.}}\quad&(1).\sum_{0\leq u \leq
a_4+a_7+a_9}\!\!\!\!\!\!(-1)^u\left[\smallmatrix a_8-a_{10}-1+u\cr
u
\endsmallmatrix
\right]e_3^{(a_2)}e_4^{(a_1+a_2)}e_3^{(a_5)}e_2^{(a_3+a_6+a_8+u)}\cr
&\quad\quad\quad\times
e_3^{(a_3+a_6)}e_1^{(a_4+a_7+a_9+a_{10})}e_2^{(a_4+a_7+a_9-u)}e_3^{(a_4)}e_4^{(a_3+a_4)}
e_3^{(a_7)}\cr &\qquad\qquad\text{\rm if}\quad a_5+a_6\geq
a_1+a_2,\quad a_8\geq a_{10}\geq a_5,\quad a_1\geq a_5,\cr
&\qquad\qquad \qquad a_7+a_9\geq a_3+a_6,\quad a_3\geq  a_7.\cr
&(2).\sum_{0\leq u\leq a_3+a_4}\!\!\!\!\!(-1)^u \left[\smallmatrix
a_1+a_2-a_5-a_6-1+u\cr u
\endsmallmatrix
\right]e_3^{(a_2)}e_4^{(a_1+a_2+u)}e_3^{(a_5)}e_2^{(a_3+a_6+a_8)}\cr
&\qquad\qquad\times
e_3^{(a_3+a_6)}e_1^{(a_4+a_7+a_9+a_{10})}e_2^{(a_4+a_7+a_9)}e_3^{(a_4)}e_4^{(a_3+a_4-u)}
e_3^{(a_7)}\cr &\qquad\qquad\text{\rm if}\quad a_7+a_9\geq
a_3+a_6,\quad a_{10}\geq a_8\geq a_5,\quad a_6\geq a_2,\cr
&\qquad\qquad\qquad \qquad a_1+a_2\geq a_5+a_6,\quad a_3\geq
a_7.\cr &\cr {\bold {29.}}\quad &(1).\;\;\sum_{0\leq u\leq
a_3+a_6}\!\!(-1)^u \left[\smallmatrix a_8-a_5-1+u\cr u
\endsmallmatrix
\right]e_3^{(a_2)}e_2^{(a_3+a_6-u)}e_1^{(a_4+a_7+a_9)}e_3^{(a_3)}\cr
&\qquad\qquad\times
e_4^{(a_1+a_2+a_3)}e_3^{(a_5+a_6)}e_2^{(a_4+a_7+a_8+a_9+u)}e_3^{(a_4+a_7)}e_4^{(a_4)}
e_1^{(a_{10})}\cr &\qquad\quad\text{\rm if}\quad a_5+a_6+a_7\geq
a_1+a_2+a_3,\quad a_1+a_2\geq a_5+a_6,\cr &\qquad\qquad\qquad
\quad a_8\geq a_5\geq a_{10},\quad a_9\geq a_6\geq a_2.\cr&\cr
{\bold {30.}}\quad &(1).\;\;\sum_{0\leq u\leq
a_2+a_3}\!\!(-1)^u\left[\smallmatrix a_5-a_1-1+u\cr u
\endsmallmatrix
\right]e_2^{(a_3)}e_3^{(a_2+a_3-u)}e_4^{(a_1+a_2+a_3)}e_2^{(a_6)}\cr
&\qquad\qquad\times
e_1^{(a_4+a_7+a_9)}e_2^{(a_4+a_7)}e_3^{(a_4+a_5+a_6+a_7+u)}e_2^{(a_8+a_9)}e_1^{(a_{10})}
e_4^{(a_4)}\cr &\qquad\quad\text{\rm if}\quad a_6+a_7\geq
a_2+a_3,\quad a_1+a_6\geq a_8+a_9,\quad a_2\geq a_6,\cr
&\qquad\qquad\qquad\qquad
 a_5\geq a_1,\quad a_9\geq a_6,\quad a_8\geq a_{10}.\cr&\cr
 {\bold {31.}}\quad &(1).\!\!\!\sum_{0\leq
u\leq a_3+a_4+a_6}\!\!\!\!(-1)^u\left[\smallmatrix a_8-a_5-1+u\cr
u
\endsmallmatrix
\right]e_1^{(a_4)}e_3^{(a_2)}e_2^{(a_3+a_4+a_6-u)}e_1^{(a_7+a_9)}\cr
&\qquad\qquad\times
e_3^{(a_3+a_4)}e_4^{(a_1+a_2+a_3+a_4)}e_3^{(a_5+a_6)}e_2^{(a_7+a_8+a_9+u)}e_1^{(a_{10})}
e_3^{(a_7)}\cr &\qquad\qquad\quad\text{\rm if}\qquad a_1+a_2\geq
a_5+a_6,\quad a_3+a_6\geq a_7+a_9,\cr &\qquad\qquad\qquad\quad
\qquad a_9\geq a_6\geq a_2,\quad a_8\geq a_5\geq a_{10}.
\endalign $$
\par {\rm 2.}\; Applying the map $\Gamma \circ\Psi\circ\Phi$ defined
in {\rm [2]} to cases {\bf 1}$-${\bf 31}, we get another $72$
polynomial elements in one variable in the canonical basis $\bold
B$.
\endproclaim

{\head{2. Proof of Theorem 1.3 }\endhead}

{\bf 2.1.} In order to prove Theorem 1.3. We firstly need the
following identity showed in [9]. Assume that $m\geq k\geq
0,\;\delta\in\bN.$ Then
$$\sum_{0\leq i\leq\delta}(-1)^i\left[\matrix k-1+i\\i
\endmatrix\right]\left[\matrix m\\\delta-i\endmatrix\right]v^{i(m-k)}=
\left[\matrix m-k\\\delta\endmatrix\right]v^{-k\delta}.\tag i$$
\par
Also, we need the following identity showed in [10]. Assume
that\newline $m\geq k\geq 0,\;\delta,n\in\bN.$ Then
$$\align
&\sum_{0\leq i\leq\delta}(-1)^i\left[\matrix k-1+i\\i
\endmatrix\right]\left[\matrix m+n\\\delta-i
\endmatrix\right]v^{i(m-k-n)}\tag ii\cr&\cr
& =\sum_{0\leq t\leq\delta,n}\left[\matrix m-k\\\delta-t\endmatrix
\right]\left[\matrix n\\t\endmatrix\right]v^{-k (\delta-t)-n\delta+t(m+n)}.
\endalign $$\par \vskip0.1cm
 {\bf 2.2.} Now we prove {\rm 1} of Theorem 1.3. It is obvious that
all the elements from case ${\bold 1}$ to case ${\bold {31}}$ in
Theorem 1.3 are fixed by the involution $\bar{\cdot}$. So we only
need to check that these elements lie in $\Cal L$. Moreover, just
like what we have done in [2] that we only prove the most
complicated case here, i.e., case ${\bold 1}$. First of all, we
observe the monomial corresponding to $(1)-(3)$ in case ${\bold
1}$. Using the commutative relations in [2, \S 1.3], we have
$$\align & e_{2}^{(a_3)}
e_{3}^{(a_2+a_3)}e_{4}^{(a_1+a_2+a_3)}e_{2}^{(a_6)}e_{3}^{(a_5+a_6)}e_{2}^{(a_8)}
e_{1}^{(a_4+a_7+a_9+a_{10})}\tag iii\cr&\cr &\quad\times
e_{2}^{(a_4+a_7+a_9)}e_{3}^{(a_4+a_7)}e_{4}^{(a_4)}
=\sum_{\omega\,\in\,\Omega}v^{{\Cal A }(\omega)}\times {\Cal
B}(\omega)\times E^{\omega},
\endalign$$
where
$$\align {\Cal A}(\omega): =
&-(a_4+a_7+a_9+a_{10}-i)(a_4+a_7+a_9-i)-(a_4+a_7-j)(i-j)\cr&\cr
&-(a_1+a_2+a_3-r)(k-r)-(a_3-k+a_6-l)(a_5+a_6-l)\cr&\cr
&-(a_3-k+a_6-l+a_8+a_4+a_7+a_9-i-m)(a_4+a_7-j-m)\cr&\cr
&-(a_3-k)(a_2+a_3-k)-(a_2+a_3-k-t)(a_1+a_2+a_3-r-t)\cr&\cr
\endalign $$$$\align
&-(a_2\!+a_3\!-\!k-\!t+\!a_5\!+a_6\!-l+\!a_4\!+a_7\!-j\!-m-\!p)
(a_4\!-n\!-s-\!p)\cr&\cr
&-(a_4-n)(j-n)-(a_4-n-s)(l+m+k-r-s)\cr&\cr &+(a_4+a_7-j-m)l+
(a_5+a_6-l+a_4+a_7-j-m)(k-r)\cr&\cr &+(a_4-n-s-p)(r+t)+pr,
\endalign$$
and
$$\align {\Cal B}(\omega): &=\left[\matrix a_8+a_4+a_7+a_9-i\\ a_8
\endmatrix\right]\left[\matrix a_3-k+a_6-l+a_8+a_4+a_7+a_9-i\\a_3-k+a_6-l
\endmatrix\right]\cr&\cr &\times\left[\matrix a_3-k+a_6\\a_6
\endmatrix\right]\left[\matrix l+m\\l
\endmatrix\right]\left[\matrix l+m+k-r\\k-r
\endmatrix\right]\left[\matrix p+t\\p
\endmatrix\right]\left[\matrix r+s\\r
\endmatrix\right]
\cr&\cr&\times\left[\matrix a_1+a_2+a_3-r-t+a_4-n-s-p\\a_4-n-s-p
\endmatrix\right]\cr&\cr & \times \left[\matrix a_5+a_6-l+a_4+a_7-j-m\\a_4+a_7-j-m
\endmatrix\right]\cr&\cr&\times\left[\matrix a_2+a_3-k-t+a_5+a_6-l+a_4+a_7-j-m\\a_2+a_3-k-t
\endmatrix\right].\endalign $$
Note that the last two factors in ${\Cal B}(\omega)$ can also be
equivalently represented as follows
$$\align &\left[\matrix a_2+a_3-k-t+a_5+a_6-l\\a_2+a_3-k-t
\endmatrix\right]\cr&\cr&\times\left[\matrix a_2+a_3-k-t+a_5+a_6-l+a_4+a_7-j-m\\a_4+a_7-j-m
\endmatrix\right].
\endalign$$
Moreover, we have
$$\align E^{\omega}: &=e_4^{(a_1+a_2+a_3-r-t+a_4-n-s-p)}e_{34}^{(p+t)}
e_{24}^{(r+s)}e_{14}^{(n)}\cr&\cr&\times
e_3^{(a_2+a_3-k-t+a_5+a_6-l+a_4+a_7-j-m-p)}
e_{23}^{(l+m+k-r-s)}e_{13}^{(j-n)}\cr&\cr&\times
e_{2}^{(a_3-k+a_6-l+a_8+a_4+a_7+a_9-i-m)}
e_{12}^{(i-j)}e_{1}^{(a_4+a_7+a_9+a_{10}-i)},\endalign$$
and
$$\Omega=\{\omega = (i,j,k,l,m,n,r,s,t,p)\}\subset {\bN}^{10},$$
where integers $\,i,j,k,l,m,n,r,s,t,p\,$ satisfies the following
inequalities:
$$\align
&\quad 0\leq i\leq a_4+a_7+a_9;\quad  0\leq j\leq a_4+a_7, i;\quad
0\leq k \leq a_3;\cr&\cr &\quad 0\leq m\leq
a_4+a_7-j,a_3-k+a_6-l+a_8+a_4+a_7+a_9-i;\cr&\cr (\flat) &\quad
0\leq l\leq a_3-k+a_6, a_5+a_6;\quad 0\leq n\leq a_4,j;\quad 0\leq
r\leq k;\cr&\cr &\quad 0\leq s\leq a_4-n,l+m+k-r;\quad 0\leq t\leq
a_2+a_3-k;\cr&\cr &\quad 0\leq p \leq
a_4-n-s,a_2+a_3-k-t+a_5+a_6-l+a_4+a_7-j-m.
\endalign
$$
\par
It should be mentioned here that in the above argument we use $10$
different parameters to describe $\Omega$ in order to simplify the
proof of polynomial elements {\bf 1.}$(1), (2), (3)$ in one
variable, however we use $14$ different parameters to describe
$\Omega$ in the proof of the monomial element {\bf 1.}$(1)$ in
[2]. It is only because the commutative order we use here is
somewhat different from that in [2]. \par Set
$$\align
&\qquad x_1=a_4+a_7+a_9-i,\qquad
x_2=a_4+a_7-j-m,\qquad\qquad\qquad\cr (\natural)&\qquad
x_3=a_3-k,\qquad\qquad\qquad\;x_4=a_3-k+a_6-l,\qquad\qquad\qquad\cr
&\qquad x_5=a_4-n-s-p,\qquad \;\;\; x_6=k-r,\qquad\qquad\qquad\cr
&\qquad x_7=a_2+a_3-k-t.
\endalign$$
Then the degree (with respect to $v$) of the coefficient $v^{{\Cal
A}(\omega)}\times {\Cal B}(\omega)$ of $E^{\omega}$ in the sum
expression of Formula (iii) is
$$D_M: =-L_M(x_1,x_2,\cdots,x_7,m,s,p)-Q_{M}(x_1,x_2,\cdots,x_7,m,s,p),$$
where $L_M(x_1,x_2,\cdots,x_7,m,s,p)$ is a linear form in
non-negative integers $x_1,x_2,\cdots,x_7,m,s,p$, and $Q_{M}(x_1,x_2,\cdots,x_7,m,s,p)$
is a unit form in non-negative integers $x_1,$ $x_2,\cdots,x_7,m,s,p$. Moreover, we have
$$\align
&L_M(x_1,x_2,\cdots,x_7,m,s,p): \cr&\cr
&=(a_{10}-a_8)x_1+(a_8+a_9-a_5-a_6)x_2+(a_2-a_6)x_3+(a_5-a_8)x_4\cr&\cr
&+(a_5+a_6+a_7-a_1-a_2-a_3)x_5+(a_1+a_2-a_5-a_6)x_6\cr&\cr
&+(a_1-a_5)x_7+(a_9-a_6)m+(a_7-a_3)s+(a_6+a_7-a_2-a_3)p,
\endalign$$
and
$$\align
& Q_{M}(x_1,x_2,\cdots,x_7,m,s,p): \cr&\cr
&=x_1^2+x_2^2+x_3^2+x_4^2+x_5^2+x_6^2+ x_7^2+m^2+ s^2+p^2\cr&\cr
&+x_2x_4+x_2m+x_3x_6 +x_3x_7+x_3s+x_3p+x_4m+x_5x_6+x_5x_7\cr&\cr
&+x_5s+x_5p+x_6x_7+x_6s+2x_6p+x_7p+sp-x_1x_4-x_1m-x_2x_6\cr&\cr
&-x_2x_7-x_2s-x_2p-x_3x_4-x_3m-x_4x_7-x_4p-x_6m-ms.
\endalign$$
\par
By the BDDP-algorithm (see [1] or [2] \S 3.3) and Lemma 3.4 in
[2], we know that the unit form $Q_{M}(x_1,x_2,\cdots,x_7,m,s,p)$
is weakly positive, i.e.,
$$Q_{M}(x_1,x_2,\cdots,x_7,m,s,p)\geq 0\quad \text {for any}\quad (x_1,x_2,\cdots,x_7,m,s,p)
\in {\bN}^{10}.$$ Therefore, we have
$$\align
&\sum_{0\leq u\leq a_4+a_7}(-1)^u\left[\smallmatrix a_5+a_6-a_8-a_9-1+u\\
u\endsmallmatrix\right]e_2^{(a_3)}e_3^{(a_2+a_3)}e_4^{(a_1+a_2+a_3)}e_2^{(a_6)}\cr&\cr
&\qquad\;\;\times
e_3^{(a_5+a_6+u)}e_2^{(a_8)}e_1^{(a_4+a_7+a_9+a_{10})}e_2^{(a_4+a_7+a_9)}
e_3^{(a_4+a_7-u)}e_4^{(a_4)}\cr &\cr
&=\sum_{\omega \in
\Omega_1}\left(\sum_{0\leq u\leq a_4+a_7-j-m} (-1)^u\left[\matrix
a_5+a_6-a_8-a_9-1+u
\\u\endmatrix\right]\right.\cr &\cr
&\qquad \times \left.\left[\matrix a_4+a_5+a_6+a_7-j-m-l\\
a_4+a_7-j-m-u\endmatrix\right]\times
v^{(a_4+a_7+a_8+a_9-j-m-l)u}\right)\tag iv\cr &\cr &\qquad \times
v^{{\Cal A}(\omega)}\times {\Cal B}_1(\omega)\times
E^{\omega}\cr\cr &=\sum_{\omega\in\Omega_1}v^{{\Cal
A}(\omega)-(a_5+a_6-a_8-a_9)(a_4+a_7-j-m)}\times {\Cal
B}_1(\omega)\cr&\cr&\qquad
\times\left[\matrix a_4+a_7+a_8+a_9-j-m-l\\
a_4+a_7-j-m\endmatrix\right]\times E^{\omega},\endalign$$ where
the last equality comes from Formula 2.1 (i), $\Omega_1$ is
obtained from $\Omega$ by replacing the defining inequality
$``0\leq l\leq a_3-k+a_6,a_5+a_6"$ in $(\flat)$ by $``0\leq l\leq
a_3-k+a_6,a_5+a_6+u"$, and ${\Cal B}_1(\omega)$ is obtained from
$\Cal B(\omega)$ by deleting the factor $\left[\matrix
a_5+a_6-l+a_4+a_7-j-m\\a_4+a_7-j-m\endmatrix \right].$ Note that
$a_5+a_6$ is always less than or equal to $a_5+a_6+u$, and
relations $(\natural)$ are independent of $u$, even if $u$ occurs
in the upper boundary of the defining inequality of $\,l$, we can
still conclude that the expressions behind the second equal sign
in (iv) are independent of $u$.
\par
Using relations $(\natural)$, we can get the degree $D_{P_1}$
(with respect to $v$) of the coefficient of $E^{\omega}$ in the
last sum expression of (iv), i.e.,
$$\align
D_{P_1}&=D_M-(a_5+x_4-x_3)x_2+(a_8+a_9-a_6+x_4-x_3)x_2\cr&\cr
&\quad-(a_5+a_6-a_8-a_9)x_2\cr&\cr
&=-L_{P_1}(x_1,x_2,\cdots,x_7,m,s,p)-Q_{P_1}(x_1,x_2,\cdots,x_7,m,s,p),
\endalign$$
where
$$\align
&\;\; L_{P_1}(x_1,x_2,\cdots,x_7,m,s,p)\cr&\cr
=&(a_{10}-a_8)x_1+(a_5+a_6-a_8-a_9)x_2+(a_2-a_6)x_3+(a_5-a_8)x_4\cr&\cr
+&(a_5+a_6+a_7-a_1-a_2-a_3)x_5+(a_1+a_2-a_5-a_6)x_6\cr&\cr
+&(a_1-a_5)x_7+(a_9-a_6)m+(a_7-a_3)s+(a_6+a_7-a_2-a_3)p,
\endalign$$
and
$$Q_{P_1}(x_1,x_2,\cdots,x_7,m,s,p)=Q_{M}(x_1,x_2,\cdots,x_7,m,s,p).$$
When
$$\align
&a_5+a_6+a_7\geq a_1+a_2+a_3,\quad a_5+a_6\geq a_8+a_9,\cr&\cr
&\qquad a_1\geq a_5,\quad a_{10}\geq a_8,\quad a_2\geq a_6,\quad
a_9\geq a_6,
\endalign$$
we have $D_{P_1}\leq 0$. Moreover, $D_{P_1}=0\Leftrightarrow
x_1=\cdots =x_7=m=s=p=0,$ and $E^{\omega}=E^A$. Therefore, we have
$$\align
\quad &\sum_{0\leq u\leq a_4+a_7}(-1)^u\left[\smallmatrix a_5+a_6-a_8-a_9-1+u\\
u\endsmallmatrix\right]
e_2^{(a_3)}e_3^{(a_2+a_3)}e_4^{(a_1+a_2+a_3)}e_2^{(a_6)}\cr
&\qquad\;\;\times
e_3^{(a_5+a_6+u)}e_2^{(a_8)}e_1^{(a_4+a_7+a_9+a_{10})}e_2^{(a_4+a_7+a_9)}
e_3^{(a_4+a_7-u)}e_4^{(a_4)}\cr&\cr &\equiv
e_{4}^{(a_1)}e_{34}^{(a_2)}e_{24}^{(a_3)}e_{14}^{(a_4)}e_{3}^{(a_5)}e_{23}^{(a_6)}
e_{13}^{(a_7)}e_{2}^{(a_8)}e_{12}^{(a_9)}e_{1}^{(a_{10})}\qquad
(\mod v^{-1}{\Cal L})\cr &\cr & \quad \text{\rm if}\quad
a_5+a_6+a_7\geq a_1+a_2+a_3,\quad a_5+a_6\geq a_8+a_9,\cr &\qquad
\qquad a_1\geq a_5,\quad a_{10}\geq a_8,\quad a_2\geq a_6,\quad
a_9\geq a_6.
\endalign$$
So we have proved $(1)$ of case $\bold 1$. Let us consider
$$\align
&\sum_{0\leq u\leq a_4}(-1)^u\left[\smallmatrix a_1+a_2+a_3-a_5-a_6-a_7-1+u\\
u\endsmallmatrix\right]
e_2^{(a_3)}e_3^{(a_2+a_3)}e_4^{(a_1+a_2+a_3+u)}\cr
&\qquad\;\;\times
e_2^{(a_6)}e_3^{(a_5+a_6)}e_2^{(a_8)}e_1^{(a_4+a_7+a_9+a_{10})}
e_2^{(a_4+a_7+a_9)}e_3^{(a_4+a_7)}e_4^{(a_4-u)}\cr&\cr &
=\sum_{\omega\in\Omega}\left(\sum_{0\leq u\leq
a_4-n-s-p}(-1)^u\left[\matrix a_1+a_2+a_3-
a_5-a_6-a_7-1+u\\u\endmatrix\right]\right.\cr&\cr
&\qquad\qquad\qquad\times \left[\matrix
a_1+a_2+a_3+a_4-r-t-n-s-p\\a_4-n-s-p-u\endmatrix\right]\cr&\cr
&\left.\qquad\qquad\qquad \times
v^{(a_4+a_5+a_6+a_7-r-t-n-s-p)u}\right) \times v^{{\Cal
A}(\omega)}\times {\Cal B}_2(\omega)\times E^{\omega}\tag v\cr&\cr
&=\sum_{\omega\in\Omega}v^{{\Cal
A}(\omega)-(a_1+a_2+a_3-a_5-a_6-a_7)(a_4-n-s-p)}\times {\Cal
B}_2(\omega)\cr&\cr& \qquad\qquad\times\left[\matrix
a_4+a_5+a_6+a_7-r-t-n-s-p\\a_4-n-s-p\endmatrix\right]\times
E^{\omega},
\endalign$$
where the last equality comes from Formula 2.1 (i), ${\Cal
B}_2(\omega)$ is obtained from $\Cal B(\omega)$ by delating the
factor $\left[\matrix
a_1+a_2+a_3-r-t+a_4-n-s-p\\a_4-n-s-p\endmatrix\right].$
\par
Using relations $(\natural)$, we can get the degree $D_{P_2}$
(with respect to $v$) of the coefficient of $E^{\omega}$ in the
last sum expression of (v), i.e.,
$$\align
D_{P_2}&=D_M-(a_1+x_6+x_7)x_5+(a_5+a_6+a_7-a_2-a_3+x_6+x_7)x_5\cr&\cr
&\quad-(a_1+a_2+a_3-a_5-a_6-a_7)x_5\cr&\cr
&=-L_{P_2}(x_1,x_2,\cdots,x_7,m,s,p)-Q_{P_2}(x_1,x_2,\cdots,x_7,m,s,p),
\endalign$$
where
$$\align
&\;\; L_{P_2}(x_1,x_2,\cdots,x_7,m,s,p)\cr&\cr
=&(a_{10}-a_8)x_1+(a_8+a_9-a_5-a_6)x_2+(a_2-a_6)x_3+(a_5-a_8)x_4\cr&\cr
+&(a_1+a_2+a_3-a_5-a_6-a_7)x_5+(a_1+a_2-a_5-a_6)x_6\cr&\cr
+&(a_1-a_5)x_7+(a_9-a_6)m+(a_7-a_3)s+(a_6+a_7-a_2-a_3)p,
\endalign$$
and
$$Q_{P_2}(x_1,x_2,\cdots,x_7,m,s,p)=Q_{M}(x_1,x_2,\cdots,x_7,m,s,p).$$
When
$$\align
&a_1+a_2+a_3\geq a_5+a_6+a_7,\quad a_8+a_9\geq a_5+a_6,\quad
a_2\geq a_6,\cr &\qquad a_6+a_7\geq a_2+a_3,\quad a_{10}\geq
a_8,\quad a_5\geq a_8,
\endalign$$
we have $D_{P_2}\leq 0$. Moreover, $D_{P_2}=0\Leftrightarrow
x_1=\cdots =x_7=m=s=p=0,$ and $E^{\omega}=E^A$. Therefore, we have
$$\align
\quad &\sum_{0\leq u\leq a_4}(-1)^u\left[\smallmatrix a_1+a_2+a_3-a_5-a_6-a_7-1+u\\
u\endsmallmatrix\right]
e_2^{(a_3)}e_3^{(a_2+a_3)}e_4^{(a_1+a_2+a_3+u)}e_2^{(a_6)}\cr
&\qquad\;\;\times
e_3^{(a_5+a_6)}e_2^{(a_8)}e_1^{(a_4+a_7+a_9+a_{10})}e_2^{(a_4+a_7+a_9)}
e_3^{(a_4+a_7-u)}e_4^{(a_4)}\cr&\cr &\equiv
e_{4}^{(a_1)}e_{34}^{(a_2)}e_{24}^{(a_3)}e_{14}^{(a_4)}e_{3}^{(a_5)}e_{23}^{(a_6)}
e_{13}^{(a_7)}e_{2}^{(a_8)}e_{12}^{(a_9)}e_{1}^{(a_{10})}\qquad
(\mod v^{-1}{\Cal L})\cr &\cr & \quad \text{\rm if}\quad
a_1+a_2+a_3\geq a_5+a_6+a_7,\quad a_8+a_9\geq a_5+a_6,\quad
a_2\geq a_6,\cr &\qquad \qquad \quad a_6+a_7\geq a_2+a_3,\quad
a_{10}\geq a_8,\quad a_5\geq a_8.
 \endalign$$ This gives (2) of case $\bold 1$. Finally, we
 consider
$$\align
&\sum_{0\leq u\leq a_2+a_3}(-1)^u\left[\smallmatrix
a_5-a_1-1+u\\u\endsmallmatrix \right]
e_2^{(a_3)}e_3^{(a_2+a_3-u)}e_4^{(a_1+a_2+a_3)}e_2^{(a_6)}\cr
&\qquad\;\times
e_3^{(a_5+a_6+u)}e_2^{(a_8)}e_1^{(a_4+a_7+a_9+a_{10})}e_2^{(a_4+a_7+a_9)}
e_3^{(a_4+a_7)}e_4^{(a_4)}\cr&\cr &
=\sum_{\omega\in\Omega_1}\left(\sum_{0\leq u\leq a_2+a_3-k-t}
(-1)^u\left[\matrix a_5-a_1-1+u\\ u\endmatrix\right]\right.\cr&\cr
&\qquad\qquad\qquad\times\left[\matrix a_2+a_3+a_5-r-t+a_6-l-k+r\\
a_2+a_3-k-t-u
\endmatrix\right]\cr&\cr
&\left.\qquad\qquad\qquad\times
v^{(a_1+a_2+a_3-a_6+l+k-2r-t)u}\right)\times v^{{\Cal A}(\omega)}
\times {\Cal B}_3(\omega)\times E^{\omega}\tag vi \cr&\cr
&=\sum_{\omega\in\Omega_2}v^{{\Cal
A}(\omega)-(a_5-a_1)(a_2+a_3-k-t-w)-(a_6-l-k+r)(a_2+a_3-k-t)}
\cr&\cr& \qquad\times v^{(a_2+a_3-k-t+a_5+a_6-l)w}\times {\Cal
B}_3(\omega)\times\left[\matrix
a_1+a_2+a_3-r-t\\a_2+a_3-k-t-w\endmatrix\right]\cr&\cr
&\qquad\times\left[\matrix a_6-l-k+r\\w\endmatrix\right]\times
E^{\omega},
\endalign$$
where the last equality comes from the Formula 2.1 (ii),
$\Omega_2$ is obtained from $\Omega_1$ by adding $``0\leq w\leq
a_2+a_3-k-t, a_6-l-k+r"$ to its defining inequalities only, and
${\Cal B}_3(\omega)$ is obtained from $\Cal B(\omega)$ by delating
the factor $\left[\matrix
a_2+a_3-k-t+a_5+a_6-l\\a_2+a_3-k-t\endmatrix\right].$ Note that
the expressions behind the second equal sign in (vi) are also
independent of $u$, just like (1) of case {\bf 1}.
\par
Using relations $(\natural)$, we can get the degree $D_{P_3}$
(with respect to $v$) of the coefficient of $E^{\omega}$ in the
last sum expression of (vi), i.e.,
$$\align
D_{P_3}&=D_M-(a_5+x_4-x_3)x_7-(a_5-a_1)(x_7-w)-(x_4-x_3-x_6)x_7\cr&\cr
&\quad+(a_5+x_7+x_4-x_3)w+(a_1+x_6+w)(x_7-w)\cr&\cr&\quad
+(x_4-x_3-x_6-w)w\cr\endalign $$$$\align
&=-L_{P_3}(x_1,x_2,\cdots,x_7,m,s,p,w)-Q_{P_3}(x_1,x_2,\cdots,x_7,m,s,p,w),
\endalign$$
where
$$\align
&L_{P_3}(x_1,x_2,\cdots,x_7,m,s,p,w)\cr&\cr =&
(a_{10}-a_8)x_1+(a_8+a_9-a_5-a_6)x_2+(a_2+a_5-a_1-a_6)x_3\cr&\cr
&+(a_1-a_8)x_4+(a_5+a_6+a_7-a_1-a_2-a_3)x_5+(a_2-a_6)x_6\cr&\cr
&+(a_5-a_1)(x_7-w)+(a_9-a_6)m+(a_7-a_3)s\cr&\cr
&+(a_6+a_7-a_2-a_3)p+(a_5-a_1)(x_4-x_3-x_6-w),
\endalign$$
and
$$\align
&Q_{P_3}(x_1,x_2,\cdots,x_7,m,s,p,w)\cr&\cr =&
Q_{M}(x_1,x_2,\cdots,x_7,m,s,p)+2(x_4-x_3-x_6-w)(x_7-w).
\endalign$$
By the definition of $\Omega_2$, we have $(x_4-x_3-x_6-w)\geq 0,$
and $(x_7-w)\geq 0$. When
$$\align
&a_8+a_9\geq a_5+a_6,\quad a_6+a_7\geq a_2+a_3,\quad a_2\geq a_6,\cr
&\qquad a_5\geq a_1\geq a_8,\quad a_{10}\geq a_8,
\endalign$$
we have $D_{P_3}\leq 0$. Moreover, $D_{P_3}=0\Leftrightarrow
x_1=\cdots =x_7=m=s=p=w=0,$ and $E^{\omega}=E^A.$ Then we have
$$\align
\quad &\sum_{0\leq u\leq a_2+a_3}(-1)^u\left[\smallmatrix
a_5-a_1-1+u\\u\endsmallmatrix\right]
e_2^{(a_3)}e_3^{(a_2+a_3-u)}e_4^{(a_1+a_2+a_3)}e_2^{(a_6)}\cr
&\qquad\;\;\times
e_3^{(a_5+a_6+u)}e_2^{(a_8)}e_1^{(a_4+a_7+a_9+a_{10})}e_2^{(a_4+a_7+a_9)}
e_3^{(a_4+a_7)}e_4^{(a_4)}\cr&\cr &\equiv
e_{4}^{(a_1)}e_{34}^{(a_2)}e_{24}^{(a_3)}e_{14}^{(a_4)}e_{3}^{(a_5)}e_{23}^{(a_6)}
e_{13}^{(a_7)}e_{2}^{(a_8)}e_{12}^{(a_9)}e_{1}^{(a_{10})}\qquad
(\mod v^{-1}{\Cal L}) \endalign$$ $$\align &\quad \text{\rm
if}\quad a_8+a_9\geq a_5+a_6,\quad a_6+a_7\geq a_2+a_3,\quad
a_2\geq a_6,\cr &\qquad \qquad \qquad a_5\geq a_1\geq a_8,\quad
a_{10}\geq a_8,
\endalign$$ and $(3)$ of case $\bold 1$ is proved.

\bigskip\par
\head {3. Polynomial elements in several variables}\endhead
Finally, we conclude our note with the following two remarks:

\proclaim{Remark 1} {\rm The proof of all the other cases in
Theorem 1.3 is quite similar to that of the above case. Also, the
proof of a polynomial element has close relations with that of the
corresponding monomial element.}\endproclaim

\proclaim{Remark 2}{\rm We see that the regions of $62$ monomial
elements and $144$ polynomial elements in one variable in the
canonical basis $\bold B$ do not fill the space ${\Bbb N}^{10}$.
Actually, we need only to fill the space ${\Bbb N}^{9}$ because
the regions what we consider are independent of $a_4$. Recently,
Marsh claims (see [M]) that there should be $672$ regions
corresponding to elements in the canonical basis $\bold B$ of the
quantized enveloping algebra for type $A_4$. So we believe that in
addition to monomial elements and polynomial elements in one
variable we have worked out, there exist many polynomial elements
in two or more variables in the canonical basis $\bold B$. And we
have been keeping on the computations of polynomial elements in
two variables in the canonical basis $\bold B$. We do have found
more than thirty polynomial elements in two independent variables
in the canonical basis $\bold B$. Here, so-called ``independent
variables" means that the summing in the following two polynomials
is independent of the order of $u$ and $w$. For example,
corresponding to monomial element {\bf 4.}(1) in Theorem 3.1 in
[2], we list the following two polynomial elements in two
independent variables $u$ and $w$:\par
$$\align
&\sum \Sb 0\leq u\leq a_7+a_9 \\ 0\leq w \leq a_2+a_3+a_4 \endSb
(-1)^{u+w}\left[\smallmatrix a_{10}-a_8-1+u \cr u
\endsmallmatrix\right]\left[\smallmatrix a_5-a_1-1+w \cr w
\endsmallmatrix\right]e_1^{(a_4)}e_2^{(a_3+a_4)}\cr
&\qquad \qquad \qquad \times e_3^{(a_2+a_3+a_4-w)}e_2^{(a_6)}
e_1^{(a_7+a_9-u)} e_4^{(a_1+a_2+a_3+a_4)}\cr&\cr &\qquad \qquad
\qquad \times
e_2^{(a_7)}e_3^{(a_5+a_6+a_7+w)}e_2^{(a_8+a_9)}e_1^{(a_{10}+u)},\tag
a \cr&\cr &\qquad \text{\rm if}\quad a_3+a_6+a_8\geq
a_7+a_9+a_{10},\quad a_1+a_6\geq a_8+a_9,\cr &\qquad\qquad \;\;
a_{10}\geq a_8,\quad a_9\geq a_6,\quad a_2\geq a_6,\quad a_5\geq
a_1,
\endalign $$
and
$$\align &\sum \Sb 0\leq u\leq a_2+a_3+a_4
\\ 0\leq w\leq a_4 \endSb
(-1)^{u+w}\left[\smallmatrix a_5-a_1-1+u\cr u
\endsmallmatrix
\right]\left[\smallmatrix a_7+a_9-a_3-a_6-1+w\cr w
\endsmallmatrix
\right]e_1^{(a_4-w)}\cr & \qquad \qquad \qquad \times
e_2^{(a_3+a_4)}e_3^{(a_2+a_3+a_4-u)} e_2^{(a_6)}
e_1^{(a_7+a_9+w)}\cr &\cr&\qquad\qquad\qquad\times
e_4^{(a_1+a_2+a_3+a_4)}e_2^{(a_7)}e_3^{(a_5+a_6+a_7+u)}e_2^{(a_8+a_9)}
e_1^{(a_{10})},\tag b\cr&\cr \endalign $$$$\align
&\qquad\quad\text{\rm if}\qquad a_7+a_9\geq a_3+a_6,\quad
a_1+a_6\geq a_8+a_9,\quad a_5\geq a_1,\cr &\qquad\qquad\qquad
\qquad a_8\geq a_{10},\quad a_2\geq a_6,\quad a_3\geq a_7.
\endalign $$}
\endproclaim

{\head{Acknowledgement}\endhead}

This work is supported in part by the National Natural Science
Foundation of China (10271088) and the Natural Science Foundation
of Henan Province (0311010100). The first named author would like
to thank Professor Nanhua Xi and Professor Kaiming Zhao for their
financial support and for their valuable advice and comments during
his visit the Morning Side Centre of Mathematics in the Academy of
Mathematics and System Sciences in Beijing in May--September 2001.
The second named author is also grateful to the Abdus Salam International
Centre for Theoretical Physics for its financial support and hospitability
during his visit. Finally, both authors would like to thank Professor
Robert Marsh for helpful communications, and for his sending us [1]
and other papers.

\end{document}
21-08-2003
\input amstex
\documentstyle{amsppt}

\def\mod{\text{\rm mod}\, }

\def\i{\text{\rm i}}
\def\ii{\text{\rm ii}}
\def\iii{\text{\rm iii}}
\def\iv{\text{\rm iv}}
\def\v{\text{\rm v}}
\def\vi{\text{\rm vi}}

\def\bN{\Bbb N}

\font\boldtitlefont=cmb10 scaled\magstep1\font\boldfont=cmb10
\font\bigmath=cmbxti10 scaled\magstep2 \font\smmath=cmr10
 
\font\scten=cmti10
\topmatter\magnification=\magstep1
\hsize14truecm
\vsize21.6truecm
\leftheadtext {Hu and Ye}
\rightheadtext {Canonical Basis for Type ${\bigmath A}_{\bigmath
4}\;\;({\hbox {\rm II}})$}
\topmatter
\title
{\boldtitlefont  CANONICAL BASIS FOR TYPE ${\hbox{\bigmath
A}}_{\hbox{\boldfont 4}}\; ({\hbox{\boldtitlefont II}})$\\
$-$ Polynomial  elements in \\ one variable}
\endtitle
\author
{Hu Yuwang$^{1}$~ and ~Ye Jiachen$^{2~3}$}
\endauthor
\affil
{\scten $^1$Department of Mathematics, Xinyang Teachers College,\\
Henan 464000, People's Republic of China\\
E-Mail: hsh0412\@sina.com\\
\smallskip
$^2$Department of Applied Mathematics, Tongji University,\\
Shanghai 200092, People's Republic of China\\
E-Mail: jcye\@sh.cnuninet.net\\
\smallskip
$^3$The Abdus Salam International Centre for Theoretical Physics,\\
34014 Trieste, Italy\\}
\endaffil
\endtopmatter
\noindent {\bf Abstract:}\quad {\smmath All the $62$ monomial
elements in the canonical basis $\bold B$ of the quantized
enveloping algebra for type $A_4$ have been determined in [2].
According to Lusztig's idea [7], the elements in the canonical
basis $\bold B$  consist of monomials and linear combinations of
monomials (for convenience, we call them polynomials). In this
note, we compute all the $144$ polynomial elements in one variable
in the canonical basis $\bold B$ of the quantized enveloping
algebra for type $A_4$ based on our joint note [2]. We conjecture
that there are other polynomial elements in two or three variables
in the canonical basis $\bold B$, which include independent
variables and dependent variables. Moreover, it is conjectured
that there are no polynomial elements in the canonical basis
$\bold B$ with four or more variables.}
\medskip
\noindent {\bf Keywords:} {\smmath Canonical basis,\quad
Linear-combination,\quad Polynomials.}
\medskip
\noindent {\bf 2000 MR Subject Classification:} 17B37\quad 20G42\quad 81R50
\vskip1cm
\document

\head {1. Polynomial elements in one variable in the canonical
basis of $U^+$}\endhead

We shall freely use the notations in [2] without further comments.

\par \vskip0.1cm

{\bf 1.1.} In [2], we have determined the $62$ monomial elements
in the canonical basis $\bold B$ of the quantized enveloping
algebra for type $A_4$. Each monomial element corresponds to a
region which consists of six independent inequalities. Here
independence of six inequalities implies that we can't deduce one
inequality from the others. Moreover, the regions can't be
described with less than six inequalities, and the interiors of
any two regions are disjoint. These regions of monomial elements
have ``nice" forms, which will help us compute polynomial elements
in the canonical basis $\bold B$.

All the $144$ polynomial elements in one variable in the canonical
basis $\bold B$ will be determined in this note based on our joint
work [2]. We have been keeping on the computations of polynomial
elements in several variables in the canonical basis $\bold B$.
And we do have found more than thirty polynomial elements in two
independent variables $u$ and $w$ in the canonical basis $\bold
B$. We conjecture that there are other polynomial elements in two
or three variables in the canonical basis $\bold B$, which include
independent variables and dependent variables. Moreover, it is
conjectured that there are no polynomial elements in the canonical
basis $\bold B$ with four or more variables.

It should be mentioned here that so-called \lq\lq independent
variables\rq\rq ~is that the summing in the polynomial is
independent of the order of variables, and \lq\lq dependent
variables\rq\rq ~implies that the summing in the polynomial is
dependent of the order of variables.

\par\vskip0.1cm

{\bf 1.2.} Let us compute polynomial elements in one variable in
the canonical basis $\bold B$ of the quantized enveloping algebra
for type $A_4$. According to Lusztig's idea, these polynomial
elements should be linear combinations of the monomial elements.
We shall describe how one can compute them and begin with monomial
element {\bf 1.}(1) in [2] as an example. Let
$A=(a_{1},a_{2},a_{3},a_{4},a_{5},a_{6},a_{7},a_{8},a_{9},a_{10})\in
{\Bbb N}^{10}  $, we have known from [2] P. 242 that in order to
make the linear form $l(x_1, x_2, \cdots , x_{14})$ and the unit
form $q(x_1, x_2, \cdots , x_{14})$ to be non-negative for any
$x=(x_1, x_2, \cdots , x_{14})\in \Bbb{N}^{14}$, the following
inequalities
$$\align
& \quad a_5+a_6+a_7\geq a_1+a_2+a_3,\qquad a_7+a_9\geq a_3+a_6,
\cr & \quad a_7+a_9+a_{10}\geq a_3+a_6+a_8, \qquad a_6+a_7\geq
a_2+a_3,\cr (\sharp)&\quad
 a_7+a_8+a_9\geq a_2+a_3+a_5 ,\qquad a_2+a_5\geq
a_6+a_8,\cr& \quad a_8+a_9\geq a_5+a_6,\qquad a_9+a_{10}\geq
a_6+a_8,\cr & \quad a_1\geq a_5,\quad a_2\geq a_6,\quad a_5\geq
a_8, \quad a_{10}\geq a_8, \quad a_9\geq a_6,
\endalign $$
equivalently,
$$\align
& (*) \quad\;\; a_5\geq a_8,\qquad \;\;a_{10}\geq a_8,\qquad
\;\;a_2\geq a_6,\qquad \;\;a_1\geq a_5,\cr &\quad \quad
\;\;a_8+a_9\geq a_5+a_6, \quad \quad \;\;a_5+a_6+a_7\geq
a_1+a_2+a_3
\endalign $$
must hold, and then the monomial
$$\align
&(**)\qquad
e_{2}^{(a_3)}e_{3}^{(a_2+a_3)}e_{4}^{(a_1+a_2+a_3)}e_{2}^{(a_6)}
e_{3}^{(a_5+a_6)}e_{2}^{(a_8)}e_{1}^{(a_4+a_7+a_9+a_{10})}\times\\
&\qquad \qquad \times e_{2}^{(a_4+a_7+a_9)}
e_{3}^{(a_4+a_7)}e_{4}^{(a_4)}
\endalign $$
belongs to $\bold B$. Because it is impossible that we deduce any
one from the other inequalities in $(*)$, the six inequalities in
$(*)$ are independent. $(*)$ is called the region of the monomial
$(**)$, and the six independent inequalities in $(*)$ are called
the defining inequalities of the region $(*)$.  \par \vskip0.1cm

We now consider the linear combination of the monomial $(**)$ that
could become a new member in $\bold B$. Because there is an
one-to-one correspondence between ${\Bbb N}^{10}$ and $\bold B$,
we have to observe the regions similar to  $(*)$. Firstly, we
reverse the first defining inequality $a_5 \geq a_8$ in $(*)$, and
get a new inequality $a_8 \geq a_5$. In order to make the other
inequalities hold in $(\sharp)$ but $a_5 \geq a_8$, we have to
replace the remaining five defining inequalities in $(*)$ by the
following five defining inequalities

$$\align
& \qquad \;\;a_{10}\geq a_8,\qquad \;\;a_2+a_5\geq a_6+a_8,\qquad
\;\;a_1\geq a_5,\cr &\quad \quad \;\;a_9\geq a_6, \quad \quad
\;\;a_5+a_6+a_7\geq a_1+a_2+a_3.
\endalign $$
The five new defining inequalities together with $ a_8\geq a_5$
form a new region, which corresponds to a linear combination of
the monomial $(**)$. Similarly, we can deal with other two cases
that the second defining inequality and the third defining
inequality in $(*)$ are reversed, respectively.  Secondly, we
reverse the fourth defining inequality $a_1 \geq a_5$ in $(*)$,
and get a new inequality $a_5 \geq a_1$. In order to make the
other inequalities hold in $(\sharp)$ but $a_1 \geq a_5$, we have
to replace the remaining five defining inequalities in $(*)$ by
the following five defining inequalities
$$\align
& \qquad \;\;a_6+a_7\geq a_2+a_3,\qquad \;\;a_8+a_9\geq
a_5+a_6,\cr &\qquad  \quad a_2\geq a_6,\quad \quad \;\;a_1\geq
a_8, \quad \quad \;\;a_{10}\geq a_8.
\endalign $$
The five new defining inequalities together with $ a_5\geq a_1$
form a new region, which corresponds to a linear combination of
the monomial $(**)$. Similarly, we can deal with other two cases
that the fifth defining inequality and the sixth defining
inequality in $(*)$ are reversed, respectively.  By many
computations, we have noticed that the coefficient of the linear
combination of the monomial $(**)$ is closely related to the
reversed inequality. Therefore, corresponding to the six defining
inequalities in $(*)$, we get the following six possibilities,
respectively.
$$\align
& \qquad(\i)\;\; a_8\geq a_5,\qquad (\ii)\;\;a_8\geq a_{10},\qquad
(\iii)\;\;a_6\geq a_2,\qquad (\iv)\;\;a_5\geq a_1,\cr &\qquad
(\v)\;\;a_5+a_6\geq a_8+a_9, \quad \quad (\vi)\;\;a_1+a_2+a_3\geq
a_5+a_6+a_7,
\endalign $$
Then we get six polynomial elements in one variable, each
corresponds to one of the above six cases. Also, the region of
each of the six polynomial elements in one variable consists of
six independent inequalities. In this way, corresponding to each
of the 62 monomial elements, we compute all polynomial elements in
one variable and their regions.

Observing the regions of every six polynomial elements in one
variable which correspond to one of the 62 monomials, we find that
at least three regions may have already occurred in $\Cal S$, the
set of the $62$ regions determined by the $62$ monomial elements,
and there may be at most another two regions, which are the same
as those of polynomial elements in one variable we already
computed before. They will not make contribution to the canonical
basis $\bold B$.

In the above example, the regions corresponding to $(\i), (\ii),
(\iii)$ have already occurred in $\Cal S$, so the corresponding
three polynomial elements in one variable don't make contribution
to the canonical basis $\bold B$, although they are combinations
of the monomials $(**)$. In the other hand, the regions
corresponding to $(\iv),(\v),(\vi)$ don't occur in $\Cal S$, and
the corresponding three polynomial elements in one variable become
new members in the canonical basis $\bold B$ because this is the
first consideration for non-monomial case.

Repeating the above procedure for all $62$ monomials one by one,
we obtain all the $144$ polynomial elements in one variable, which
belong to the canonical basis $\bold B$. The main result
concerning with the polynomial elements in one variable in the
canonical basis $\bold B$ of the quantized enveloping algebra for
type $A_4$ is the following theorem.

\proclaim {Theorem 1.3} Let
$A=(a_1,a_2,a_3,a_4,a_5,a_6,a_7,a_8,a_9,a_{10})\in {\bN}^{10}$.
Then \par {\rm 1.}\; corresponding to $62$ equivalence classes for
$\sim $, the $72$ polynomial elements in one variable in the
canonical basis $\bold B$
$$\theta(A)=\theta(a_1,a_2,a_3,a_4,a_5,a_6,a_7,a_8,a_9,a_{10})\in
{\bold B}$$ are given by the following

$$\align {\bold 1.}\quad
&(1).\;\;\sum_{0\leq u\leq
a_4+a_7}\!\!\!(-1)^u\left[\!\smallmatrix a_5+a_6-a_8-a_9-1+u\cr u
\endsmallmatrix
\!\right]e_2^{(a_3)}e_3^{(a_2+a_3)}e_4^{(a_1+a_2+a_3)}e_2^{(a_6)}\cr
&\qquad\qquad\times
e_3^{(a_5+a_6+u)}e_2^{(a_8)}e_1^{(a_4+a_7+a_9+a_{10})}e_2^{(a_4+a_7+a_9)}
e_3^{(a_4+a_7-u)}e_4^{(a_4)}\cr &\qquad\qquad\text{\rm if}\qquad
a_5+a_6+a_7\geq a_1+a_2+a_3,\quad a_5+a_6\geq a_8+a_9,\cr
&\qquad\qquad\qquad \qquad  a_9\geq a_6,\quad a_2\geq a_6,\quad
a_1\geq a_5,\quad a_{10}\geq a_8.\cr&\cr &(2).\;\;\sum_{0\leq
u\leq a_4}\!\!(-1)^u \!\left[\!\smallmatrix
a_1+a_2+a_3-a_5-a_6-a_7-1+u\cr u
\endsmallmatrix
\!\right]e_2^{(a_3)}e_3^{(a_2+a_3)}e_4^{(a_1+a_2+a_3+u)}\cr &
\qquad\qquad\times
e_2^{(a_6)}e_3^{(a_5+a_6)}e_2^{(a_8)}e_1^{(a_4+a_7+a_9+a_{10})}e_2^{(a_4+a_7+a_9)}
e_3^{(a_4+a_7)}e_4^{(a_4-u)}\cr &\quad\quad\text{\rm if}\quad
a_1+a_2+a_3\geq a_5+a_6+a_7,\quad a_8+a_9\geq a_5+a_6,\quad
a_2\geq a_6,\cr & \qquad\qquad\qquad\qquad a_6+a_7\geq
a_2+a_3,\quad  a_5\geq a_8,\quad a_{10}\geq a_8.\cr
&(3).\;\;\sum_{0\leq u\leq a_2+a_3}\!\!(-1)^u\left[\smallmatrix
a_5-a_1-1+u\cr u
\endsmallmatrix
\right]e_2^{(a_3)}e_3^{(a_2+a_3-u)}e_4^{(a_1+a_2+a_3)}e_2^{(a_6)}\cr
&\qquad\qquad\times
e_3^{(a_5+a_6+u)}e_2^{(a_8)}e_1^{(a_4+a_7+a_9+a_{10})}e_2^{(a_4+a_7+a_9)}
e_3^{(a_4+a_7)}e_4^{(a_4)}\cr &\qquad\qquad\text{\rm if}\quad
a_6+a_7\geq a_2+a_3,\quad a_8+a_9\geq a_5+a_6,\quad a_2\geq
a_6,\cr &\qquad\qquad \qquad\qquad\qquad a_5\geq a_1\geq a_8,\quad
a_{10}\geq a_8. \cr&\cr {\bold 2.}\quad &(1).\;\;\sum_{0\leq u\leq
a_4+a_7+a_9}\!\!\!(-1)^u\left[\smallmatrix a_{10}-a_8-1+u\cr u
\endsmallmatrix
\right]e_3^{(a_2)}e_2^{(a_3+a_6)}e_1^{(a_4+a_7+a_9-u)}e_3^{(a_3)}\cr
&\qquad\qquad\quad\times
e_2^{(a_4+a_7)}e_3^{(a_4)}e_4^{(a_1+a_2+a_3+a_4)}e_3^{(a_5+a_6+a_7)}e_2^{(a_8+a_9)}
e_1^{(a_{10}+u)}\cr &\qquad\qquad\text{\rm if}\qquad
a_1+a_2+a_3\geq a_5+a_6+a_7,\quad a_5+a_6\geq a_8+a_9,\cr
&\quad\qquad \qquad\qquad\qquad a_9\geq a_6\geq a_2,\quad
a_{10}\geq a_8,\quad a_7\geq a_3.\cr &\cr&(2).\;\;\sum_{0\leq
u\leq a_4+a_7}\!\!\!(-1)^u\left[\smallmatrix
a_8+a_9-a_5-a_6-1+u\cr u
\endsmallmatrix
\right]e_3^{(a_2)}e_2^{(a_3+a_6)}e_1^{(a_4+a_7+a_9)}e_3^{(a_3)}\cr
&\quad\quad\quad\times
e_2^{(a_4+a_7-u)}e_3^{(a_4)}e_4^{(a_1+a_2+a_3+a_4)}e_3^{(a_5+a_6+a_7)}e_2^{(a_8+a_9+u)}
e_1^{(a_{10})}\cr &\qquad\qquad\text{\rm if}\qquad a_1+a_2+a_3\geq
a_5+a_6+a_7,\quad a_8+a_9\geq a_5+a_6,\cr &\quad\qquad
\qquad\qquad\qquad a_5\geq a_8\geq a_{10},\quad a_6\geq a_2,\quad
a_7\geq a_3.\cr&\cr &(3).\;\;\sum_{0\leq u\leq
a_4+a_7}\!\!\!(-1)^u\left[\smallmatrix a_6-a_9-1+u\cr u
\endsmallmatrix
\right]e_3^{(a_2)}e_2^{(a_3+a_6+u)}e_1^{(a_4+a_7+a_9)}e_3^{(a_3)}\cr
&\qquad\qquad\times
e_2^{(a_4+a_7-u)}e_3^{(a_4)}e_4^{(a_1+a_2+a_3+a_4)}e_3^{(a_5+a_6+a_7)}e_2^{(a_8+a_9)}
e_1^{(a_{10})}\cr &\qquad\qquad\text{\rm if}\qquad a_1+a_2+a_3\geq
a_5+a_6+a_7,\quad a_5\geq a_8\geq a_{10},\cr &\qquad\qquad\qquad
\qquad\quad a_6\geq a_9\geq a_2,\quad a_7\geq a_3.\cr &\cr
\endalign$$$$\align {\bold 3.}\quad &(1).\;\;\sum_{0\leq u\leq
a_3+a_4+a_6}\!\!\!\!\!\!\!\!\!\!(-1)^u\left[\smallmatrix
a_8-a_5-1+u\cr u
\endsmallmatrix
\right]e_1^{(a_4)}e_3^{(a_2)}e_2^{(a_3+a_4+a_6-u)}e_4^{(a_1+a_2)}\cr
&\qquad\qquad\times
e_1^{(a_7+a_9)}e_3^{(a_3+a_4+a_5+a_6)}e_2^{(a_7+a_8+a_9+u)}e_4^{(a_3+a_4)}e_1^{(a_{10})}
e_3^{(a_7)}\cr &\qquad\qquad\text{\rm if}\qquad a_5+a_6\geq
a_1+a_2,\quad a_3+a_6\geq a_7+a_9,\quad a_1\geq a_5,\cr
&\qquad\qquad \qquad\qquad\qquad\qquad  a_8\geq a_5\geq
a_{10},\quad a_9 \geq a_6.\cr&\cr
 &(2).\;\;\sum_{0\leq u
\leq a_7}\!(-1)^u\left[\smallmatrix a_5+a_6-a_8-a_9-1+u\cr u
\endsmallmatrix
\right]e_1^{(a_4)}e_3^{(a_2)}e_2^{(a_3+a_4+a_6)}e_4^{(a_1+a_2)}\cr
&\qquad\qquad\times
e_1^{(a_7+a_9)}e_3^{(a_3+a_4+a_5+a_6+u)}e_2^{(a_7+a_8+a_9)}e_4^{(a_3+a_4)}e_1^{(a_{10})}
e_3^{(a_7-u)}\cr &\qquad\qquad\text{\rm if}\qquad a_5+a_6\geq
a_8+a_9\geq a_1+a_2,\quad a_3+a_6\geq a_7+a_9,\cr
&\qquad\qquad\qquad \qquad\qquad\qquad  a_8\geq a_{10},\quad
a_1\geq a_5,\quad a_9\geq a_6.\cr &(3).\;\;\sum_{0\leq u\leq
a_2}\!(-1)^u\left[\smallmatrix a_5-a_1-1+u\cr u
\endsmallmatrix
\right]e_1^{(a_4)}e_3^{(a_2-u)}e_2^{(a_3+a_4+a_6)}e_4^{(a_1+a_2)}e_1^{(a_7+a_9)}\cr
&\qquad\qquad\qquad\qquad\times
e_3^{(a_3+a_4+a_5+a_6+u)}e_2^{(a_7+a_8+a_9)}e_4^{(a_3+a_4)}e_1^{(a_{10})}e_3^{(a_7)}\cr
&\qquad\qquad\text{\rm if}\qquad  a_8+a_9\geq a_5+a_6,\quad
a_3+a_6\geq a_7+a_9,\cr &\qquad\qquad\qquad\qquad   a_8\geq
a_1\geq a_5\geq a_{10},\quad a_6\geq a_2.\cr&\cr {\bold 4.}\quad
&(1).\;\;\sum_{0\leq u\leq
a_2+a_3+a_4}\!\!\!\!\!\!\!\!(-1)^u\left[\smallmatrix a_5-a_1-1+u
\cr u
\endsmallmatrix
\right]e_1^{(a_4)}e_2^{(a_3+a_4)}e_3^{(a_2+a_3+a_4-u)}e_2^{(a_6)}\cr
&\qquad\qquad\times
e_1^{(a_7+a_9)}e_4^{(a_1+a_2+a_3+a_4)}e_2^{(a_7)}e_3^{(a_5+a_6+a_7+u)}e_2^{(a_8+a_9)}
e_1^{(a_{10})}\cr &\qquad\qquad\text{\rm if}\qquad a_3+a_6\geq
a_7+a_9,\quad a_1+a_6\geq a_8+a_9,\quad a_2\geq a_6,\cr
&\qquad\qquad\qquad \qquad a_8\geq a_{10},\quad a_5\geq a_1,\quad
a_9\geq a_6 .\cr&\cr
  &(2).\;\;\sum_{0\leq u\leq
a_7+a_9}\!\!\!\!\!\!\!\!(-1)^u\left[\smallmatrix a_{10}-a_8-1+u
\cr u
\endsmallmatrix
\right]e_1^{(a_4)}e_2^{(a_3+a_4)}e_3^{(a_2+a_3+a_4)}e_2^{(a_6)}e_1^{(a_7+a_9-u)}\cr
&\qquad\qquad\qquad\qquad\qquad\times
e_4^{(a_1+a_2+a_3+a_4)}e_2^{(a_7)}e_3^{(a_5+a_6+a_7)}e_2^{(a_8+a_9)}e_1^{(a_{10}+u)}\cr
&\qquad\qquad\text{\rm if}\qquad a_3+a_6+a_8\geq
a_7+a_9+a_{10},\quad a_5+a_6\geq a_8+a_9,\cr &\qquad\qquad
\qquad\qquad a_{10}\geq a_8,\quad a_9\geq a_6,\quad a_2\geq
a_6,\quad a_1\geq a_5.\cr&\cr &(3).\;\;\sum_{0\leq u\leq
a_4}\!\!\!(-1)^u\left[\smallmatrix a_7+a_9-a_3-a_6-1+u\cr u
\endsmallmatrix
\right]e_1^{(a_4-u)}e_2^{(a_3+a_4)}e_3^{(a_2+a_3+a_4)}e_2^{(a_6)}\cr
&\qquad\qquad\qquad\times
e_1^{(a_7+a_9+u)}e_4^{(a_1+a_2+a_3+a_4)}e_2^{(a_7)}e_3^{(a_5+a_6+a_7)}e_2^{(a_8+a_9)}
e_1^{(a_{10})}\cr &\qquad\qquad\text{\rm if}\qquad a_7+a_9\geq
a_3+a_6,\quad a_5+a_6\geq a_8+a_9,\quad a_1\geq a_5,\cr
&\qquad\qquad\qquad \qquad a_8\geq a_{10},\quad a_2\geq a_6,\quad
a_3\geq a_7 .\cr&\cr \endalign $$$$\align {\bold 5.}\quad
&(1).\;\;\sum_{0\leq u\leq a_3+a_4}\!\!\!(-1)^u\left[\smallmatrix
a_5+a_6-a_1-a_2-1+u\cr u
\endsmallmatrix
\right]e_1^{(a_4)}e_3^{(a_2)}e_2^{(a_3+a_4+a_6)}e_1^{(a_7+a_9)}\cr
&\qquad\qquad\times
e_3^{(a_3+a_4-u)}e_4^{(a_1+a_2+a_3+a_4)}e_2^{(a_7)}e_3^{(a_5+a_6+a_7+u)}e_2^{(a_8+a_9)}
e_1^{(a_{10})}\cr &\qquad\qquad\text{\rm if}\qquad a_5+a_6\geq
a_1+a_2\geq a_8+a_9,\quad a_3+a_6\geq a_7+a_9,\cr
&\qquad\qquad\qquad \qquad  a_8\geq a_{10},\quad a_9\geq a_6,\quad
a_1\geq a_5 .\cr &\cr&(2).\;\;\sum_{0\leq u\leq
a_7+a_9}\!\!\!(-1)^u \left[\smallmatrix a_{10}-a_8-1+u\cr u
\endsmallmatrix
\right]e_1^{(a_4)}e_3^{(a_2)}e_2^{(a_3+a_4+a_6)}e_1^{(a_7+a_9-u)}\cr
&\qquad\qquad\times
e_3^{(a_3+a_4)}e_4^{(a_1+a_2+a_3+a_4)}e_2^{(a_7)}e_3^{(a_5+a_6+a_7)}e_2^{(a_8+a_9)}
e_1^{(a_{10}+u)}\cr &\qquad\qquad\text{\rm if}\qquad
a_3+a_6+a_8\geq a_7+a_9+a_{10},\quad a_9\geq a_6\geq a_2,\cr
&\qquad\qquad \qquad\qquad  a_1+a_2\geq a_5+a_6\geq a_8+a_9,\quad
a_{10}\geq a_8.\cr &(3).\;\;\sum_{0\leq u\leq
a_7}\!\!\!(-1)^u\left[\smallmatrix a_6-a_9-1+u\cr u
\endsmallmatrix
\right]e_1^{(a_4)}e_3^{(a_2)}e_2^{(a_3+a_4+a_6+u)}e_1^{(a_7+a_9)}e_3^{(a_3+a_4)}\cr
&\qquad\qquad\qquad\qquad\times
e_4^{(a_1+a_2+a_3+a_4)}e_2^{(a_7-u)}e_3^{(a_5+a_6+a_7)}e_2^{(a_8+a_9)}e_1^{(a_{10})}\cr
&\qquad\qquad\text{\rm if}\qquad a_1+a_2\geq a_5+a_6,\quad a_6\geq
a_9\geq a_2,\quad a_3\geq a_7,\cr &\qquad\qquad\qquad
\qquad\qquad\qquad a_5\geq a_8\geq a_{10}.\cr &\cr
 {\bold 6.}\quad
&(1).\;\;\sum_{0\leq u\leq a_6}(-1)^u\left[\smallmatrix
a_8-a_5-1+u\cr u
\endsmallmatrix
\right]e_1^{(a_4)}e_2^{(a_3+a_4)}e_3^{(a_2+a_3+a_4)}e_2^{(a_6-u)}\cr
&\qquad\quad\times
e_4^{(a_1+a_2+a_3+a_4)}e_3^{(a_5+a_6)}e_1^{(a_7+a_9)}e_2^{(a_7+a_8+a_9+u)}e_3^{(a_7)}
e_1^{(a_{10})}\cr &\qquad\text{\rm if}\quad a_3+a_6\geq
a_7+a_9,\quad a_2+a_5 \geq a_6+a_8,\quad a_1\geq a_5,\cr
&\qquad\qquad\qquad\qquad a_8\geq a_5\geq a_{10},\quad a_9\geq
a_6.\cr&\cr &(2).\;\;\sum_{0\leq u\leq
a_4}\!\!\!(-1)^u\left[\smallmatrix a_7+a_9-a_3-a_6-1+u\cr u
\endsmallmatrix
\right]e_1^{(a_4-u)}e_2^{(a_3+a_4)}e_3^{(a_2+a_3+a_4)}e_2^{(a_6)}\cr
&\qquad\qquad\times
e_4^{(a_1+a_2+a_3+a_4)}e_3^{(a_5+a_6)}e_1^{(a_7+a_9+u)}e_2^{(a_7+a_8+a_9)}e_3^{(a_7)}
e_1^{(a_{10})}\cr &\quad\text{\rm if}\quad a_3+a_5+a_6\geq
a_7+a_8+a_9,\quad a_7+a_9\geq a_3+a_6,\quad a_1\geq a_5,\cr
&\qquad\qquad\qquad \qquad a_8+a_9\geq a_5+a_6,\quad a_8\geq
a_{10},\quad a_2\geq a_6.\cr&\cr &(3).\;\;\sum_{0\leq u\leq
a_2+a_3+a_4}\!\!\!(-1)^u\left[\smallmatrix a_5-a_1-1+u\cr u
\endsmallmatrix
\right]e_1^{(a_4)}e_2^{(a_3+a_4)}e_3^{(a_2+a_3+a_4-u)}e_2^{(a_6)}\cr
&\qquad\qquad\qquad\times
e_4^{(a_1+a_2+a_3+a_4)}e_3^{(a_5+a_6+u)}e_1^{(a_7+a_9)}e_2^{(a_7+a_8+a_9)}e_3^{(a_7)}
e_1^{(a_{10})}\cr &\qquad\qquad\text{\rm if}\qquad a_8+a_9\geq
a_5+a_6,\quad a_3+a_6\geq a_7+a_9,\cr &\qquad\qquad\qquad\qquad
a_5\geq a_1\geq a_8\geq a_{10},\quad  a_2\geq a_6.\cr
&\cr\endalign$$

$$\align
 {\bold 7.}\quad &(1).\;\;\sum_{0\leq u\leq
a_7}\!\!\!(-1)^u\left[\smallmatrix a_5+a_6-a_8-a_9-1+u\cr u
\endsmallmatrix
\right]e_3^{(a_2)}e_2^{(a_3+a_6)}e_1^{(a_4+a_7+a_9)}e_2^{(a_4)}\cr
&\qquad\qquad\times
e_4^{(a_1+a_2)}e_3^{(a_3+a_4+a_5+a_6+u)}e_2^{(a_7+a_8+a_9)}e_4^{(a_3+a_4)}e_1^{(a_{10})}
e_3^{(a_7-u)}\cr &\qquad\qquad\text{\rm if}\qquad a_5+a_6\geq
a_8+a_9\geq a_1+a_2,\quad a_7+a_9\geq a_3+a_6,\cr
&\qquad\qquad\qquad \qquad  a_1\geq a_5,\quad a_3\geq a_7,\quad
a_8\geq a_{10} .\cr &\cr &(2).\;\;\sum_{0\leq u\leq
a_2}\!\!\!(-1)^u \left[\smallmatrix a_5-a_1-1+u\cr u
\endsmallmatrix
\right]e_3^{(a_2-u)}e_2^{(a_3+a_6)}e_1^{(a_4+a_7+a_9)}e_2^{(a_4)}e_4^{(a_1+a_2)}\cr
&\qquad\qquad\qquad\times
e_3^{(a_3+a_4+a_5+a_6+u)}e_2^{(a_7+a_8+a_9)}e_4^{(a_3+a_4)}e_1^{(a_{10})}e_3^{(a_7)}\cr
&\qquad\text{\rm if}\qquad a_1+a_3+a_6\geq a_7+a_8+a_9,\quad
a_7+a_9\geq a_3+a_6,\quad a_5\geq a_1,\cr &\qquad\qquad\qquad
a_8+a_9\geq a_5+a_6,\quad a_6\geq a_2,\quad a_8\geq a_{10} .\cr
&(3).\;\;\sum_{0\leq u\leq a_4+a_7+a_9}\!\!\!(-1)^u
\left[\smallmatrix a_{10}-a_8-1+u\cr u
\endsmallmatrix
\right]e_3^{(a_2)}e_2^{(a_3+a_6)}e_1^{(a_4+a_7+a_9-u)}e_2^{(a_4)}\cr
&\qquad\qquad\qquad\times
e_4^{(a_1+a_2)}e_3^{(a_3+a_4+a_5+a_6)}e_2^{(a_7+a_8+a_9)}e_4^{(a_3+a_4)}e_1^{(a_{10}+u)}
e_3^{(a_7)}\cr
&\qquad\text{\rm if}\qquad a_3+a_5+a_6\geq a_7+a_8+a_9,\quad
a_7+a_9\geq a_3+a_6,\quad a_1\geq a_5,\cr &\qquad\qquad\qquad
a_8+a_9\geq a_5+a_6\geq a_1+a_2,\quad a_{10}\geq a_8.\cr &\cr
 {\bold 8.}\quad &(1).\sum_{0\leq u\leq a_4}\!\!\!(-1)^u
\left[\smallmatrix a_7+a_9+a_{10}-a_3-a_6-a_8-1+u\cr u
\endsmallmatrix
\right]e_3^{(a_2)}e_1^{(a_4-u)}e_4^{(a_1+a_2)}\cr &\quad \times
e_2^{(a_3+a_4+a_6)}e_3^{(a_3+a_4+a_5+a_6)}e_4^{(a_3+a_4)}e_2^{(a_8)}e_1^{(a_7+a_9+a_{10}+u)}
e_2^{(a_7+a_9)}e_3^{(a_7)}\cr &\qquad\text{\rm if}\qquad
a_7+a_9+a_{10}\geq a_3+a_6+a_8,\quad a_3+a_6\geq a_7+a_9,\cr
&\qquad\qquad\qquad a_8+a_9\geq a_5+a_6\geq a_1+a_2,\quad a_1\geq
a_5\geq a_8.\cr &\cr&(2).\;\;\sum_{0\leq u\leq
a_2}(-1)^u\left[\smallmatrix a_5-a_1-1+u\cr u
\endsmallmatrix
\right]e_3^{(a_2-u)}e_1^{(a_4)}e_4^{(a_1+a_2)}e_2^{(a_3+a_4+a_6)}\cr
&\quad\quad\times
e_3^{(a_3+a_4+a_5+a_6+u)}e_4^{(a_3+a_4)}e_2^{(a_8)}e_1^{(a_7+a_9+a_{10})}e_2^{(a_7+a_9)}
e_3^{(a_7)}\cr &\qquad\text{\rm if}\qquad a_3+a_6+a_8\geq
a_7+a_9+a_{10},\quad a_8+a_9\geq a_5+a_6,\cr &\qquad\qquad \qquad
a_5\geq a_1\geq a_8,\quad a_6\geq a_2,\quad a_{10}\geq a_8
.\cr&\cr &(3).\;\;\sum_{0\leq u\leq
a_7}\!\!\!(-1)^u\left[\smallmatrix a_5+a_6-a_8-a_9-1+u\cr u
\endsmallmatrix
\right]e_3^{(a_2)}e_1^{(a_4)}e_4^{(a_1+a_2)}e_2^{(a_3+a_4+a_6)}\cr
&\qquad\quad\times
e_3^{(a_3+a_4+a_5+a_6+u)}e_4^{(a_3+a_4)}e_2^{(a_8)}e_1^{(a_7+a_9+a_{10})}e_2^{(a_7+a_9)}
e_3^{(a_7-u)}\cr &\qquad\text{\rm if}\qquad a_3+a_6+a_8\geq
a_7+a_9+a_{10},\quad a_1\geq a_5,\quad a_{10}\geq a_8,\cr &
\qquad\qquad\qquad a_5+a_6\geq a_8+a_9\geq a_1+a_2,\quad a_9\geq
a_6.\cr &\cr\endalign $$

$$\align {\bold 9.}\quad &(1).\!\!\sum_{0\leq
u\leq a_6+a_7}\!\!\!(-1)^u \left[\smallmatrix a_8-a_5-1+u\cr u
\endsmallmatrix
\right]e_1^{(a_4)}e_2^{(a_3+a_4)}e_3^{(a_2+a_3+a_4)}e_1^{(a_7)}e_2^{(a_6+a_7-u)}\cr
&\qquad\qquad\qquad\quad\quad\times
e_1^{(a_9)}e_4^{(a_1+a_2+a_3+a_4)}e_3^{(a_5+a_6+a_7)}e_2^{(a_8+a_9+u)}e_1^{(a_{10})}\cr
&\qquad\qquad\qquad\text{\rm if}\qquad a_2+a_5\geq a_6+a_8,\quad
a_1\geq a_5\geq a_{10},\cr &\qquad\qquad \qquad\qquad\qquad
a_8\geq a_5,\quad a_3\geq a_7,\quad a_6\geq a_9.\cr &\cr
&(2).\sum_{0\leq u\leq a_2+a_3+a_4}\!\!\!(-1)^u \left[\smallmatrix
a_5-a_1-1+u\cr u
\endsmallmatrix
\right]e_1^{(a_4)}e_2^{(a_3+a_4)}e_3^{(a_2+a_3+a_4-u)}e_1^{(a_7)}\cr
&\qquad\qquad\times
e_2^{(a_6+a_7)}e_1^{(a_9)}e_4^{(a_1+a_2+a_3+a_4)}e_3^{(a_5+a_6+a_7+u)}e_2^{(a_8+a_9)}
e_1^{(a_{10})}\cr
&\qquad\qquad\text{\rm if}\qquad a_5\geq a_1\geq a_8\geq
 a_{10},\quad a_2\geq a_6\geq a_9,\quad a_3\geq a_7.\cr
&(3).\sum_{0\leq u\leq a_3+a_4}\!\!\!(-1)^u\left[\smallmatrix
a_6-a_2-1+u\cr u
\endsmallmatrix
\right]e_1^{(a_4)}e_2^{(a_3+a_4-u)}e_3^{(a_2+a_3+a_4)}e_1^{(a_7)}\cr
&\qquad\qquad\times
e_2^{(a_6+a_7+u)}e_1^{(a_9)}e_4^{(a_1+a_2+a_3+a_4)}e_3^{(a_5+a_6+a_7)}e_2^{(a_8+a_9)}
e_1^{(a_{10})}\cr &\qquad\qquad\qquad\text{\rm if}\qquad
a_1+a_2\geq a_5+a_6,\quad a_5\geq a_8\geq a_{10},\cr &\qquad\qquad
\qquad\qquad\qquad a_6\geq a_2\geq a_9,\quad a_3\geq a_7.\cr &\cr
{\bold {10.}}\quad &(1).\;\;\sum_{0\leq u\leq a_3+a_4}\!\!(-1)^u
\left[\smallmatrix a_6-a_2-1+u\cr u
\endsmallmatrix
\right]e_1^{(a_4)}e_2^{(a_3+a_4-u)}e_3^{(a_2+a_3+a_4)}e_1^{(a_7)}\cr
&\qquad\quad\times
e_2^{(a_6+a_7+u)}e_4^{(a_1+a_2+a_3+a_4)}e_3^{(a_5+a_6+a_7)}e_2^{(a_8)}e_1^{(a_9+a_{10})}
e_2^{(a_9)}\cr &\quad\qquad\qquad\text{\rm if}\quad a_1+a_2\geq
a_5+a_6,\quad a_2+a_8\geq a_9+a_{10},\cr &\qquad\qquad\quad \qquad
a_{10}\geq a_8,\quad a_5\geq a_8,\quad a_6\geq a_2,\quad a_3\geq
a_7.\cr &\cr&(2).\;\;\sum_{0\leq u\leq a_7}\!\!(-1)^u
\left[\smallmatrix a_9+a_{10}-a_6-a_8-1+u\cr u
\endsmallmatrix
\right]e_1^{(a_4)}e_2^{(a_3+a_4)}e_3^{(a_2+a_3+a_4)}e_1^{(a_7-u)}\cr
&\qquad\qquad\qquad\times
e_2^{(a_6+a_7)}e_4^{(a_1+a_2+a_3+a_4)}e_3^{(a_5+a_6+a_7)}e_2^{(a_8)}e_1^{(a_9+a_{10}+u)}
e_2^{(a_9)}\cr &\quad\qquad\qquad\text{\rm if}\quad
a_3+a_6+a_8\geq a_7+a_9+a_{10},\quad a_9+a_{10}\geq a_6+a_8,\cr
&\qquad\qquad \qquad\quad\qquad a_1\geq a_5\geq a_8,\quad a_2\geq
a_6\geq a_9 .\cr&\cr &(3).\;\;\sum_{0\leq u\leq
a_2+a_3+a_4}\!\!(-1)^u\left[\smallmatrix a_5-a_1-1+u\cr u
\endsmallmatrix
\right]e_1^{(a_4)}e_2^{(a_3+a_4)}e_3^{(a_2+a_3+a_4-u)}e_1^{(a_7)}\cr
&\qquad\qquad\qquad\times
e_2^{(a_6+a_7)}e_4^{(a_1+a_2+a_3+a_4)}e_3^{(a_5+a_6+a_7+u)}e_2^{(a_8)}e_1^{(a_9+a_{10})}
e_2^{(a_9)}\cr &\quad\qquad\qquad\text{\rm if}\quad a_6+a_8\geq
a_9+a_{10},\quad a_5\geq a_1\geq a_8,\cr &\qquad\qquad\qquad
\quad\qquad a_{10}\geq a_8,\quad a_3\geq a_7,\quad a_2\geq a_6
.\cr &\cr\endalign $$

$$\align{\bold {11.}}\quad
&(1).\;\;\sum_{0\leq u\leq a_3+a_6}\!\!(-1)^u\left[\smallmatrix
a_5-a_8-1+u\cr u
\endsmallmatrix
\right]e_4^{(a_1)}e_3^{(a_2+a_5+u)}e_2^{(a_3+a_6+a_8)}e_4^{(a_2)}\cr
&\qquad\qquad\times
e_3^{(a_3+a_6-u)}e_1^{(a_4+a_7+a_9+a_{10})}e_2^{(a_4+a_7+a_9)}e_3^{(a_4)}e_4^{(a_3+a_4)}
e_3^{(a_7)}\cr &\quad\qquad\qquad\qquad\text{\rm if}\quad
a_7+a_8+a_9\geq a_3+a_5+a_6,\quad a_5\geq a_8,\cr
&\qquad\qquad\qquad\quad \qquad a_{10}\geq a_8\geq a_1,\quad
a_6\geq a_2,\quad a_3\geq a_7 .\cr &\cr&(2).\;\;\sum_{0\leq u\leq
a_3+a_4}\!\!(-1)^u \left[\smallmatrix a_2-a_6-1+u\cr u
\endsmallmatrix
\right]e_4^{(a_1)}e_3^{(a_2+a_5)}e_2^{(a_3+a_6+a_8)}e_4^{(a_2+u)}\cr
&\qquad\qquad\times
e_3^{(a_3+a_6)}e_1^{(a_4+a_7+a_9+a_{10})}e_2^{(a_4+a_7+a_9)}e_3^{(a_4)}e_4^{(a_3+a_4-u)}
e_3^{(a_7)}\cr &\quad\qquad\text{\rm if}\quad a_7+a_9\geq
a_3+a_6,\quad a_6+a_8\geq a_2+a_5,\quad a_3\geq a_7,\cr
&\qquad\qquad\qquad \quad\qquad a_{10}\geq a_8,\quad  a_2\geq
a_6,\quad a_5\geq a_1 .\cr &(3).\;\;\sum_{0\leq u\leq
a_4+a_7+a_9}\!\!(-1)^u\left[\smallmatrix a_8-a_{10}-1+u\cr u
\endsmallmatrix
\right]e_4^{(a_1)}e_3^{(a_2+a_5)}e_2^{(a_3+a_6+a_8+u)}e_4^{(a_2)}\cr
&\qquad\qquad\qquad\times
e_3^{(a_3+a_6)}e_1^{(a_4+a_7+a_9+a_{10})}e_2^{(a_4+a_7+a_9-u)}e_3^{(a_4)}e_4^{(a_3+a_4)}
e_3^{(a_7)}\cr
&\quad\qquad\qquad\qquad\text{\rm if}\quad a_7+a_9\geq
a_3+a_6,\quad a_6\geq a_2,\cr &\qquad\qquad\qquad\quad\qquad
a_8\geq a_{10}\geq a_5\geq a_1,\quad  a_3\geq a_7.\cr &\cr
 {\bold {12.}}\quad &(1).\sum_{0\leq u\leq a_4+a_7}\!\!\!\!\!\!(-1)^u
\left[\smallmatrix a_5+a_6-a_8-a_9-1+u\cr u
\endsmallmatrix
\right]e_3^{(a_2)}e_2^{(a_3+a_6)}e_1^{(a_4+a_7+a_9)}e_4^{(a_1+a_2)}\cr
&\qquad\qquad\times
e_3^{(a_3+a_5+a_6+u)}e_2^{(a_4+a_7+a_8+a_9)}e_4^{(a_3)}e_3^{(a_4+a_7-u)}e_4^{(a_4)}
e_1^{(a_{10})}\cr &\qquad\qquad\qquad\text{\rm if}\quad
a_5+a_6\geq a_8+a_9\geq a_1+a_2,\quad a_7\geq a_3,\cr
&\qquad\qquad\qquad\qquad a_8\geq a_{10},\quad  a_1\geq a_5,\quad
a_9\geq a_6.\cr&\cr &(2).\;\;\sum_{0\leq u\leq
a_2}\!\!(-1)^u\left[\smallmatrix a_5-a_1-1+u\cr u
\endsmallmatrix
\right]e_3^{(a_2-u)}e_2^{(a_3+a_6)}e_1^{(a_4+a_7+a_9)}e_4^{(a_1+a_2)}\cr
&\qquad\qquad\times
e_3^{(a_3+a_5+a_6+u)}e_2^{(a_4+a_7+a_8+a_9)}e_4^{(a_3)}e_3^{(a_4+a_7)}e_4^{(a_4)}
e_1^{(a_{10})}\cr &\qquad\qquad\qquad\text{\rm if}\quad
a_8+a_9\geq a_5+a_6,\quad a_7\geq a_3,\cr
&\qquad\qquad\qquad\qquad a_5\geq a_1\geq a_8\geq a_{10},\quad
a_6\geq a_2.\cr &\cr &(3).\;\;\sum_{0\leq u\leq a_3+a_6}\!\!(-1)^u
\left[\smallmatrix a_8-a_5-1+u\cr u
\endsmallmatrix
\right]e_3^{(a_2)}e_2^{(a_3+a_6-u)}e_1^{(a_4+a_7+a_9)}e_4^{(a_1+a_2)}\cr
&\qquad\qquad\times
e_3^{(a_3+a_5+a_6)}e_2^{(a_4+a_7+a_8+a_9+u)}e_4^{(a_3)}e_3^{(a_4+a_7)}e_4^{(a_4)}
e_1^{(a_{10})}\cr &\qquad\qquad\qquad\text{\rm if}\quad
a_5+a_6\geq a_1+a_2,\quad a_1\geq a_5\geq a_{10},\cr
&\qquad\qquad\qquad \qquad a_8\geq a_5,\quad  a_7\geq a_3,\quad
a_9\geq a_6.\cr &\cr \endalign $$$$\align{\bold {13.}}\quad
&(1).\;\;\sum_{0\leq u\leq
a_4+a_6+a_7}\!\!(-1)^u\left[\smallmatrix a_8-a_5-1+u\cr u
\endsmallmatrix
\right]e_2^{(a_3)}e_3^{(a_2+a_3)}e_4^{(a_1+a_2+a_3)}e_1^{(a_4+a_7)}\cr
&\qquad\qquad\times
e_2^{(a_4+a_6+a_7-u)}e_3^{(a_4+a_5+a_6+a_7)}e_1^{(a_9)}e_2^{(a_8+a_9+u)}e_1^{(a_{10})}
e_4^{(a_4)}\cr &\qquad\qquad\qquad\text{\rm if}\quad
a_5+a_6+a_7\geq a_1+a_2+a_3,\quad  a_2+a_5\geq a_6+a_8,\cr
&\qquad\qquad\qquad \qquad a_1\geq a_5\geq a_{10},\quad  a_6\geq
a_9,\quad a_8\geq a_5 .\cr &\cr&(2).\;\;\sum_{0\leq u\leq
a_3}\!\!(-1)^u \left[\smallmatrix a_6-a_2-1+u\cr u
\endsmallmatrix
\right]e_2^{(a_3-u)}e_3^{(a_2+a_3)}e_4^{(a_1+a_2+a_3)}e_1^{(a_4+a_7)}\cr
&\qquad\qquad\times
e_2^{(a_4+a_6+a_7+u)}e_3^{(a_4+a_5+a_6+a_7)}e_1^{(a_9)}e_2^{(a_8+a_9)}e_1^{(a_{10})}
e_4^{(a_4)}\cr &\qquad\qquad\qquad\text{\rm if}\quad
a_5+a_6+a_7\geq a_1+a_2+a_3,\quad  a_1+a_2\geq a_5+a_6,\cr
&\qquad\qquad\qquad \qquad a_5\geq a_8\geq a_{10},\quad  a_6\geq
a_2\geq a_9.\cr &(3).\;\;\sum_{0\leq u\leq a_2+a_3}\!\!(-1)^u
\left[\smallmatrix a_5-a_1-1+u\cr u
\endsmallmatrix
\right]e_2^{(a_3)}e_3^{(a_2+a_3-u)}e_4^{(a_1+a_2+a_3)}e_1^{(a_4+a_7)}\cr
&\qquad\qquad\times
e_2^{(a_4+a_6+a_7)}e_3^{(a_4+a_5+a_6+a_7+u)}e_1^{(a_9)}e_2^{(a_8+a_9)}e_1^{(a_{10})}
e_4^{(a_4)}\cr
&\qquad\qquad\qquad\text{\rm if}\quad a_6+a_7\geq
a_2+a_3,\quad a_2\geq a_6\geq a_9,\cr &\qquad\qquad\qquad
\qquad a_5\geq a_1\geq a_8\geq a_{10}.\cr &\cr
 {\bold {14.}}\quad &(1).\sum_{0\leq u\leq
a_4+a_7+a_9}\!\!\!\!\!\!\!\!\!\!\!(-1)^u\left[\smallmatrix
a_{10}-a_8-1+u\cr u
\endsmallmatrix
\right]e_3^{(a_2)}e_2^{(a_3+a_6)}e_1^{(a_4+a_7+a_9-u)}e_2^{(a_4)}\cr
&\qquad\qquad\times
e_3^{(a_3+a_4)}e_4^{(a_1+a_2+a_3+a_4)}e_3^{(a_5+a_6)}e_2^{(a_7+a_8+a_9)}e_1^{(a_{10}+u)}
e_3^{(a_7)}\cr &\quad\text{\rm if}\quad a_3+a_5+a_6\geq
a_7+a_8+a_9,\quad a_7+a_9\geq a_3+a_6,\quad a_6\geq a_2,\cr
&\qquad\qquad \quad   a_8+a_9\geq a_5+a_6,\quad a_1+a_2\geq
a_5+a_6 ,\quad a_{10}\geq a_8.\cr &\cr&(2).\;\;\sum_{0\leq u\leq
a_4}\!\!\!(-1)^u\left[\smallmatrix a_7+a_8+a_9-a_3-a_5-a_6-1+u \cr
u
\endsmallmatrix
\right]e_3^{(a_2)}e_2^{(a_3+a_6)}e_1^{(a_4+a_7+a_9)}\cr &\qquad
\quad\times
e_2^{(a_4-u)}e_3^{(a_3+a_4)}e_4^{(a_1+a_2+a_3+a_4)}e_3^{(a_5+a_6)}e_2^{(a_7+a_8+a_9+u)}
e_1^{(a_{10})}e_3^{(a_7)}\cr &\quad\text{\rm if}\quad
a_7+a_8+a_9\geq a_3+a_5+a_6,\quad a_1+a_2\geq a_5+a_6,\quad
a_3\geq a_7,\cr &\qquad\qquad \qquad\quad   a_5\geq a_8\geq
a_{10},\quad a_6\geq a_2.\cr&\cr {\bold
{15.}}\quad&(1).\sum_{0\leq u\leq a_4}\!\!\!(-1)^u
\left[\smallmatrix a_7+a_9+a_{10}-a_3-a_6-a_8-1+u\cr u
\endsmallmatrix
\right]e_3^{(a_2)}e_1^{(a_4-u)}e_2^{(a_3+a_4+a_6)}\cr &\quad
\quad\times
e_3^{(a_3+a_4)}e_4^{(a_1+a_2+a_3+a_4)}e_3^{(a_5+a_6)}e_2^{(a_8)}e_1^{(a_7+a_9+a_{10}+u)}
e_2^{(a_7+a_9)}e_3^{(a_7)}\cr &\quad\text{\rm if}\quad
a_7+a_9+a_{10}\geq a_3+a_6+a_8,\quad a_1+a_2\geq a_5+a_6,\quad
a_6\geq a_2,\cr &\qquad\qquad\qquad \quad a_8+a_9\geq
a_5+a_6,\quad a_3+a_6\geq a_7+a_9,\quad a_5\geq a_8.\cr &\cr
\endalign $$$$\align &(2).\;\;\sum_{0\leq u\leq
a_3+a_4+a_6}\!\!\!(-1)^u \left[\smallmatrix a_8-a_5-1+u\cr u
\endsmallmatrix
\right]e_3^{(a_2)}e_1^{(a_4)}e_2^{(a_3+a_4+a_6-u)}e_3^{(a_3+a_4)}\cr
&\qquad\qquad\qquad\times
e_4^{(a_1+a_2+a_3+a_4)}e_3^{(a_5+a_6)}e_2^{(a_8+u)}e_1^{(a_7+a_9+a_{10})}e_2^{(a_7+a_9)}
e_3^{(a_7)}\cr &\qquad\text{\rm if}\qquad a_3+a_6+a_8\geq
a_7+a_9+a_{10},\quad a_1+a_2\geq a_5+a_6,\cr &\qquad\qquad
\qquad\quad a_9\geq a_6\geq a_2,\quad a_{10}\geq a_8\geq
a_5.\cr&\cr  {\bold {16.}}\quad &(1).\sum_{0\leq u\leq
a_3+a_4}\!\!\!(-1)^u \left[\smallmatrix a_5+a_6-a_1-a_2-1+u\cr u
\endsmallmatrix
\right]e_3^{(a_2)}e_2^{(a_3+a_6)}e_1^{(a_4+a_7+a_9)}e_2^{(a_4)}\cr
&\qquad\quad\times
e_3^{(a_3+a_4-u)}e_4^{(a_1+a_2+a_3+a_4)}e_2^{(a_7)}e_3^{(a_5+a_6+a_7+u)}e_2^{(a_8+a_9)}
e_1^{(a_{10})}\cr &\qquad\text{\rm if}\qquad a_5+a_6\geq
a_1+a_2\geq a_8+a_9,\quad a_7+a_9\geq a_3+a_6,\cr
&\qquad\qquad\qquad \quad a_8\geq a_{10},\quad a_3\geq a_7,\quad
a_1\geq a_5.\cr &(2).\sum_{0\leq u\leq a_4+a_7+a_9}\!\!\!(-1)^u
\left[\smallmatrix a_{10}-a_8-1+u\cr u
\endsmallmatrix
\right]e_3^{(a_2)}e_2^{(a_3+a_6)}e_1^{(a_4+a_7+a_9-u)}e_2^{(a_4)}\cr
&\qquad\quad\times
e_3^{(a_3+a_4)}e_4^{(a_1+a_2+a_3+a_4)}e_2^{(a_7)}e_3^{(a_5+a_6+a_7)}e_2^{(a_8+a_9)}
e_1^{(a_{10}+u)}\cr &\qquad\text{\rm if}\qquad a_1+a_2\geq
a_5+a_6\geq a_8+a_9,\quad a_7+a_9\geq a_3+a_6,\cr
&\qquad\qquad\qquad \quad a_{10}\geq a_8,\quad a_3\geq a_7,\quad
a_6\geq a_2.\cr&\cr {\bold {17.}}\quad&(1).\sum_{0\leq u\leq
a_4}\!\!\!(-1)^u \left[\smallmatrix a_1+a_2+a_3-a_5-a_6-a_7-1+u\cr
u
\endsmallmatrix
\right]e_2^{(a_3)}e_3^{(a_2+a_3)}e_4^{(a_1+a_2+a_3+u)}\cr &\qquad
\quad\times
e_3^{(a_5)}e_2^{(a_6+a_8)}e_1^{(a_4+a_7+a_9+a_{10})}e_3^{(a_6)}e_2^{(a_4+a_7+a_9)}
e_3^{(a_4+a_7)}e_4^{(a_4-u)}\cr &\quad\text{\rm if}\quad
a_1+a_2+a_3\geq a_5+a_6+a_7,\quad a_2+a_5\geq a_6+a_8,\quad
a_9\geq a_6,\cr &\qquad\qquad\qquad \quad a_6+a_7\geq
a_2+a_3,\quad  a_{10}\geq a_8\geq  a_5.\cr
&\cr&(2).\;\;\sum_{0\leq u\leq a_4+a_7+a_9}\!\!\!(-1)^u
\left[\smallmatrix a_8-a_{10}-1+u\cr u
\endsmallmatrix
\right]e_2^{(a_3)}e_3^{(a_2+a_3)}e_4^{(a_1+a_2+a_3)}e_3^{(a_5)}\cr
&\qquad\times
e_2^{(a_6+a_8+u)}e_1^{(a_4+a_7+a_9+a_{10})}e_3^{(a_6)}e_2^{(a_4+a_7+a_9-u)}e_3^{(a_4+a_7)}
e_4^{(a_4)}\cr &\qquad\text{\rm if}\qquad a_5+a_6+a_7\geq
a_1+a_2+a_3,\quad a_2+a_5\geq a_6+a_8,\cr &\qquad\qquad\qquad\quad
a_1\geq a_5,\quad a_8\geq a_{10}\geq a_5,\quad a_9\geq a_6.\cr
&\cr  {\bold {18.}}\quad &(1).\;\;\sum_{0\leq u \leq
a_3}\!\!(-1)^u\left[\smallmatrix a_6-a_2-1+u\cr u
\endsmallmatrix
\right]e_2^{(a_3-u)}e_3^{(a_2+a_3)}e_4^{(a_1+a_2+a_3)}e_1^{(a_4+a_7)}\cr
&\qquad\qquad\qquad\times
e_2^{(a_4+a_6+a_7+u)}e_3^{(a_4+a_5+a_6+a_7)}e_2^{(a_8)}e_1^{(a_9+a_{10})}e_4^{(a_4)}
e_2^{(a_9)}\cr
&\quad\text{\rm if}\quad a_5+a_6+a_7\geq a_1+a_2+a_3,\quad
a_2+a_8\geq a_9+a_{10},\quad a_5\geq a_8,\cr &\qquad\qquad
\qquad a_1+a_2\geq a_5+a_6,\quad a_{10}\geq a_8,\quad a_6\geq a_2
.\cr &\cr \endalign $$$$\align &(2).\;\;\sum_{0\leq u\leq a_2+a_3}\!\!(-1)^u
\left[\smallmatrix a_5-a_1-1+u\cr u
\endsmallmatrix
\right]e_2^{(a_3)}e_3^{(a_2+a_3-u)}e_4^{(a_1+a_2+a_3)}e_1^{(a_4+a_7)}\cr
&\qquad\qquad\qquad\times
e_2^{(a_4+a_6+a_7)}e_3^{(a_4+a_5+a_6+a_7+u)}e_2^{(a_8)}e_1^{(a_9+a_{10})}e_4^{(a_4)}
e_2^{(a_9)}\cr &\qquad\qquad\text{\rm if}\qquad a_6+a_7\geq
a_2+a_3,\quad a_6+a_8\geq a_9+a_{10},\cr &\qquad\qquad\qquad\qquad
a_5\geq a_1\geq a_8,\quad a_2\geq a_6,\quad a_{10}\geq a_8.\cr
&\cr {\bold {19.}}\quad &(1).\sum_{0\leq u \leq
a_6+a_7}\!\!\!\!\!\!\!\!(-1)^u\left[\smallmatrix a_5-a_8-1+u \cr u
\endsmallmatrix
\right]e_2^{(a_3)}e_4^{(a_1)}e_1^{(a_4+a_7)}e_2^{(a_4)}e_3^{(a_2+a_3+a_4+a_5+u)}\cr
&\qquad\qquad\qquad\times
e_2^{(a_6+a_7+a_8)}e_1^{(a_9+a_{10})}e_4^{(a_2+a_3+a_4)}e_3^{(a_6+a_7-u)}e_2^{(a_9)}\cr
&\quad\qquad\qquad\text{\rm if}\quad a_2+a_3\geq a_6+a_7,\quad
a_6+a_8\geq a_9+a_{10},\quad a_5\geq a_8,\cr &
\qquad\quad\qquad\qquad\qquad  a_{10}\geq a_8\geq a_1,\quad
a_7\geq a_3.\cr &(2).\;\;\sum_{0\leq u\leq a_9}\!\!(-1)^u
\left[\smallmatrix a_8-a_{10}-1+u\cr u
\endsmallmatrix
\right]e_2^{(a_3)}e_4^{(a_1)}e_1^{(a_4+a_7)}e_2^{(a_4)}e_3^{(a_2+a_3+a_4+a_5)}\cr
&\qquad\qquad\qquad\times
e_2^{(a_6+a_7+a_8+u)}e_1^{(a_9+a_{10})}e_4^{(a_2+a_3+a_4)}e_3^{(a_6+a_7)}e_2^{(a_9-u)}\cr
&\quad\qquad\qquad\text{\rm if}\quad a_2+a_3+a_5\geq
a_6+a_7+a_8,\quad a_6\geq a_9,\cr &\qquad\quad\qquad\qquad \qquad
a_8\geq  a_{10}\geq a_5\geq a_1,\quad a_7\geq a_3.\cr &\cr {\bold
{20.}}\quad &(1).\sum_{0\leq u\leq a_4+a_7+a_9}\!\!(-1)^u
\left[\smallmatrix a_{10}-a_8-1+u\cr u
\endsmallmatrix
\right]e_2^{(a_3)}e_3^{(a_2+a_3)}e_2^{(a_6)}e_1^{(a_4+a_7+a_9-u)}\cr
&\quad\quad\quad\times
e_2^{(a_4+a_7)}e_3^{(a_4)}e_4^{(a_1+a_2+a_3+a_4)}e_3^{(a_5+a_6+a_7)}e_2^{(a_8+a_9)}
e_1^{(a_{10}+u)}\cr &\qquad\text{\rm if}\quad a_1+a_2+a_3\geq
a_5+a_6+a_7,\quad a_6+a_7\geq a_2+a_3,\quad a_2\geq a_6,\cr
&\qquad\qquad\qquad a_5+a_6\geq a_8+a_9,\quad a_{10}\geq a_8,\quad
a_9\geq a_6.\cr &\cr &(2).\sum_{0\leq u\leq
a_4+a_7}\!\!(-1)^u\left[\smallmatrix a_8+a_9-a_5-a_6-1+u\cr u
\endsmallmatrix
\right]e_2^{(a_3)}e_3^{(a_2+a_3)}e_2^{(a_6)}e_1^{(a_4+a_7+a_9)}\cr
&\qquad\quad\quad\times
e_2^{(a_4+a_7-u)}e_3^{(a_4)}e_4^{(a_1+a_2+a_3+a_4)}e_3^{(a_5+a_6+a_7)}e_2^{(a_8+a_9+u)}
e_1^{(a_{10})}\cr &\qquad\text{\rm if}\quad a_1+a_2+a_3\geq
a_5+a_6+a_7,\quad a_6+a_7\geq a_2+a_3,\quad a_2\geq a_6,\cr
&\qquad\qquad\qquad a_8+a_9\geq a_5+a_6,\quad a_5\geq a_8\geq
a_{10}.\cr &\cr {\bold {21.}}\quad &(1).\!\!\!\sum_{0\leq u\leq
a_4}\!\!\!\!(-1)^u\!\left[\!\smallmatrix
a_1+a_2+a_3-a_5-a_6-a_7-1+u\cr u
\endsmallmatrix
\!\right]\!\!
e_3^{(a_2)}e_2^{(a_3+a_6)}e_3^{(a_3)}e_4^{(a_1+a_2+a_3+u)}\cr
&\quad\qquad\qquad\times
e_3^{(a_5+a_6)}e_2^{(a_8)}e_1^{(a_4+a_7+a_9+a_{10})}e_2^{(a_4+a_7+a_9)}e_3^{(a_4+a_7)}
e_4^{(a_4-u)}\cr
&\qquad\text{\rm if}\quad a_1+a_2+a_3\geq a_5+a_6+a_7,\quad
a_8+a_9\geq a_5+a_6,\quad a_5\geq a_8,\cr &\qquad\qquad\qquad
a_7\geq a_3,\quad a_6\geq a_2,\quad a_{10}\geq a_8.\cr &\cr
\endalign $$$$\align &(2).\sum_{0\leq u\leq a_4+a_7}\!\!(-1)^u\left[\smallmatrix
a_5+a_6-a_8-a_9-1+u\cr u
\endsmallmatrix
\right]e_3^{(a_2)}e_2^{(a_3+a_6)}e_3^{(a_3)}e_4^{(a_1+a_2+a_3)}\cr
&\quad\quad\times
e_3^{(a_5+a_6+u)}e_2^{(a_8)}e_1^{(a_4+a_7+a_9+a_{10})}e_2^{(a_4+a_7+a_9)}e_3^{(a_4+a_7-u)}
e_4^{(a_4)}\cr &\qquad\qquad\qquad\text{\rm if}\quad
a_5+a_6+a_7\geq a_1+a_2+a_3,\quad a_{10}\geq a_8,\cr
&\qquad\qquad\qquad\qquad a_1+a_2\geq a_5+a_6\geq a_8+a_9,\quad
a_9\geq a_6\geq a_2.\cr &\cr  {\bold {22.}}\quad
&(1).\;\;\sum_{0\leq u\leq a_3}\!\!(-1)^u\left[\smallmatrix
a_5+a_6-a_1-a_2-1+u\cr u
\endsmallmatrix
\right]e_3^{(a_2)}e_2^{(a_3+a_6)}e_1^{(a_4+a_7+a_9)}e_3^{(a_3-u)}\cr
&\qquad\qquad\times
e_4^{(a_1+a_2+a_3)}e_2^{(a_4+a_7)}e_3^{(a_4+a_5+a_6+a_7+u)}e_2^{(a_8+a_9)}e_4^{(a_4)}
e_1^{(a_{10})}\cr &\qquad\qquad\qquad\text{\rm if}\quad
a_5+a_6\geq a_1+a_2 \geq a_8+a_9,\quad a_7\geq a_3,\cr
&\qquad\qquad\qquad\qquad a_9\geq a_6,\quad a_8\geq a_{10},\quad
a_1\geq a_5.\cr &(2).\;\;\sum_{0\leq u\leq
a_4+a_7}\!\!(-1)^u\left[\smallmatrix a_6-a_9-1+u\cr u
\endsmallmatrix
\right]e_3^{(a_2)}e_2^{(a_3+a_6+u)}e_1^{(a_4+a_7+a_9)}e_3^{(a_3)}\cr
&\qquad\qquad\times
e_4^{(a_1+a_2+a_3)}e_2^{(a_4+a_7-u)}e_3^{(a_4+a_5+a_6+a_7)}e_2^{(a_8+a_9)}e_4^{(a_4)}
e_1^{(a_{10})}\cr &\qquad\qquad\qquad\text{\rm if}\quad
a_5+a_6+a_7\geq a_1+a_2+a_3,\quad a_1+a_2\geq a_5+a_6,\cr
&\qquad\qquad\qquad \qquad a_5\geq  a_8\geq a_{10},\quad a_6\geq
a_9\geq a_2.\cr &\cr {\bold {23.}}\quad &(1).\sum_{0\leq u\leq
a_2+a_3}\!\!(-1)^u \left[\smallmatrix a_5-a_1-1+u\cr u
\endsmallmatrix
\right]e_2^{(a_3)}e_3^{(a_2+a_3-u)}e_4^{(a_1+a_2+a_3)}e_2^{(a_6)}\cr
&\quad\quad\times
e_1^{(a_4+a_7+a_9)}e_3^{(a_5+a_6+u)}e_2^{(a_4+a_7+a_8+a_9)}e_3^{(a_4+a_7)}e_1^{(a_{10})}
e_4^{(a_4)}\cr &\qquad\qquad\qquad\text{\rm if}\quad a_6+a_7\geq
a_2+a_3,\quad a_8+a_9\geq a_5+a_6,\cr &\qquad\qquad\qquad \qquad
a_5\geq a_1\geq  a_8\geq a_{10},\quad a_2\geq a_6.\cr &\cr
&(2).\;\;\sum_{0\leq u\leq a_6}\!\!(-1)^u\left[\smallmatrix
a_8-a_5-1+u\cr u
\endsmallmatrix
\right]e_2^{(a_3)}e_3^{(a_2+a_3)}e_4^{(a_1+a_2+a_3)}e_2^{(a_6-u)}\cr
&\quad\quad\times
e_1^{(a_4+a_7+a_9)}e_3^{(a_5+a_6)}e_2^{(a_4+a_7+a_8+a_9+u)}e_3^{(a_4+a_7)}e_1^{(a_{10})}
e_4^{(a_4)}\cr &\quad\quad\text{\rm if}\quad a_5+a_6+a_7\geq
a_1+a_2+a_3,\quad a_2+a_5\geq a_6+a_8,\cr &\qquad\qquad\qquad
\qquad a_1\geq a_5\geq a_{10},\quad a_9\geq a_6,\quad a_8\geq a_5
.\cr&\cr  {\bold {24.}}\quad &(1).\!\!\!\sum_{0\leq u\leq
a_3+a_4}\!\!\!\!(-1)^u\!\!\left[\!\smallmatrix
a_1+a_2-a_5-a_6-1+u\cr u
\endsmallmatrix
\!\right]\!e_3^{(a_2)}e_2^{(a_3+a_6)}e_4^{(a_1+a_2+u)}e_3^{(a_3+a_5+a_6)}\cr
&\qquad\qquad\qquad\times
e_2^{(a_8)}e_1^{(a_4+a_7+a_9+a_{10})}e_2^{(a_4+a_7+a_9)}e_3^{(a_4)}e_4^{(a_3+a_4-u)}
e_3^{(a_7)}\cr &\qquad\quad\text{\rm if}\quad a_7+a_8+a_9\geq
a_3+a_5+a_6,\quad a_1+a_2\geq a_5+a_6,\cr &\qquad\qquad\qquad
a_3\geq a_7,\quad a_5\geq a_8,\quad  a_{10}\geq a_8,\quad a_6\geq
a_2.\cr&\cr \endalign $$$$\align &(2).\;\;\sum_{0\leq u\leq
a_2}\!\!(-1)^u \left[\smallmatrix a_5-a_1-1+u\cr u
\endsmallmatrix
\right]e_3^{(a_2-u)}e_2^{(a_3+a_6)}e_4^{(a_1+a_2)}e_3^{(a_3+a_5+a_6+u)}\cr
&\qquad\qquad\qquad\times
e_2^{(a_8)}e_1^{(a_4+a_7+a_9+a_{10})}e_2^{(a_4+a_7+a_9)}e_3^{(a_4)}e_4^{(a_3+a_4)}
e_3^{(a_7)}\cr &\qquad\qquad\qquad\text{\rm if}\quad
a_7+a_8+a_9\geq a_3+a_5+a_6,\quad a_5\geq a_1\geq a_8,\cr
&\qquad\qquad\qquad \qquad\qquad a_3\geq a_7,\quad  a_{10}\geq
a_8,\quad a_6\geq a_2.\cr&\cr {\bold {25.}}\quad
&(1).\!\!\!\sum_{0\leq u\leq a_4+a_7}\!\!\!\!\!\!\!(-1)^u\!
\left[\!\smallmatrix a_9+a_{10}-a_6-a_8-1+u\cr u
\endsmallmatrix
\!\right]\!e_2^{(a_3)}e_3^{(a_2+a_3)}e_1^{(a_4+a_7-u)}e_2^{(a_4+a_6+a_7)}\cr
&\qquad\qquad\times
e_3^{(a_4)}e_4^{(a_1+a_2+a_3+a_4)}e_3^{(a_5+a_6+a_7)}e_2^{(a_8)}e_1^{(a_9+a_{10}+u)}
e_2^{(a_9)}\cr &\qquad\quad\text{\rm if}\quad a_1+a_2+a_3\geq
a_5+a_6+a_7,\quad a_9+a_{10}\geq a_6+a_8,\cr &\qquad\qquad \qquad
a_6+a_7\geq a_2+a_3,\quad a_2\geq  a_6\geq a_9,\quad a_5\geq
a_8.\cr &(2).\;\;\sum_{0\leq u\leq a_3}\!\!(-1)^u
\left[\smallmatrix a_6-a_2-1+u\cr u
\endsmallmatrix
\right]e_2^{(a_3-u)}e_3^{(a_2+a_3)}e_1^{(a_4+a_7)}e_2^{(a_4+a_6+a_7+u)}\cr
&\qquad\qquad\times
e_3^{(a_4)}e_4^{(a_1+a_2+a_3+a_4)}e_3^{(a_5+a_6+a_7)}e_2^{(a_8)}e_1^{(a_9+a_{10})}
e_2^{(a_9)}\cr &\qquad\quad\text{\rm if}\quad a_1+a_2+a_3\geq
a_5+a_6+a_7,\quad a_2+a_8\geq a_9+a_{10},\cr &\qquad\qquad \qquad
a_7\geq a_3,\quad a_6\geq  a_2,\quad a_{10}\geq a_8,\quad a_5\geq
a_8.\cr&\cr  {\bold {26.}}\quad &(1).\!\!\!\sum_{0\leq u\leq
a_9}\!\!(-1)^u\left[\!\smallmatrix a_8-a_{10}-1+u\cr u
\endsmallmatrix
\!\right]\!e_1^{(a_4)}e_2^{(a_3+a_4)}e_3^{(a_2+a_3+a_4)}e_4^{(a_1+a_2+a_3+a_4)}\cr
&\qquad\qquad\times
e_3^{(a_5)}e_1^{(a_7)}e_2^{(a_6+a_7+a_8+u)}e_1^{(a_9+a_{10})}e_3^{(a_6+a_7)}e_2^{(a_9-u)}\cr
&\qquad\qquad\qquad\text{\rm if}\quad a_2+a_5\geq a_6+a_8,\quad
a_8\geq a_{10}\geq a_5,\cr &\qquad\qquad\qquad \qquad  a_1\geq
a_5,\quad a_3\geq a_7,\quad a_6\geq a_9.\cr&\cr &(2).\sum_{0\leq
u\leq a_3+a_4}\!\!(-1)^u\left[\smallmatrix a_6+a_8-a_2-a_5-1+u\cr
u
\endsmallmatrix
\right]\!e_1^{(a_4)}e_2^{(a_3+a_4-u)}e_3^{(a_2+a_3+a_4)}\cr
&\qquad\times
e_4^{(a_1+a_2+a_3+a_4)}e_3^{(a_5)}e_1^{(a_7)}e_2^{(a_6+a_7+a_8+u)}e_1^{(a_9+a_{10})}
e_3^{(a_6+a_7)}e_2^{(a_9)}\cr &\qquad\quad\text{\rm if}\quad
a_6+a_8\geq a_2+a_5\geq a_9+a_{10},\quad a_1\geq a_5,\cr
&\qquad\qquad\qquad a_{10}\geq a_8,\quad a_3\geq a_7,\quad a_2\geq
a_6.\cr&\cr {\bold {27.}}\quad &(1).\;\;\sum_{0\leq u\leq
a_3+a_4+a_6}\!\!(-1)^u\left[\smallmatrix a_5-a_8-1+u\cr u
\endsmallmatrix
\right]e_4^{(a_1)}e_3^{(a_2+a_5+u)}e_2^{(a_3+a_6+a_8)}\cr &\qquad
\times
e_1^{(a_4+a_7+a_9+a_{10})}e_2^{(a_4)}e_4^{(a_2)}e_3^{(a_3+a_4+a_6-u)}e_4^{(a_3+a_4)}
e_2^{(a_7+a_9)}e_3^{(a_7)}\cr
&\qquad\text{\rm if}\quad a_7+a_9+a_{10}\geq
a_3+a_6+a_8,\quad a_3+a_6\geq a_7+a_9,\cr &\quad\qquad\qquad
 a_8+a_9\geq a_5+a_6,\quad a_5\geq a_8\geq a_1,\quad a_6\geq
a_2.\cr&\cr \endalign $$$$\align &(2).\!\!\!\!\!\sum_{0\leq u\leq
a_7}\!\!(-1)^u \left[\smallmatrix a_6-a_9-1+u\cr u
\endsmallmatrix
\right]e_4^{(a_1)}e_3^{(a_2+a_5)}e_2^{(a_3+a_6+a_8)}e_1^{(a_4+a_7+a_9+a_{10})}\cr
&\qquad\qquad\times
e_2^{(a_4)}e_4^{(a_2)}e_3^{(a_3+a_4+a_6+u)}e_4^{(a_3+a_4)}e_2^{(a_7+a_9)}e_3^{(a_7-u)}\cr
&\qquad\qquad\qquad\text{\rm if}\quad a_7+a_9+a_{10}\geq
a_3+a_6+a_8,\quad a_3\geq a_7,\cr &\qquad\qquad\qquad\qquad
\qquad a_6\geq a_9\geq a_2,\quad a_8\geq  a_5\geq a_1.\cr&\cr
 {\bold {28.}}\quad&(1).\sum_{0\leq u \leq
a_4+a_7+a_9}\!\!\!\!\!\!(-1)^u\left[\smallmatrix a_8-a_{10}-1+u\cr
u
\endsmallmatrix
\right]e_3^{(a_2)}e_4^{(a_1+a_2)}e_3^{(a_5)}e_2^{(a_3+a_6+a_8+u)}\cr
&\quad\quad\quad\times
e_3^{(a_3+a_6)}e_1^{(a_4+a_7+a_9+a_{10})}e_2^{(a_4+a_7+a_9-u)}e_3^{(a_4)}e_4^{(a_3+a_4)}
e_3^{(a_7)}\cr &\qquad\qquad\text{\rm if}\quad a_5+a_6\geq
a_1+a_2,\quad a_8\geq a_{10}\geq a_5,\quad a_1\geq a_5,\cr
&\qquad\qquad \qquad a_7+a_9\geq a_3+a_6,\quad a_3\geq  a_7.\cr
&(2).\sum_{0\leq u\leq a_3+a_4}\!\!\!\!\!(-1)^u \left[\smallmatrix
a_1+a_2-a_5-a_6-1+u\cr u
\endsmallmatrix
\right]e_3^{(a_2)}e_4^{(a_1+a_2+u)}e_3^{(a_5)}e_2^{(a_3+a_6+a_8)}\cr
&\qquad\qquad\times
e_3^{(a_3+a_6)}e_1^{(a_4+a_7+a_9+a_{10})}e_2^{(a_4+a_7+a_9)}e_3^{(a_4)}e_4^{(a_3+a_4-u)}
e_3^{(a_7)}\cr &\qquad\qquad\text{\rm if}\quad a_7+a_9\geq
a_3+a_6,\quad a_{10}\geq a_8\geq a_5,\quad a_6\geq a_2,\cr
&\qquad\qquad\qquad \qquad a_1+a_2\geq a_5+a_6,\quad a_3\geq
a_7.\cr &\cr {\bold {29.}}\quad &(1).\;\;\sum_{0\leq u\leq
a_3+a_6}\!\!(-1)^u \left[\smallmatrix a_8-a_5-1+u\cr u
\endsmallmatrix
\right]e_3^{(a_2)}e_2^{(a_3+a_6-u)}e_1^{(a_4+a_7+a_9)}e_3^{(a_3)}\cr
&\qquad\qquad\times
e_4^{(a_1+a_2+a_3)}e_3^{(a_5+a_6)}e_2^{(a_4+a_7+a_8+a_9+u)}e_3^{(a_4+a_7)}e_4^{(a_4)}
e_1^{(a_{10})}\cr &\qquad\quad\text{\rm if}\quad a_5+a_6+a_7\geq
a_1+a_2+a_3,\quad a_1+a_2\geq a_5+a_6,\cr &\qquad\qquad\qquad
\quad a_8\geq a_5\geq a_{10},\quad a_9\geq a_6\geq a_2.\cr&\cr
{\bold {30.}}\quad &(1).\;\;\sum_{0\leq u\leq
a_2+a_3}\!\!(-1)^u\left[\smallmatrix a_5-a_1-1+u\cr u
\endsmallmatrix
\right]e_2^{(a_3)}e_3^{(a_2+a_3-u)}e_4^{(a_1+a_2+a_3)}e_2^{(a_6)}\cr
&\qquad\qquad\times
e_1^{(a_4+a_7+a_9)}e_2^{(a_4+a_7)}e_3^{(a_4+a_5+a_6+a_7+u)}e_2^{(a_8+a_9)}e_1^{(a_{10})}
e_4^{(a_4)}\cr &\qquad\quad\text{\rm if}\quad a_6+a_7\geq
a_2+a_3,\quad a_1+a_6\geq a_8+a_9,\quad a_2\geq a_6,\cr
&\qquad\qquad\qquad\qquad
 a_5\geq a_1,\quad a_9\geq a_6,\quad a_8\geq a_{10}.\cr&\cr
 {\bold {31.}}\quad &(1).\!\!\!\sum_{0\leq
u\leq a_3+a_4+a_6}\!\!\!\!(-1)^u\left[\smallmatrix a_8-a_5-1+u\cr
u
\endsmallmatrix
\right]e_1^{(a_4)}e_3^{(a_2)}e_2^{(a_3+a_4+a_6-u)}e_1^{(a_7+a_9)}\cr
&\qquad\qquad\times
e_3^{(a_3+a_4)}e_4^{(a_1+a_2+a_3+a_4)}e_3^{(a_5+a_6)}e_2^{(a_7+a_8+a_9+u)}e_1^{(a_{10})}
e_3^{(a_7)}\cr &\qquad\qquad\quad\text{\rm if}\qquad a_1+a_2\geq
a_5+a_6,\quad a_3+a_6\geq a_7+a_9,\cr &\qquad\qquad\qquad\quad
\qquad a_9\geq a_6\geq a_2,\quad a_8\geq a_5\geq a_{10}.
\endalign $$
\par {\rm 2.}\; Applying the map $\Gamma \circ\Psi\circ\Phi$ defined
in {\rm [2]} to cases {\bf 1}$-${\bf 31}, we get another $72$
polynomial elements in one variable in the canonical basis $\bold
B$.
\endproclaim

{\head{2. Proof of Theorem 1.3 }\endhead}

{\bf 2.1.} In order to prove Theorem 1.3. We firstly need the
following identity showed in [9]. Assume that $m\geq k\geq
0,\;\delta\in\bN.$ Then
$$\sum_{0\leq i\leq\delta}(-1)^i\left[\matrix k-1+i\\i
\endmatrix\right]\left[\matrix m\\\delta-i\endmatrix\right]v^{i(m-k)}=
\left[\matrix m-k\\\delta\endmatrix\right]v^{-k\delta}.\tag i$$
\par
Also, we need the following identity showed in [10]. Assume
that\newline $m\geq k\geq 0,\;\delta,n\in\bN.$ Then
$$\align
&\sum_{0\leq i\leq\delta}(-1)^i\left[\matrix k-1+i\\i
\endmatrix\right]\left[\matrix m+n\\\delta-i
\endmatrix\right]v^{i(m-k-n)}\tag ii\cr&\cr
& =\sum_{0\leq t\leq\delta,n}\left[\matrix m-k\\\delta-t\endmatrix
\right]\left[\matrix n\\t\endmatrix\right]v^{-k (\delta-t)-n\delta+t(m+n)}.
\endalign $$\par \vskip0.1cm
 {\bf 2.2.} Now we prove {\rm 1} of Theorem 1.3. It is obvious that
all the elements from case ${\bold 1}$ to case ${\bold {31}}$ in
Theorem 1.3 are fixed by the involution $\bar{\cdot}$. So we only
need to check that these elements lie in $\Cal L$. Moreover, just
like what we have done in [2] that we only prove the most
complicated case here, i.e., case ${\bold 1}$. First of all, we
observe the monomial corresponding to $(1)-(3)$ in case ${\bold
1}$. Using the commutative relations in [2, \S 1.3], we have
$$\align & e_{2}^{(a_3)}
e_{3}^{(a_2+a_3)}e_{4}^{(a_1+a_2+a_3)}e_{2}^{(a_6)}e_{3}^{(a_5+a_6)}e_{2}^{(a_8)}
e_{1}^{(a_4+a_7+a_9+a_{10})}\tag iii\cr&\cr &\quad\times
e_{2}^{(a_4+a_7+a_9)}e_{3}^{(a_4+a_7)}e_{4}^{(a_4)}
=\sum_{\omega\,\in\,\Omega}v^{{\Cal A }(\omega)}\times {\Cal
B}(\omega)\times E^{\omega},
\endalign$$
where
$$\align {\Cal A}(\omega): =
&-(a_4+a_7+a_9+a_{10}-i)(a_4+a_7+a_9-i)-(a_4+a_7-j)(i-j)\cr&\cr
&-(a_1+a_2+a_3-r)(k-r)-(a_3-k+a_6-l)(a_5+a_6-l)\cr&\cr
&-(a_3-k+a_6-l+a_8+a_4+a_7+a_9-i-m)(a_4+a_7-j-m)\cr&\cr
&-(a_3-k)(a_2+a_3-k)-(a_2+a_3-k-t)(a_1+a_2+a_3-r-t)\cr&\cr
\endalign $$$$\align
&-(a_2\!+a_3\!-\!k-\!t+\!a_5\!+a_6\!-l+\!a_4\!+a_7\!-j\!-m-\!p)
(a_4\!-n\!-s-\!p)\cr&\cr
&-(a_4-n)(j-n)-(a_4-n-s)(l+m+k-r-s)\cr&\cr &+(a_4+a_7-j-m)l+
(a_5+a_6-l+a_4+a_7-j-m)(k-r)\cr&\cr &+(a_4-n-s-p)(r+t)+pr,
\endalign$$
and
$$\align {\Cal B}(\omega): &=\left[\matrix a_8+a_4+a_7+a_9-i\\ a_8
\endmatrix\right]\left[\matrix a_3-k+a_6-l+a_8+a_4+a_7+a_9-i\\a_3-k+a_6-l
\endmatrix\right]\cr&\cr &\times\left[\matrix a_3-k+a_6\\a_6
\endmatrix\right]\left[\matrix l+m\\l
\endmatrix\right]\left[\matrix l+m+k-r\\k-r
\endmatrix\right]\left[\matrix p+t\\p
\endmatrix\right]\left[\matrix r+s\\r
\endmatrix\right]
\cr&\cr&\times\left[\matrix a_1+a_2+a_3-r-t+a_4-n-s-p\\a_4-n-s-p
\endmatrix\right]\cr&\cr & \times \left[\matrix a_5+a_6-l+a_4+a_7-j-m\\a_4+a_7-j-m
\endmatrix\right]\cr&\cr&\times\left[\matrix a_2+a_3-k-t+a_5+a_6-l+a_4+a_7-j-m\\a_2+a_3-k-t
\endmatrix\right].\endalign $$
Note that the last two factors in ${\Cal B}(\omega)$ can also be
equivalently represented as follows
$$\align &\left[\matrix a_2+a_3-k-t+a_5+a_6-l\\a_2+a_3-k-t
\endmatrix\right]\cr&\cr&\times\left[\matrix a_2+a_3-k-t+a_5+a_6-l+a_4+a_7-j-m\\a_4+a_7-j-m
\endmatrix\right].
\endalign$$
Moreover, we have
$$\align E^{\omega}: &=e_4^{(a_1+a_2+a_3-r-t+a_4-n-s-p)}e_{34}^{(p+t)}
e_{24}^{(r+s)}e_{14}^{(n)}\cr&\cr&\times
e_3^{(a_2+a_3-k-t+a_5+a_6-l+a_4+a_7-j-m-p)}
e_{23}^{(l+m+k-r-s)}e_{13}^{(j-n)}\cr&\cr&\times
e_{2}^{(a_3-k+a_6-l+a_8+a_4+a_7+a_9-i-m)}
e_{12}^{(i-j)}e_{1}^{(a_4+a_7+a_9+a_{10}-i)},\endalign$$
and
$$\Omega=\{\omega = (i,j,k,l,m,n,r,s,t,p)\}\subset {\bN}^{10},$$
where integers $\,i,j,k,l,m,n,r,s,t,p\,$ satisfies the following
inequalities:
$$\align
&\quad 0\leq i\leq a_4+a_7+a_9;\quad  0\leq j\leq a_4+a_7, i;\quad
0\leq k \leq a_3;\cr&\cr &\quad 0\leq m\leq
a_4+a_7-j,a_3-k+a_6-l+a_8+a_4+a_7+a_9-i;\cr&\cr (\flat) &\quad
0\leq l\leq a_3-k+a_6, a_5+a_6;\quad 0\leq n\leq a_4,j;\quad 0\leq
r\leq k;\cr&\cr &\quad 0\leq s\leq a_4-n,l+m+k-r;\quad 0\leq t\leq
a_2+a_3-k;\cr&\cr &\quad 0\leq p \leq
a_4-n-s,a_2+a_3-k-t+a_5+a_6-l+a_4+a_7-j-m.
\endalign
$$
\par
It should be mentioned here that in the above argument we use $10$
different parameters to describe $\Omega$ in order to simplify the
proof of polynomial elements {\bf 1.}$(1), (2), (3)$ in one
variable, however we use $14$ different parameters to describe
$\Omega$ in the proof of the monomial element {\bf 1.}$(1)$ in
[2]. It is only because the commutative order we use here is
somewhat different from that in [2]. \par Set
$$\align
&\qquad x_1=a_4+a_7+a_9-i,\qquad
x_2=a_4+a_7-j-m,\qquad\qquad\qquad\cr (\natural)&\qquad
x_3=a_3-k,\qquad\qquad\qquad\;x_4=a_3-k+a_6-l,\qquad\qquad\qquad\cr
&\qquad x_5=a_4-n-s-p,\qquad \;\;\; x_6=k-r,\qquad\qquad\qquad\cr
&\qquad x_7=a_2+a_3-k-t.
\endalign$$
Then the degree (with respect to $v$) of the coefficient $v^{{\Cal
A}(\omega)}\times {\Cal B}(\omega)$ of $E^{\omega}$ in the sum
expression of Formula (iii) is
$$D_M: =-L_M(x_1,x_2,\cdots,x_7,m,s,p)-Q_{M}(x_1,x_2,\cdots,x_7,m,s,p),$$
where $L_M(x_1,x_2,\cdots,x_7,m,s,p)$ is a linear form in
non-negative integers $x_1,x_2,\cdots,x_7,m,s,p$, and $Q_{M}(x_1,x_2,\cdots,x_7,m,s,p)$
is a unit form in non-negative integers $x_1,$ $x_2,\cdots,x_7,m,s,p$. Moreover, we have
$$\align
&L_M(x_1,x_2,\cdots,x_7,m,s,p): \cr&\cr
&=(a_{10}-a_8)x_1+(a_8+a_9-a_5-a_6)x_2+(a_2-a_6)x_3+(a_5-a_8)x_4\cr&\cr
&+(a_5+a_6+a_7-a_1-a_2-a_3)x_5+(a_1+a_2-a_5-a_6)x_6\cr&\cr
&+(a_1-a_5)x_7+(a_9-a_6)m+(a_7-a_3)s+(a_6+a_7-a_2-a_3)p,
\endalign$$
and
$$\align
& Q_{M}(x_1,x_2,\cdots,x_7,m,s,p): \cr&\cr
&=x_1^2+x_2^2+x_3^2+x_4^2+x_5^2+x_6^2+ x_7^2+m^2+ s^2+p^2\cr&\cr
&+x_2x_4+x_2m+x_3x_6 +x_3x_7+x_3s+x_3p+x_4m+x_5x_6+x_5x_7\cr&\cr
&+x_5s+x_5p+x_6x_7+x_6s+2x_6p+x_7p+sp-x_1x_4-x_1m-x_2x_6\cr&\cr
&-x_2x_7-x_2s-x_2p-x_3x_4-x_3m-x_4x_7-x_4p-x_6m-ms.
\endalign$$
\par
By the BDDP-algorithm (see [1] or [2] \S 3.3) and Lemma 3.4 in
[2], we know that the unit form $Q_{M}(x_1,x_2,\cdots,x_7,m,s,p)$
is weakly positive, i.e.,
$$Q_{M}(x_1,x_2,\cdots,x_7,m,s,p)\geq 0\quad \text {for any}\quad (x_1,x_2,\cdots,x_7,m,s,p)
\in {\bN}^{10}.$$ Therefore, we have
$$\align
&\sum_{0\leq u\leq a_4+a_7}(-1)^u\left[\smallmatrix a_5+a_6-a_8-a_9-1+u\\
u\endsmallmatrix\right]e_2^{(a_3)}e_3^{(a_2+a_3)}e_4^{(a_1+a_2+a_3)}e_2^{(a_6)}\cr&\cr
&\qquad\;\;\times
e_3^{(a_5+a_6+u)}e_2^{(a_8)}e_1^{(a_4+a_7+a_9+a_{10})}e_2^{(a_4+a_7+a_9)}
e_3^{(a_4+a_7-u)}e_4^{(a_4)}\cr &\cr
&=\sum_{\omega \in
\Omega_1}\left(\sum_{0\leq u\leq a_4+a_7-j-m} (-1)^u\left[\matrix
a_5+a_6-a_8-a_9-1+u
\\u\endmatrix\right]\right.\cr &\cr
&\qquad \times \left.\left[\matrix a_4+a_5+a_6+a_7-j-m-l\\
a_4+a_7-j-m-u\endmatrix\right]\times
v^{(a_4+a_7+a_8+a_9-j-m-l)u}\right)\tag iv\cr &\cr &\qquad \times
v^{{\Cal A}(\omega)}\times {\Cal B}_1(\omega)\times
E^{\omega}\cr\cr &=\sum_{\omega\in\Omega_1}v^{{\Cal
A}(\omega)-(a_5+a_6-a_8-a_9)(a_4+a_7-j-m)}\times {\Cal
B}_1(\omega)\cr&\cr&\qquad
\times\left[\matrix a_4+a_7+a_8+a_9-j-m-l\\
a_4+a_7-j-m\endmatrix\right]\times E^{\omega},\endalign$$ where
the last equality comes from Formula 2.1 (i), $\Omega_1$ is
obtained from $\Omega$ by replacing the defining inequality
$``0\leq l\leq a_3-k+a_6,a_5+a_6"$ in $(\flat)$ by $``0\leq l\leq
a_3-k+a_6,a_5+a_6+u"$, and ${\Cal B}_1(\omega)$ is obtained from
$\Cal B(\omega)$ by deleting the factor $\left[\matrix
a_5+a_6-l+a_4+a_7-j-m\\a_4+a_7-j-m\endmatrix \right].$ Note that
$a_5+a_6$ is always less than or equal to $a_5+a_6+u$, and
relations $(\natural)$ are independent of $u$, even if $u$ occurs
in the upper boundary of the defining inequality of $\,l$, we can
still conclude that the expressions behind the second equal sign
in (iv) are independent of $u$.
\par
Using relations $(\natural)$, we can get the degree $D_{P_1}$
(with respect to $v$) of the coefficient of $E^{\omega}$ in the
last sum expression of (iv), i.e.,
$$\align
D_{P_1}&=D_M-(a_5+x_4-x_3)x_2+(a_8+a_9-a_6+x_4-x_3)x_2\cr&\cr
&\quad-(a_5+a_6-a_8-a_9)x_2\cr&\cr
&=-L_{P_1}(x_1,x_2,\cdots,x_7,m,s,p)-Q_{P_1}(x_1,x_2,\cdots,x_7,m,s,p),
\endalign$$
where
$$\align
&\;\; L_{P_1}(x_1,x_2,\cdots,x_7,m,s,p)\cr&\cr
=&(a_{10}-a_8)x_1+(a_5+a_6-a_8-a_9)x_2+(a_2-a_6)x_3+(a_5-a_8)x_4\cr&\cr
+&(a_5+a_6+a_7-a_1-a_2-a_3)x_5+(a_1+a_2-a_5-a_6)x_6\cr&\cr
+&(a_1-a_5)x_7+(a_9-a_6)m+(a_7-a_3)s+(a_6+a_7-a_2-a_3)p,
\endalign$$
and
$$Q_{P_1}(x_1,x_2,\cdots,x_7,m,s,p)=Q_{M}(x_1,x_2,\cdots,x_7,m,s,p).$$
When
$$\align
&a_5+a_6+a_7\geq a_1+a_2+a_3,\quad a_5+a_6\geq a_8+a_9,\cr&\cr
&\qquad a_1\geq a_5,\quad a_{10}\geq a_8,\quad a_2\geq a_6,\quad
a_9\geq a_6,
\endalign$$
we have $D_{P_1}\leq 0$. Moreover, $D_{P_1}=0\Leftrightarrow
x_1=\cdots =x_7=m=s=p=0,$ and $E^{\omega}=E^A$. Therefore, we have
$$\align
\quad &\sum_{0\leq u\leq a_4+a_7}(-1)^u\left[\smallmatrix a_5+a_6-a_8-a_9-1+u\\
u\endsmallmatrix\right]
e_2^{(a_3)}e_3^{(a_2+a_3)}e_4^{(a_1+a_2+a_3)}e_2^{(a_6)}\cr
&\qquad\;\;\times
e_3^{(a_5+a_6+u)}e_2^{(a_8)}e_1^{(a_4+a_7+a_9+a_{10})}e_2^{(a_4+a_7+a_9)}
e_3^{(a_4+a_7-u)}e_4^{(a_4)}\cr&\cr &\equiv
e_{4}^{(a_1)}e_{34}^{(a_2)}e_{24}^{(a_3)}e_{14}^{(a_4)}e_{3}^{(a_5)}e_{23}^{(a_6)}
e_{13}^{(a_7)}e_{2}^{(a_8)}e_{12}^{(a_9)}e_{1}^{(a_{10})}\qquad
(\mod v^{-1}{\Cal L})\cr &\cr & \quad \text{\rm if}\quad
a_5+a_6+a_7\geq a_1+a_2+a_3,\quad a_5+a_6\geq a_8+a_9,\cr &\qquad
\qquad a_1\geq a_5,\quad a_{10}\geq a_8,\quad a_2\geq a_6,\quad
a_9\geq a_6.
\endalign$$
So we have proved $(1)$ of case $\bold 1$. Let us consider
$$\align
&\sum_{0\leq u\leq a_4}(-1)^u\left[\smallmatrix a_1+a_2+a_3-a_5-a_6-a_7-1+u\\
u\endsmallmatrix\right]
e_2^{(a_3)}e_3^{(a_2+a_3)}e_4^{(a_1+a_2+a_3+u)}\cr
&\qquad\;\;\times
e_2^{(a_6)}e_3^{(a_5+a_6)}e_2^{(a_8)}e_1^{(a_4+a_7+a_9+a_{10})}
e_2^{(a_4+a_7+a_9)}e_3^{(a_4+a_7)}e_4^{(a_4-u)}\cr&\cr &
=\sum_{\omega\in\Omega}\left(\sum_{0\leq u\leq
a_4-n-s-p}(-1)^u\left[\matrix a_1+a_2+a_3-
a_5-a_6-a_7-1+u\\u\endmatrix\right]\right.\cr&\cr
&\qquad\qquad\qquad\times \left[\matrix
a_1+a_2+a_3+a_4-r-t-n-s-p\\a_4-n-s-p-u\endmatrix\right]\cr&\cr
&\left.\qquad\qquad\qquad \times
v^{(a_4+a_5+a_6+a_7-r-t-n-s-p)u}\right) \times v^{{\Cal
A}(\omega)}\times {\Cal B}_2(\omega)\times E^{\omega}\tag v\cr&\cr
&=\sum_{\omega\in\Omega}v^{{\Cal
A}(\omega)-(a_1+a_2+a_3-a_5-a_6-a_7)(a_4-n-s-p)}\times {\Cal
B}_2(\omega)\cr&\cr& \qquad\qquad\times\left[\matrix
a_4+a_5+a_6+a_7-r-t-n-s-p\\a_4-n-s-p\endmatrix\right]\times
E^{\omega},
\endalign$$
where the last equality comes from Formula 2.1 (i), ${\Cal
B}_2(\omega)$ is obtained from $\Cal B(\omega)$ by delating the
factor $\left[\matrix
a_1+a_2+a_3-r-t+a_4-n-s-p\\a_4-n-s-p\endmatrix\right].$
\par
Using relations $(\natural)$, we can get the degree $D_{P_2}$
(with respect to $v$) of the coefficient of $E^{\omega}$ in the
last sum expression of (v), i.e.,
$$\align
D_{P_2}&=D_M-(a_1+x_6+x_7)x_5+(a_5+a_6+a_7-a_2-a_3+x_6+x_7)x_5\cr&\cr
&\quad-(a_1+a_2+a_3-a_5-a_6-a_7)x_5\cr&\cr
&=-L_{P_2}(x_1,x_2,\cdots,x_7,m,s,p)-Q_{P_2}(x_1,x_2,\cdots,x_7,m,s,p),
\endalign$$
where
$$\align
&\;\; L_{P_2}(x_1,x_2,\cdots,x_7,m,s,p)\cr&\cr
=&(a_{10}-a_8)x_1+(a_8+a_9-a_5-a_6)x_2+(a_2-a_6)x_3+(a_5-a_8)x_4\cr&\cr
+&(a_1+a_2+a_3-a_5-a_6-a_7)x_5+(a_1+a_2-a_5-a_6)x_6\cr&\cr
+&(a_1-a_5)x_7+(a_9-a_6)m+(a_7-a_3)s+(a_6+a_7-a_2-a_3)p,
\endalign$$
and
$$Q_{P_2}(x_1,x_2,\cdots,x_7,m,s,p)=Q_{M}(x_1,x_2,\cdots,x_7,m,s,p).$$
When
$$\align
&a_1+a_2+a_3\geq a_5+a_6+a_7,\quad a_8+a_9\geq a_5+a_6,\quad
a_2\geq a_6,\cr &\qquad a_6+a_7\geq a_2+a_3,\quad a_{10}\geq
a_8,\quad a_5\geq a_8,
\endalign$$
we have $D_{P_2}\leq 0$. Moreover, $D_{P_2}=0\Leftrightarrow
x_1=\cdots =x_7=m=s=p=0,$ and $E^{\omega}=E^A$. Therefore, we have
$$\align
\quad &\sum_{0\leq u\leq a_4}(-1)^u\left[\smallmatrix a_1+a_2+a_3-a_5-a_6-a_7-1+u\\
u\endsmallmatrix\right]
e_2^{(a_3)}e_3^{(a_2+a_3)}e_4^{(a_1+a_2+a_3+u)}e_2^{(a_6)}\cr
&\qquad\;\;\times
e_3^{(a_5+a_6)}e_2^{(a_8)}e_1^{(a_4+a_7+a_9+a_{10})}e_2^{(a_4+a_7+a_9)}
e_3^{(a_4+a_7-u)}e_4^{(a_4)}\cr&\cr &\equiv
e_{4}^{(a_1)}e_{34}^{(a_2)}e_{24}^{(a_3)}e_{14}^{(a_4)}e_{3}^{(a_5)}e_{23}^{(a_6)}
e_{13}^{(a_7)}e_{2}^{(a_8)}e_{12}^{(a_9)}e_{1}^{(a_{10})}\qquad
(\mod v^{-1}{\Cal L})\cr &\cr & \quad \text{\rm if}\quad
a_1+a_2+a_3\geq a_5+a_6+a_7,\quad a_8+a_9\geq a_5+a_6,\quad
a_2\geq a_6,\cr &\qquad \qquad \quad a_6+a_7\geq a_2+a_3,\quad
a_{10}\geq a_8,\quad a_5\geq a_8.
 \endalign$$ This gives (2) of case $\bold 1$. Finally, we
 consider
$$\align
&\sum_{0\leq u\leq a_2+a_3}(-1)^u\left[\smallmatrix
a_5-a_1-1+u\\u\endsmallmatrix \right]
e_2^{(a_3)}e_3^{(a_2+a_3-u)}e_4^{(a_1+a_2+a_3)}e_2^{(a_6)}\cr
&\qquad\;\times
e_3^{(a_5+a_6+u)}e_2^{(a_8)}e_1^{(a_4+a_7+a_9+a_{10})}e_2^{(a_4+a_7+a_9)}
e_3^{(a_4+a_7)}e_4^{(a_4)}\cr&\cr &
=\sum_{\omega\in\Omega_1}\left(\sum_{0\leq u\leq a_2+a_3-k-t}
(-1)^u\left[\matrix a_5-a_1-1+u\\ u\endmatrix\right]\right.\cr&\cr
&\qquad\qquad\qquad\times\left[\matrix a_2+a_3+a_5-r-t+a_6-l-k+r\\
a_2+a_3-k-t-u
\endmatrix\right]\cr&\cr
&\left.\qquad\qquad\qquad\times
v^{(a_1+a_2+a_3-a_6+l+k-2r-t)u}\right)\times v^{{\Cal A}(\omega)}
\times {\Cal B}_3(\omega)\times E^{\omega}\tag vi \cr&\cr
&=\sum_{\omega\in\Omega_2}v^{{\Cal
A}(\omega)-(a_5-a_1)(a_2+a_3-k-t-w)-(a_6-l-k+r)(a_2+a_3-k-t)}
\cr&\cr& \qquad\times v^{(a_2+a_3-k-t+a_5+a_6-l)w}\times {\Cal
B}_3(\omega)\times\left[\matrix
a_1+a_2+a_3-r-t\\a_2+a_3-k-t-w\endmatrix\right]\cr&\cr
&\qquad\times\left[\matrix a_6-l-k+r\\w\endmatrix\right]\times
E^{\omega},
\endalign$$
where the last equality comes from the Formula 2.1 (ii),
$\Omega_2$ is obtained from $\Omega_1$ by adding $``0\leq w\leq
a_2+a_3-k-t, a_6-l-k+r"$ to its defining inequalities only, and
${\Cal B}_3(\omega)$ is obtained from $\Cal B(\omega)$ by delating
the factor $\left[\matrix
a_2+a_3-k-t+a_5+a_6-l\\a_2+a_3-k-t\endmatrix\right].$ Note that
the expressions behind the second equal sign in (vi) are also
independent of $u$, just like (1) of case {\bf 1}.
\par
Using relations $(\natural)$, we can get the degree $D_{P_3}$
(with respect to $v$) of the coefficient of $E^{\omega}$ in the
last sum expression of (vi), i.e.,
$$\align
D_{P_3}&=D_M-(a_5+x_4-x_3)x_7-(a_5-a_1)(x_7-w)-(x_4-x_3-x_6)x_7\cr&\cr
&\quad+(a_5+x_7+x_4-x_3)w+(a_1+x_6+w)(x_7-w)\cr&\cr&\quad
+(x_4-x_3-x_6-w)w\cr\endalign $$$$\align
&=-L_{P_3}(x_1,x_2,\cdots,x_7,m,s,p,w)-Q_{P_3}(x_1,x_2,\cdots,x_7,m,s,p,w),
\endalign$$
where
$$\align
&L_{P_3}(x_1,x_2,\cdots,x_7,m,s,p,w)\cr&\cr =&
(a_{10}-a_8)x_1+(a_8+a_9-a_5-a_6)x_2+(a_2+a_5-a_1-a_6)x_3\cr&\cr
&+(a_1-a_8)x_4+(a_5+a_6+a_7-a_1-a_2-a_3)x_5+(a_2-a_6)x_6\cr&\cr
&+(a_5-a_1)(x_7-w)+(a_9-a_6)m+(a_7-a_3)s\cr&\cr
&+(a_6+a_7-a_2-a_3)p+(a_5-a_1)(x_4-x_3-x_6-w),
\endalign$$
and
$$\align
&Q_{P_3}(x_1,x_2,\cdots,x_7,m,s,p,w)\cr&\cr =&
Q_{M}(x_1,x_2,\cdots,x_7,m,s,p)+2(x_4-x_3-x_6-w)(x_7-w).
\endalign$$
By the definition of $\Omega_2$, we have $(x_4-x_3-x_6-w)\geq 0,$
and $(x_7-w)\geq 0$. When
$$\align
&a_8+a_9\geq a_5+a_6,\quad a_6+a_7\geq a_2+a_3,\quad a_2\geq a_6,\cr
&\qquad a_5\geq a_1\geq a_8,\quad a_{10}\geq a_8,
\endalign$$
we have $D_{P_3}\leq 0$. Moreover, $D_{P_3}=0\Leftrightarrow
x_1=\cdots =x_7=m=s=p=w=0,$ and $E^{\omega}=E^A.$ Then we have
$$\align
\quad &\sum_{0\leq u\leq a_2+a_3}(-1)^u\left[\smallmatrix
a_5-a_1-1+u\\u\endsmallmatrix\right]
e_2^{(a_3)}e_3^{(a_2+a_3-u)}e_4^{(a_1+a_2+a_3)}e_2^{(a_6)}\cr
&\qquad\;\;\times
e_3^{(a_5+a_6+u)}e_2^{(a_8)}e_1^{(a_4+a_7+a_9+a_{10})}e_2^{(a_4+a_7+a_9)}
e_3^{(a_4+a_7)}e_4^{(a_4)}\cr&\cr &\equiv
e_{4}^{(a_1)}e_{34}^{(a_2)}e_{24}^{(a_3)}e_{14}^{(a_4)}e_{3}^{(a_5)}e_{23}^{(a_6)}
e_{13}^{(a_7)}e_{2}^{(a_8)}e_{12}^{(a_9)}e_{1}^{(a_{10})}\qquad
(\mod v^{-1}{\Cal L}) \endalign$$ $$\align &\quad \text{\rm
if}\quad a_8+a_9\geq a_5+a_6,\quad a_6+a_7\geq a_2+a_3,\quad
a_2\geq a_6,\cr &\qquad \qquad \qquad a_5\geq a_1\geq a_8,\quad
a_{10}\geq a_8,
\endalign$$ and $(3)$ of case $\bold 1$ is proved.

\bigskip\par
\head {3. Polynomial elements in several variables}\endhead
Finally, we conclude our note with the following two remarks:

\proclaim{Remark 1} {\rm The proof of all the other cases in
Theorem 1.3 is quite similar to that of the above case. Also, the
proof of a polynomial element has close relations with that of the
corresponding monomial element.}\endproclaim

\proclaim{Remark 2}{\rm We see that the regions of $62$ monomial
elements and $144$ polynomial elements in one variable in the
canonical basis $\bold B$ do not fill the space ${\Bbb N}^{10}$.
Actually, we need only to fill the space ${\Bbb N}^{9}$ because
the regions what we consider are independent of $a_4$. Recently,
Marsh claims (see [M]) that there should be $672$ regions
corresponding to elements in the canonical basis $\bold B$ of the
quantized enveloping algebra for type $A_4$. So we believe that in
addition to monomial elements and polynomial elements in one
variable we have worked out, there exist many polynomial elements
in two or more variables in the canonical basis $\bold B$. And we
have been keeping on the computations of polynomial elements in
two variables in the canonical basis $\bold B$. We do have found
more than thirty polynomial elements in two independent variables
in the canonical basis $\bold B$. Here, so-called ``independent
variables" means that the summing in the following two polynomials
is independent of the order of $u$ and $w$. For example,
corresponding to monomial element {\bf 4.}(1) in Theorem 3.1 in
[2], we list the following two polynomial elements in two
independent variables $u$ and $w$:\par
$$\align
&\sum \Sb 0\leq u\leq a_7+a_9 \\ 0\leq w \leq a_2+a_3+a_4 \endSb
(-1)^{u+w}\left[\smallmatrix a_{10}-a_8-1+u \cr u
\endsmallmatrix\right]\left[\smallmatrix a_5-a_1-1+w \cr w
\endsmallmatrix\right]e_1^{(a_4)}e_2^{(a_3+a_4)}\cr
&\qquad \qquad \qquad \times e_3^{(a_2+a_3+a_4-w)}e_2^{(a_6)}
e_1^{(a_7+a_9-u)} e_4^{(a_1+a_2+a_3+a_4)}\cr&\cr &\qquad \qquad
\qquad \times
e_2^{(a_7)}e_3^{(a_5+a_6+a_7+w)}e_2^{(a_8+a_9)}e_1^{(a_{10}+u)},\tag
a \cr&\cr &\qquad \text{\rm if}\quad a_3+a_6+a_8\geq
a_7+a_9+a_{10},\quad a_1+a_6\geq a_8+a_9,\cr &\qquad\qquad \;\;
a_{10}\geq a_8,\quad a_9\geq a_6,\quad a_2\geq a_6,\quad a_5\geq
a_1,
\endalign $$
and
$$\align &\sum \Sb 0\leq u\leq a_2+a_3+a_4
\\ 0\leq w\leq a_4 \endSb
(-1)^{u+w}\left[\smallmatrix a_5-a_1-1+u\cr u
\endsmallmatrix
\right]\left[\smallmatrix a_7+a_9-a_3-a_6-1+w\cr w
\endsmallmatrix
\right]e_1^{(a_4-w)}\cr & \qquad \qquad \qquad \times
e_2^{(a_3+a_4)}e_3^{(a_2+a_3+a_4-u)} e_2^{(a_6)}
e_1^{(a_7+a_9+w)}\cr &\cr&\qquad\qquad\qquad\times
e_4^{(a_1+a_2+a_3+a_4)}e_2^{(a_7)}e_3^{(a_5+a_6+a_7+u)}e_2^{(a_8+a_9)}
e_1^{(a_{10})},\tag b\cr&\cr \endalign $$$$\align
&\qquad\quad\text{\rm if}\qquad a_7+a_9\geq a_3+a_6,\quad
a_1+a_6\geq a_8+a_9,\quad a_5\geq a_1,\cr &\qquad\qquad\qquad
\qquad a_8\geq a_{10},\quad a_2\geq a_6,\quad a_3\geq a_7.
\endalign $$}
\endproclaim

{\head{Acknowledgement}\endhead}

This work is supported in part by the National Natural Science
Foundation of China (10271088) and the Natural Science Foundation
of Henan Province (0311010100). The first named author would like
to thank Professor Nanhua Xi and Professor Kaiming Zhao for their
financial support and for their valuable advice and comments during
his visit the Morning Side Centre of Mathematics in the Academy of
Mathematics and System Sciences in Beijing in May--September 2001.
The second named author is also grateful to the Abdus Salam International
Centre for Theoretical Physics for its financial support and hospitability
during his visit. Finally, both authors would like to thank Professor
Robert Marsh for helpful communications, and for his sending us [1]
and other papers.

\end{document}